\documentclass[a4paper,11pt]{article}
 \usepackage[top=3.5cm,bottom=3.5cm,left=2.5cm,right=2.5cm,a4paper]{geometry}
\usepackage{siunitx}
\usepackage{color}
\usepackage[colorlinks]{hyperref}
\hypersetup{
    colorlinks=true, %% false ==> Boxed links; true ==> colored links
    linkcolor={blue!50!black},
    citecolor={green!60!black},
    urlcolor= {red!80!black}
}
\usepackage{amsmath,amsthm,stmaryrd}
\usepackage{amssymb}
\usepackage{upgreek}
\usepackage{mdframed,empheq}
\usepackage{float,moredefs}
\usepackage{graphicx,pgfplots}
\usepackage[]{caption}
\usepackage{enumitem,longtable,setspace}
\usetikzlibrary{fit,matrix,arrows,decorations.pathmorphing,patterns}
\usepgfplotslibrary{groupplots}
\usetikzlibrary{shapes,arrows,positioning,fit}
\usepackage[textsize=small,colorinlistoftodos]{todonotes}
\pgfplotsset{compat=newest}
\usepackage[]{algorithm,algorithmic}
\usepackage{mathbbol}
\usepackage{subfig}
\usepackage[]{algorithm,algorithmic}
\usepackage{comment}

\numberwithin{equation}{section}
 \usepackage{multirow}
\makeatletter
\renewcommand*\env@matrix[1][\arraystretch]{%
  \edef\arraystretch{#1}%
  \hskip -\arraycolsep
  \let\@ifnextchar\new@ifnextchar
  \array{*\c@MaxMatrixCols c}}
\makeatother

\newcommand{\flo}[1]{{\color{black}{#1}}}
\newcommand{\modif}  [1]{{\color{black}{#1}}}

\newcommand{\modifspace}[1]{{\color{black}{#1}}}

\usepackage{scalerel,stackengine}
\stackMath
\newcommand\reallywidehat[1]{%
\savestack{\tmpbox}{\stretchto{%
  \scaleto{%
    \scalerel*[\widthof{\ensuremath{#1}}]{\kern-.6pt\bigwedge\kern-.6pt}%
    {\rule[-\textheight/2]{1ex}{\textheight}}%WIDTH-LIMITED BIG WEDGE
  }{\textheight}% 
}{0.5ex}}%
\stackon[1pt]{#1}{\tmpbox}%
}

\newcommand*\samethanks[1][\value{footnote}]{\footnotemark[#1]}
\newcommand{\RR}{\mathbb{R}}

\newcommand{\bx}{\mathbf{x}}

\newlength{\plotwidth}  \newlength{\modelwidth}
\newlength{\plotheight} \newlength{\modelheight}
\newcommand{\datafile}{}\newcommand{\modelfile}{}

\newcommand{\dataA}{} \newcommand{\legendA}{}
\newcommand{\dataB}{} \newcommand{\legendB}{}
\newcommand{\dataC}{} \newcommand{\legendC}{}
\newcommand{\dataD}{} \newcommand{\legendD}{}
\newcommand{\dataE}{} \newcommand{\legendE}{}

\newcommand{\myxlabel}{}
\newcommand{\myylabel}{}
\newcommand{\myblue}  {blue!90!white}
\newcommand{\mygreen} {green!50!black}
\newcommand{\myred}   {red!40!white}

\newcommand{\dofedge}  {\hat{\mathsf{n}}}  % dof of trace
\newcommand{\dofpe}     {\mathsf{n}_e}    % dof of w
\newcommand{\dofve}     {\mathsf{m}_e}    % dof component of \nabla w
\newcommand{\nfe}     {\mathsf{n}_{\mathrm{face}}^e} % number of edge of element K^e
\newcommand{\dofglobal}{\hat{\mathrm{n}}}  % size of global unknown

\newcommand{\by}{\mathbf{y}}
\newcommand{\Real}{\mathrm{Re}}
\newcommand{\displacement}{\mathbf{u}}
\newcommand{\velocity}{\mathbf{v}}
\newcommand{\tensS}{\boldsymbol{\sigma}}
\newcommand{\tensE}{\boldsymbol{\epsilon}}
\newcommand{\tensC}{\boldsymbol{C}}
\newcommand{\tensCS}{\boldsymbol{S}}
\newcommand{\Gambold}{\boldsymbol{\Gamma}}

\newcommand{\n}{\boldsymbol{\nu}}
\newcommand{\ii}{\mathrm{i}}

\newcommand{\err}   {\mathfrak{e}}
\newcommand{\errall}{\mathfrak{E}}
\newcommand{\cbtitle}{}
\newcommand{\experiment}{}
\newcommand{\legendpos} {}

\newcommand{\qPomega}{} 
\newcommand{\qSomega}{} 

\newcommand{\formulationU}{(\displacement,\tensS)_{\tensCS}}

\newcommand{\lambdaV}{\boldsymbol{\lambda}_{\velocity}}
\newcommand{\lambdaU}{\boldsymbol{\lambda}_{\displacement}}

\newcommand{\lambdaUh}{\boldsymbol{\lambda}_{\displacement h}}

\newcommand{\tauU}{\boldsymbol{\tau}_{\displacement}}

\newcommand{\identity}{\mathbb{Id}}
\newcommand{\cp}  {\mathrm{c}_{\mathrm{P}}}
\newcommand{\cs}  {\mathrm{c}_{\mathrm{S}}}
\newcommand{\cqp} {\mathrm{c}_{\mathrm{qP}}}
\newcommand{\cqsa}{\mathrm{c}_{\mathrm{qS1}}}
\newcommand{\cqsb}{\mathrm{c}_{\mathrm{qS2}}}
\newcommand{\csh} {\mathrm{c}_{\mathrm{sH}}}

\newcommand{\tauSuIda}{\boldsymbol{\tau}_{\displacement \mathrm{I}}}  
\newcommand{\tauSuIdRe}{\boldsymbol{\tau}_{\displacement \mathbb{R}\pm}}  
\newcommand{\tauSuIdIm}{\boldsymbol{\tau}_{\displacement \mathbb{I}\pm}}  
\newcommand{\tauSuIdRep}{\boldsymbol{\tau}_{\displacement \mathbb{R}+}}
\newcommand{\tauSuIdRem}{\boldsymbol{\tau}_{\displacement \mathbb{R}-}}  
\newcommand{\tauSuIdImp}{\boldsymbol{\tau}_{\displacement \mathbb{I}+}}
\newcommand{\tauSuIdImm}{\boldsymbol{\tau}_{\displacement \mathbb{I}-}}  
\newcommand{\tauSuIdb}{\boldsymbol{\tau}_{\displacement \mathrm{I}_s}}
\newcommand{\tauSuKCa}{\boldsymbol{\tau}_{\displacement \Gambold}}     
\newcommand{\tauSuKCb}{\boldsymbol{\tau}_{\displacement \Gambold_{p}}} 
\newcommand{\tauSuKCunit}{\boldsymbol{\tau}_{\displacement \Gambold_{1}}} 
\newcommand{\tauSuMGa}{\boldsymbol{\tau}_{\displacement \mathrm{G}}}  
\newcommand{\tauSuMGunit}{\boldsymbol{\tau}_{\displacement \mathrm{G}_{1}}}

\newcommand{\plottauSuIda}{}
\newcommand{\plottauSuIdb}{}
\newcommand{\plottauSuKCa}{}
\newcommand{\plottauSuKCb}{}
\newcommand{\plottauSuMGa}{}
\newcommand{\plottauSvIda}{}
\newcommand{\plottauSvIdb}{}
\newcommand{\plottauSvKCb}{}
\newcommand{\plottauSvMGa}{}
\newcommand{\plottauSuMGb}{}

\newcommand{\myorderplot}{}
\newcommand{\myfreqplot} {}

\newtheorem{remark}     {Remark}

\usepackage[capitalize,nameinlink]{cleveref}
\crefname{section}   {Section}   {Sections}
\crefname{subsection}{Subsection}{Subsections}
\Crefname{section}   {Section}   {Sections}
\Crefname{subsection}{Subsection}{Subsections}
\Crefname{figure}    {Figure}    {Figures}

\crefformat{equation}{\textup{#2(#1)#3}}
\crefrangeformat{equation}{\textup{#3(#1)#4--#5(#2)#6}}
\crefmultiformat{equation}{\textup{#2(#1)#3}}{ and \textup{#2(#1)#3}}
{, \textup{#2(#1)#3}}{, and \textup{#2(#1)#3}}
\crefrangemultiformat{equation}{\textup{#3(#1)#4--#5(#2)#6}}%
{ and \textup{#3(#1)#4--#5(#2)#6}}{, \textup{#3(#1)#4--#5(#2)#6}}{, and \textup{#3(#1)#4--#5(#2)#6}}
\Crefformat{equation}{#2Equation~\textup{(#1)}#3}
\Crefrangeformat{equation}{Equations~\textup{#3(#1)#4--#5(#2)#6}}
\Crefmultiformat{equation}{Equations~\textup{#2(#1)#3}}{ and \textup{#2(#1)#3}}
{, \textup{#2(#1)#3}}{, and \textup{#2(#1)#3}}
\Crefrangemultiformat{equation}{Equations~\textup{#3(#1)#4--#5(#2)#6}}%
{ and \textup{#3(#1)#4--#5(#2)#6}}{, \textup{#3(#1)#4--#5(#2)#6}}{, and \textup{#3(#1)#4--#5(#2)#6}}

\crefdefaultlabelformat{#2\textup{#1}#3}

\crefname{proposition}{Proposition}{Propositions}
\Crefname{proposition}{Proposition}{Propositions}
\crefname{definition} {Definition} {Definitions}
\Crefname{definition} {Definition} {Definitions}
\crefname{theorem}    {Theorem}    {Theorems}
\Crefname{theorem}    {Theorem}    {Theorems}
\crefname{remark}     {Remark}     {Remarks}
\Crefname{remark}     {Remark}     {Remarks}
\crefname{assumption} {Assumption} {Assumptions}
\Crefname{assumption} {Assumption} {Assumptions}
\title{
       Numerical investigation of stabilization in the Hybridizable Discontinuous Galerkin 
       method for linear anisotropic elastic equation
       }

\author{
Ha Pham\thanks{Project-Team Makutu, Inria, University of Pau and Pays de l'Adour, TotalEnergies, CNRS UMR 5142}
\and
Florian Faucher\samethanks[1]
\and
H\'el\`ene Barucq\samethanks[1]
}

\date{\today}

\begin{document}
\maketitle

% ============================================================================
\begin{abstract}

This work \flo{is concerned with implementing} the hybridizable discontinuous 
Galerkin (HDG) method to solve the linear anisotropic elastic equation 
in the frequency domain. 
First-order formulation with the compliance tensor and Voigt notation
are employed to provide a compact description of the discretized problem 
and flexibility with highly heterogeneous media.
We further focus on the question of optimal choice\flo{s} of stabilization in 
the definition of HDG numerical traces. 
For this purpose, we construct a hybridized Godunov-upwind flux for 
anisotropic % elasticity 
\flo{elastic media} possessing three distinct wavespeeds.
This stabilization removes the need to choose \flo{a} scaling factor, % s,
contrary to \flo{the} identity and Kelvin--Christoffel based stabilizations which
are popular choices in the literature.
We carry out comparisons among these families
for isotropic and anisotropic material, with constant background 
and highly heterogeneous ones, in two and three dimensions.
\flo{These experiments} establish the optimality of the Godunov stabilization which 
can be used as a reference choice for \flo{a} generic material \flo{in which} different
types of waves propagate.

% This work concerns the implementation of the hybridizable discontinuous 
% Galerkin (HDG) method to solve the linear anisotropic elastic equation 
% in the frequency domain. 
% Voigt notation and a first-order formulation with the compliance tensor 
% are employed to provide compact description of the discretized problem 
% and flexibility with highly heterogeneous media.
% We further focus on the question of optimal choice of stabilization in 
% the definition of HDG numerical traces. 
% For this purpose, we construct a hybridized Godunov-upwind flux for 
% anisotropic elasticity possessing three distinct wavespeeds.
% This stabilization removes the need to choose scaling factors and can 
% be used as a reference choice for all materials.
% Its optimality is established by comparing with identity and 
% Kelvin--Christoffel based stabilizations in two and three dimensions, 
% for isotropic and anisotropic material, with constant background and 
% highly heterogeneous ones.  
  
\end{abstract}
% ============================================================================

% \tableofcontents 

% ------        main sections ------------------------------------------
% ============================================================================
%
\section{Introduction}
%
% ============================================================================

In this work, we implement the Hybridizable Discontinuous Galerkin (HDG) method to solve 
the linear anisotropic elastic equation in \flo{the} frequency domain and construct robust stabilization
\flo{in the definition of the} numerical traces \flo{required by the method}. 
The investigation has \flo{in mind} applications which are geared towards 
 wave simulation and full waveform inversion (FWI) in heterogeneous Earth subsurfaces.
To discretize with the HDG method, the elastic equation with its constitutive law is reformulated
as a first-order system with \flo{a} pair of unknowns consisting of displacement $\mathbf{u}$ and stress $\boldsymbol{\sigma}$.
This \flo{system} is then written as a union of local problems defined on each mesh element with prescribed data given by
the HDG global variable $\lambdaU$ representing the trace of $\mathbf{u}$. 
The local problems are linked by transmission conditions and boundary conditions 
defined on \flo{the} interfaces and exterior boundary, respectively.
Together with a definition of \flo{the} numerical trace of the traction 
\flo{denoted by} $\widehat{\boldsymbol{\sigma}}$ 
on the boundary of each element, the transmission conditions combined with the inversion of 
local problems result in a global problem in $\lambdaU$ defined on the skeleton of the 
mesh. A choice of stabilization operator denoted by $\tauU$ is required to define $\widehat{\boldsymbol{\sigma}}$, and
the judicial choice of which is important to ensure the robustness and accuracy of the 
discretization method.

Works implementing HDG method in linear elasticity are e.g., \cite{bonnasse2018hybridizable} 
for time-harmonic elasticity, \cite{sevilla2018superconvergent,giacomini2019discontinuous} 
for static isotropic elasticity, 
\cite{soon2009hybridizable,fu2015analysis,cockburn2013superconvergent,du2020new} 
for static elasticity including anisotropy, and
\cite{terrana2018spectral,nguyen2011high,fernandez2018hybridized} 
for time-domain isotropic elasticity. 
These \flo{works}, including ours, follow the approach of the original HDG method devised by 
Cockburn and co-authors for second-order elliptic equation\flo{s} \cite{cockburn2008superconvergent,cockburn2009unified}.
\flo{The} HDG method retains the attractive features of the Discontinuous Galerkin family, such as 
being locally conservative, flexible with h-p adaptivity strategy  
\modif{and complex geometry \cite{fernandez2018hybridized}}. 
Additionally, it provides optimal convergence with 
\modif{$L^2$-based approximation} with equal degree \flo{of polynomials} for 
the primal variable and its derivatives. Finally, it is highly parallelizable
and amenable to hybrid computing architecture, cf. \cite{fabien2020gpu}.
For more detailed discussions of the method, we refer to 
\cite{arnold2002unified,cockburn2009unified,cockburn2016static,soon2009hybridizable,sevilla2018superconvergent,giacomini2019discontinuous,
cockburn2023hybridizable},
and for its implementation to various equations, \cite{fernandez2018hybridized,giacomini2021hdglab}. 
Other methods related to HDG are Weak Galerkin or Hybrid High Order, cf. \cite{wang2018weak} and \cite{di2015hybrid,cockburn2016bridging}  for elasticity.

In view of application\flo{s} in seismology, in particular with FWI, 
\flo{strong features of the HDG method are its flexibility with} 
% can handle 
material heterogeneity and \modif{unstructured meshes}, 
and \flo{that it can} provide accurate and stable performance at 
tenable computational costs.
\flo{The} Finite Difference method can outperform Finite Element methods on 
media with low-variation \flo{properties}; however, its need for structured grids
creates problems for backgrounds with strong variation and topography (roughness of boundary), 
and \flo{with implementation of} boundary conditions, cf. \cite[Section 4.2]{tournier20223d}.
In our work, the above advantages of \flo{the} HDG method are further exploited by employing a first-order formulation 
in $(\mathbf{u},\boldsymbol{\sigma})$ with \flo{the} constitutive law written in terms of the compliance tensor.
This allows for parameter variation within a mesh cell and, 
as shown in our numerical investigation, leads to \flo{a} better 
representation of complex media and \flo{a more} numerically stable resolution for these media. 
Another feature of our implementation is that the discretized problem is compactly described using Voigt representation,
which maximizes matrix operations and provides \flo{an} efficient book-keeping 
of physical parameters for any generic anisotropic material. 
These features are important not only in code building for the forward 
but also for \flo{the} inverse solver.

In this work, we will focus on the question of choosing a robust 
choice for the stabilization matrix $\tauU$. 
Most of \flo{the} above-cited references which employ \flo{the} HDG method for elasticity
work with \flo{an} identity-based stabilization, i.e., choosing $\tauU=\tau \identity$ 
with $\identity$ the identity matrix and $\tau$ a scalar scaling factor.
This is with the exception of \cite{soon2009hybridizable,fu2015analysis,cockburn2013superconvergent} 
which investigate \flo{stabilization of the form} $\tauU = \tau \boldsymbol{\Gamma}$ 
where $\boldsymbol{\Gamma}$ 
is the Kelvin--Christoffel (KC) matrix for elastostatics,
and \cite{terrana2018spectral} which constructs the hybridized Godunov-upwind flux for time-domain isotropic elasticity.
For static isotropic elasticity, 
numerical investigations on choices \flo{for the} scaling factor $\tau$ 
are carried out for identity-based stabilizations in  \cite{sevilla2018superconvergent}, and 
in \cite{soon2009hybridizable,fu2015analysis} for  KC-based stabilizations. 
We also refer to  \cite{gopalakrishnan2015stabilization} for 
a study of scaling factors \flo{for} identity-based stabilization\flo{s} for \flo{the} 
Helmholtz equation considering complex-valued scaling factors.
On the other hand, there are no such investigations for elastodynamics, with isotropy or anisotropy.

We start our investigation by extending the hybridization technique of upwind fluxes 
for isotropy in \cite{terrana2018spectral} to anisotropy under the assumption 
of three distinct wave speeds. 
The ingredients for its construction in anisotropy, which remain the same as in isotropy, 
consist of \flo{the} exact solutions of the Riemann problem and their Rankine--Huguniot jump conditions. 
They are also employed in 
\cite{wilcox2010high,tie2018unified,tie2020systematic,zhan2018exact} to construct upwind 
flux for DG implementation in \flo{the} time domain, which depend however
on information from both sides of an interface, 
 in contrast with  
the hybridized HDG fluxes in our work and \cite{terrana2018spectral},
which are defined from only one-sided data. 
Hybridization of \flo{the} Godunov flux is also discussed in a series of works 
with different formalisms \cite{bui2015godunov,bui2015rankine,bui2016construction} 
and in which the Rankine--Hugoniot jump conditions are recognized 
in sequel works (\cite{bui2015rankine,bui2016construction}) to 
\flo{provide} a more natural and direct way to devise \flo{the} HDG flux. 
Also employing the notion of Riemann problem\flo{s} to devise \flo{the} 
HDG scheme for compressible fluid flow are, e.g., 
\cite{vila2021hybridisable,nguyen2012hybridizable} and the references therein,
in which hybridization is carried out with Lax--Friedrichs, Roe and Harten--Lax--Van Leer 
approximate Riemann solver\flo{s}. 

In the second part of this work, 
% in comparing all three families of stabilization:
% identity-based, KC-based and hybridized Godunov, 
% cf. \cref{subsection:stabilization}, for various 
% settings and waves, 
% 
\modif{we compare three families of stabilization:
identity-based, KC-based, and hybridized Godunov 
(cf. \cref{subsection:stabilization}) for various 
settings and types of waves}.
We show that the Godunov stabilization as derived is optimal, while the other two families
require correct scaling factors which might not be universal in order to match the 
performance \flo{of the Godunov stabilization}.
The comparison is also carried out \flo{for} setting \flo{supporting mixed types of wavefields 
in} highly varying backgrounds, to conclude on the versatility of \flo{the} Godunov stabilization  
and the fact that it provides a reference choice for generic anisotropic material
and general simulations.

The organization of this article is as follows. We gather in \cref{Genot::sec}
the notations needed in later sections, in particular for \flo{the} discretization 
and in the construction of the hybridized flux.
In \cref{Formulation::sec}, a first-order formulation for the elastic 
equation, its HDG formulation, and \flo{corresponding} discretized problem are given. 
The construction for \flo{elastic} anisotropy is carried out in \cref{godunov_flux::sec}, and 
numerical experiments are displayed in \cref{section:numerics}. 
% \flo{Regarding} Voigt notations, we only introduce the quantities and identities needed,
% and refer to our work \cite{Pham2023hdgRR} for more details on their derivation. 

% ---------------------------------------------------------------
%
\section{Notations}\label{Genot::sec}
%
% ---------------------------------------------------------------

We introduce here the quantities that appear in the rest of the paper. 
In particular, the quantities with Voigt notation (\cref{VI::subsec,subsection:KVmatrix})
are employed to describe the discretized problem \flo{in} \cref{discretization::sec}, 
and to construct the hybridized Godunov-upwind flux in \cref{discretization::sec}.
% \flo{For further discussions, we refer to our work \cite{Pham2023hdgRR}}.

\subsection{General notations}

We denote by $\mathcal{S}_2$ the set of $3\times 3$ symmetric 
matrices and $\mathcal{S}_4$ the set of fourth order tensors 
bearing the symmetry of the elastic tensor:
\begin{equation}
\begin{aligned}
\mathcal{S}_2 & := \{ \boldsymbol{w} = (w_{ij})_{i,j=x,y,z}\in \mathbb{C}^{3\times 3} \,\big| \, w_{ij} = w_{ji}\} \,,\\[0.3em]
 \mathcal{S}_4 &:= \left\{ \mathbf{T}=(t_{ijkl})_{i,j,k,l =x,y,z}  \in 
 \mathbb{C}^{3 \times 3 \times 3 \times 3} \, \big|\,   t_{ijkl}
  = t_{jikl} = t_{ijlk} = t_{klij}  \right\} \,.
\end{aligned} 
\end{equation}

%Tensor-valued functions whose components taking value in a function space $\mathcal{V}$ are written as,
%\begin{equation} \label{eq:definition-Vfunc}
%\text{vector-valued}\hspace*{0.3cm}  \mathcal{V}^n \,, \hspace*{0.5cm} 
%\text{matrix-valued}\hspace*{0.3cm} \mathcal{V}^{n \times n } \,, \hspace*{0.5cm} \substack{\text{symmetric}\\\text{ matrix-valued}} \hspace*{0.3cm}\mathcal{V}^{n\times n}_\mathrm{sym}\,.
%\end{equation}

We list here the operations employed on the following objects:
\flo{for} vectors $ \mathbf{v}=(v_i)_{i=1}^3,  \mathbf{w}=(w_i)_{i=1}^3$, 
matrices $ \boldsymbol{\chi}=(\chi_{kl})_{k,l=1}^3,  \tilde{\boldsymbol{\chi}}=(\tilde{\chi}_{kl})_{k,l=1}^3$, and 
fourth-order tensor $ \mathbf{T} = (t_{ijkl})_{i,j,k,l=1}^3$,

\begin{subequations}
\begin{minipage}{0.45\textwidth}
\begin{equation}  \mathbf{v}\cdot \mathbf{w} = \sum_{j=1}^3 v_j w_j\,, \label{vdotprod} \end{equation} 
\begin{equation}  \boldsymbol{\chi} : \tilde{\boldsymbol{\chi}}   =  \sum_{i,j=1}^3 \chi_{ij} \, \tilde{\chi}_{ij},\label{contrac::def}\end{equation}
\begin{equation}
(\velocity\otimes \mathbf{w})_{ij}  = v_i \, w_j\,,\label{tensorprod}
\end{equation}
\begin{equation}
\mathbf{v}\cdot \boldsymbol{\chi} \cdot \mathbf{w} = \sum_{i,j=1}^3 v_i \,\chi_{ij} \, w_j, \label{bilinform}
\end{equation}
\end{minipage}\hspace*{0.1cm}
\begin{minipage}{0.5\textwidth}
\begin{equation}   \boldsymbol{\chi} \mathbf{v} = \sum_{j=1}^3 \chi_{ij}\, v_j ,\end{equation}
\begin{equation}
(\mathbf{T} \boldsymbol{\chi})_{ij}  = \sum_{k,l=1}^3 t_{ijkl} \, \chi_{kl}, \label{prod_symten_mat}\end{equation}
\begin{equation} \label{eq:def-odot}
 \mathbf{v}\odot\mathbf{w} = \dfrac{\mathbf{v}\otimes \mathbf{w} + \mathbf{w}\otimes \mathbf{v}}{2},
\end{equation}
\begin{equation}
\boldsymbol{\chi} : \mathbf{S} : \tilde{\boldsymbol{\chi}} = \sum_{i,j,k,l=1}^3 \chi_{ij} S_{ijkl} \tilde{\chi}_{kl}.
\end{equation}
\end{minipage}
\end{subequations}\newline
%We will employ the bilinear forms with vectors and matrices, 
%\begin{equation}
%\mathbf{v}\cdot \boldsymbol{\chi} \cdot \mathbf{w} = \sum_{i,j=1}^3 v_i \,\chi_{ij} \, w_j\,,
%\hspace*{1.50cm} \boldsymbol{\chi} : \mathbf{S} : \tilde{\boldsymbol{\chi}} = \sum_{i,j,k,l=1}^3 \chi_{ij} S_{ijkl} \tilde{\chi}_{kl}.
%\end{equation}
We  also employ the following differential operators: 
\begin{equation} \label{eq:symmetric-gradient}
(\nabla \mathbf{w})_{ij} =   \partial_j w_i , \hspace*{1.2cm}  \nabla^\mathfrak{s} \mathbf{w} = \dfrac{\nabla \mathbf{w} + (\nabla \mathbf{w})^t }{2}, \hspace*{1.2cm} (\nabla\cdot \boldsymbol{\chi})_i =  \sum_{j=1}^3 \partial_j \chi_{ij}  \,.
\end{equation}

\subsection{Voigt identification}\label{VI::subsec}

We define the following mapping $\mathit{i}_\mathrm{sm}$ 
which identifies the indices of a symmetric matrix \flo{with}
integer values, 
\begin{equation}\label{ismmap::def}
\mathit{i}_\mathrm{sm} :\hspace*{0.2cm} \begin{array}{ l l l }
  \hspace*{1.1cm} 
  \text{ `xx'} \mapsto 1 , \hspace*{0.2cm}&  \hspace*{0.8cm}\text{ `yy'} \mapsto 2, \hspace*{0.2cm} & \hspace*{0.8cm}\text{ `zz'} \mapsto 3,\\
\hspace*{0.2cm} \text{ \modif{`zy'}},\text{ `yz'}  \mapsto 4, \hspace*{0.2cm}&\text{ `xz'},\text{ `zx'}  \mapsto 5,
\hspace*{0.2cm}&\text{ \modif{`xy'}},\text{ \modif{`yx'}}  \mapsto 6.
\end{array}
\end{equation}
% In this way, the  upper-diagonal indices of a symmetric matrix
% are identified with the lower-diagonal ones. 
\flo{We} denote by $ \mathcal{I}_{\mathrm{sm}}$ the \emph{ordered} 
set consisting of the diagonal and lower-diagonal (equivalently 
upper-diagonal) indices,
% Its members are identified with integer value by $\mathit{i}_\mathrm{sm}$,
\begin{equation}\label{Ism::def}
  \begin{aligned}\mathcal{I}_{\mathrm{sm}} 
:= \left( \text{ `xx'\, , \,\,  `yy'\,, \,\, `zz'\,, \,\,  `yz'\,, \,\, `xz'\,, \,\,  `xy'} \right)\\
=  \left( \text{ `xx'\, , \,\,  `yy'\,, \,\, `zz'\,, \,\,  `zy'\,, \,\, `zx'\,, \,\,  `yx'} \right)\end{aligned} \hspace*{0.7cm} 
\overset{\substack{\text{identified}\\[0.5em]}}{\underset{\substack{\\[0.1em]\text{ with}}}{\longleftrightarrow}}  \hspace*{0.7cm} (1, \ldots, 6)\,.
\end{equation}

% \paragraph{Voigt identifications} 
In Voigt notation, a symmetric $3\times 3$ matrix $\boldsymbol{w}$ is identified 
with a vector of length $6$, denoted by $\overrightarrow{\boldsymbol{w}}$, by employing identification \cref{ismmap::def},
\begin{equation}\label{Voigt-sm}
\begin{aligned}
\mathcal{S}_2\, \ni\, \boldsymbol{w} =  \begin{psmallmatrix} w_{xx}  &  w_{xy}  &  w_{ xz }  \\[0.2em]
w_{yx}  &  w_{ yy }  & w_{ yz} \\[0.2em]
w_{zx }  &  w_{zy }  & w_{ zz } 
\end{psmallmatrix}
\hspace*{0.2cm} \underset{\text{Voigt identification}}{\longleftrightarrow} \hspace*{0.2cm} \overrightarrow{\boldsymbol{w}} &= % \small
\begin{pmatrix} w_{xx} \, , \,  w_{yy} \, , \,   w_{ zz}\, , \, w_{yz}\, , \, 
w_{xz}\, , \,  w_{xy} \end{pmatrix}^t \\
&= (w_{\mathit{i}_\mathrm{sm}(\mathfrak{J})})_{\mathfrak{J}\in \mathcal{I}_{\mathrm{sm}}}.
\end{aligned}
\end{equation}
A fourth order symmetric tensor in $S_4$ is identified with a symmetric matrix of size $6\times 6$ as
\begin{equation}\label{Voigt-s4tensor}
\mathcal{S}_4 \, \ni \, \mathbf{T} = (\mathsf{t}_{ijkl})_{i,j,k,l=1}^3
\hspace*{0.2cm}  \underset{\text{Voigt identification}}{\longleftrightarrow}  \hspace*{0.2cm}  \overset{=}{\mathbf{T}}=\left( \mathsf{t}_{\mathfrak{J}\mathfrak{J}'} \right)_{\mathfrak{J},\mathfrak{J}' \in \mathcal{I}_{\mathrm{sm}}}
 =\left( \mathsf{t}_{\mathit{i}_\mathrm{sm}(\mathfrak{J}) \mathit{i}_\mathrm{sm}(\mathfrak{J}')} \right)_{\mathfrak{J},\mathfrak{J}' \in \mathcal{I}_{\mathrm{sm}}}  \,.
\end{equation}
We will also work with the following identifications, 
  \begin{equation}\label{dCd::def}
\overrightarrow {\boldsymbol{w} }^\dagger  :=   D^{\dagger}\, \overrightarrow {\boldsymbol{w} }\,, \hspace*{0.8cm} \overset{=}{\mathbf{T}}^\dagger  = \overset{=}{\mathbf{T}} \, D^\dagger \,, \hspace*{0.8cm}
 ^\dagger\!\overset{=}{\mathbf{T}}\,\! ^\dagger :=  D^\dagger \overset{=}{\mathbf{T}} \, D^\dagger ,
 \end{equation}
 %with diagonal matrix $D^\dagger$,
\begin{equation}\label{Ddagger}
 \modif{\text{with diagonal matrix}} \hspace*{0.2cm} D^\dagger :=\begin{pmatrix} \identity_{3\times3} & \boldsymbol{0}_{3\times 3} \\
  \boldsymbol{0}_{3\times 3} &  \textcolor{black}{2} \, \identity_{3\times3}  \end{pmatrix} .
\end{equation}

\modif{The definition in \cref{dCd::def} and matrix $D^\dagger$ arise when tensor operations are rewritten with Voigt identification. For example, the contraction \cref{contrac::def} between two symmetric matrices $\boldsymbol{\tau_1}, \boldsymbol{\tau_2} \in S_2$ can be written as,
\begin{equation}\label{contractionV}
\boldsymbol{\tau_1} : \boldsymbol{\tau_2}   = 
\overrightarrow{\boldsymbol{\tau_1}} \cdot \overrightarrow{\boldsymbol{\tau_2}  }^\dagger
 = \overrightarrow{\boldsymbol{\tau_1}}^\dagger \cdot \overrightarrow{\boldsymbol{\tau_2}  }
  = \overrightarrow{\boldsymbol{\tau_1}}\cdot D^\dagger \cdot \overrightarrow{\boldsymbol{\tau_2},  }
\end{equation}
where the operation in the second and third expression are vector dot products \cref{vdotprod}, 
while the fourth expression is the bilinear form defined in \cref{bilinform}. The product \cref{prod_symten_mat} between a symmetric tensor $\tensC \in S_4$ and a matrix $\boldsymbol{\chi}\in S_2$ can be written as, 
\begin{subequations}\label{4tensormat}
\begin{align}
\overrightarrow{\tensC\boldsymbol{\chi}} \,\,&=\,\,   \overset{=}{\tensC} \,  D^\dagger \,  \overrightarrow{\boldsymbol{\chi}} \,\, =\,\, \overset{=}{\tensC}^\dagger \overrightarrow{\boldsymbol{\chi}}
 \,\, = \,\, \overset{=}{\tensC} \, \overrightarrow{\boldsymbol{\chi}}^\dagger \,; \label{4tensormat_V}\\[.2em]
\overrightarrow{\tensC\boldsymbol{\chi}}^{\dagger} \,\, &= \,\,
D^{\dagger}\, \overset{=}{\tensC} \,  D^\dagger \,  \overrightarrow{\boldsymbol{\chi}}\,\,  \overset{\cref{dCd::def}}{=}\, \, ^\dagger\!\overset{=}{\tensC}\,\! ^\dagger\,  \overrightarrow{\boldsymbol{\chi}} 
 \,\, = \,\, \hspace*{0.1cm} D^{\dagger} \, \overset{=}{\tensC} \, \overrightarrow{\boldsymbol{\chi}}^\dagger\,.\label{4tensormat_K}
\end{align}
\end{subequations} 
For other identities with Voigt notation, we refer to \cite[Section 3]{Pham2023hdgRR}.}

%% ------------------------------------------------
%\subsection{Constitutive law in linear elasticity} 
%% ------------------------------------------------
\subsection{\modif{Useful facts in linear elasticity}}
%We note the following useful facts regarding the tensors 
%appearing in the elastic equation, cf. \cref{1stsystem_variants:u}. 

\modif{We note the following facts which will be employed in the discretization of  the elastic equation, cf. \cref{1stsystem_variants:u}}.
The stiffness tensor $\mathbf{C} =  (C_{ijkl})_{i,j,k,l=1}^3$, compliance tensor $
\tensCS =  (S_{ijkl})_{i,j,k,l=1}^3$, the strain matrix $\boldsymbol{\epsilon}= (\epsilon_{IJ})_{I,J=x,y,z}$, 
and stress matrix $\boldsymbol{\sigma}= (\sigma_{IJ})_{I,J=x,y,z}$, are such that,
\begin{equation}
\mathbf{C}, \, \tensCS \in\mathcal{S}_4, \hspace*{1cm} \mathbf{C}=\tensCS^{-1},\hspace*{1cm} 
\boldsymbol{\epsilon}, \, \boldsymbol{\sigma} \in \mathcal{S}_2.
\end{equation}
Their relations are given by the constitutive law of linear elasticity,
\begin{equation}\label{constitutivelaw_gen}
\boldsymbol{\sigma} = \mathbf{C} \boldsymbol{\epsilon}, \hspace*{1cm} \text{or equivalently} \hspace*{1cm}
\boldsymbol{\epsilon} = \mathbf{S} \boldsymbol{\sigma}.
\end{equation}
\noindent \textbf{\modif{Fact} 1.} We write the unit vectors 
in Cartesian basis of $\RR^3$ and $\mathbb{R}^6$, respectively as,
\begin{equation}\label{unit6v::def}
\mathbf{e}_I \in \mathbb{R}^3,  \hspace*{0.1cm}  \text{with } I=x,y,z \,, \hspace*{0.4cm} \text{and} \qquad
\hat{\mathbb{e}}_\mathfrak{J} \in \mathbb{R}^6,  \hspace*{0.1cm}  \text{with }  \mathfrak{J} \in \mathcal{I}_{\mathrm{sm}} \, .
\end{equation}
The latter basis is indexed by $\mathcal{I}_{\mathrm{sm}}$ \cref{Ism::def}. 
We define the following symmetric matrix using \cref{eq:def-odot},
\begin{equation}\label{hsfE_to_smat}
   \overset{=}{\mathbb{e}}_\mathfrak{J} \,  := \, 
   \begin{dcases}  \hspace*{0.2cm} \mathbf{e}_I \odot \mathbf{e}_J\,, \hspace*{0.1cm}  & \hspace*{0.1cm} I=J, \\
   2 \,    \mathbf{e}_I \odot \mathbf{e}_J\,, \hspace*{0.1cm}  & \hspace*{0.1cm} I\neq J. \end{dcases} 
 \end{equation}
A basis for symmetric matrices is given by $\{\overset{=}{\mathbb{e}}_\mathfrak{J}\}$, and in the Voigt identification \cref{Voigt-sm}, by $\{\hat{\mathbb{e}}_\mathfrak{J}\} $. In particular for the stress tensor, $\boldsymbol{\sigma}=(\sigma_{IJ})_{I,J=1}^3\in \mathcal{S}_2$, we have
\begin{equation}\label{decomsigma_v}
\boldsymbol{\sigma} =  \sum_{\mathfrak{J} \in \mathcal{I}_{\mathrm{sm}}}   \sigma_{\mathfrak{J}}  \,  \overset{=}{\mathbb{e}}_{\mathfrak{J}}\,, \qquad
\overrightarrow{\boldsymbol{\sigma}} =  \sum_{\mathfrak{J} \in \mathcal{I}_{\mathrm{sm}}}   \sigma_{\mathfrak{J}}  \,  \hat{\mathbb{e}}_{\mathfrak{J}} \, .
\end{equation}

\noindent \textbf{\modif{Fact} 2.}  
\modif{Employing identities in \cref{4tensormat}}, the constitutive laws \cref{constitutivelaw_gen} can be written with Voigt notations \cref{Voigt-sm} and \cref{Voigt-s4tensor} as,
\begin{equation}
\boldsymbol{\sigma} = \mathbf{C} \boldsymbol{\epsilon} \hspace*{0.1cm}\Leftrightarrow\hspace*{0.1cm}
\overrightarrow{\boldsymbol{\sigma}} =  \overset{=}{\tensC}   D^\dagger \overrightarrow{\boldsymbol{\varepsilon}}\,,
\hspace*{1cm} \boldsymbol{\epsilon} = \mathbf{S} \boldsymbol{\sigma} \hspace*{0.1cm}\Leftrightarrow\hspace*{0.1cm}
 \overrightarrow{\boldsymbol{\epsilon}} =  \overset{=}{\tensCS}  D^\dagger \overrightarrow{\boldsymbol{\sigma}}\,.
\end{equation}
The inverse relation between the stiffness tensor $\tensC$ and compliance tensor $\tensCS$ can be written as,
\begin{equation}\label{rel_voigt_inverse}
 \overset{=}{\tensCS}\, = \, (^\dagger\!\overset{=}{\tensC}\,\! ^\dagger)^{-1}\,, \hspace*{1cm} ^\dagger\!\overset{=}{\tensCS}\,\! ^\dagger \,\,  = \,\, \big( \overset{=}{\tensC} \big)^{-1} \,.
\end{equation}
This is useful for computing the compliance tensor \modif{of a} generic anisotropic elastic material.

% ------------------------------------
\subsection{Kelvin--Christoffel (KC) matrix}
\label{subsection:KVmatrix}
% ------------------------------------

The $3\times 3$ Kevin-Christoffel matrix in direction $\n$ is defined as, 
cf. \cite{carcione2007wave},
\begin{equation}\label{KCmat_def}
\boldsymbol{\Gamma}(\boldsymbol{\nu}) := \boldsymbol{\nu}\cdot \tensC \cdot \boldsymbol{\nu}, 
\quad \text{with} \quad  \left(\boldsymbol{\nu}\cdot \tensC \cdot \boldsymbol{\nu}\right)_{jk}
 =   \sum_{i,l=1}^3   \nu_i \, c_{i jkl} \, \nu_l \, .
\end{equation}
From the symmetry of the stiffness tensor $\tensC \in \mathcal{S}_4$, \flo{the} matrix
$\boldsymbol{\Gamma}(\boldsymbol{\nu})$ is symmetric 
and \modif{positive definite}, and thus \flo{it} is diagonalizable with positive eigenvalues, 
cf. \cite{Pham2023hdgRR}.
The definition of the KC matrix can be rewritten in terms of the Voigt notation 
of $\tensC$ as follows, % which  is useful to compute the KC matrix for a generic stiffness tensor,
\begin{equation}\label{KCmat::Vrel}
  \boldsymbol{\Gamma}(\boldsymbol{\nu}) = \textcolor{black}{\mathbb{A}(\boldsymbol{\nu})} \,  \textcolor{black}{^\dagger\!  
  \overset{=}{\tensC}\,\! ^\dagger }\, \textcolor{black}{\mathbb{A}^t(\boldsymbol{\nu})},
\end{equation}
where, for a vector $\boldsymbol{\xi} \in \mathbb{C}^3$, we define
\begin{equation}\label{matbbA}\mathbb{A}(\boldsymbol{\xi}) := \sum_{I=x,y,z} \xi_I\, \mathbb{A}_I \,, \hspace*{0.8cm}
\mathbb{A}^\dagger(\boldsymbol{\xi}) := \sum_{I=x,y,z}\xi_I\, \mathbb{A}^\dagger_I\,,
\end{equation}
with the elementary matrices 
\begin{equation}\label{bbAmat::def}
%\small
\mathbb{A}_x =   % \small 
\begin{psmallmatrix}
   1 & 0 & 0  & 0 & 0 & 0  \\
   0 & 0 & 0  & 0 & 0 & \tfrac{1}{2} \\ 
   0 & 0 & 0 & 0 &  \tfrac{1}{2} & 0 
\end{psmallmatrix},\hspace*{0.2cm} % \normalsize
\mathbb{A}_y
 = % \small 
 \begin{psmallmatrix}
   0 & 0 & 0  & 0 & 0 & \tfrac{1}{2}  \\
   0 & 1 & 0  & 0 & 0 &  0 \\ 
   0 & 0 & 0 & \tfrac{1}{2} &  0 & 0 
\end{psmallmatrix},\hspace*{0.3cm} % \normalsize
 \mathbb{A}_z = % \small 
 \begin{psmallmatrix}
   0 & 0 & 0  & 0 &  \tfrac{1}{2} & 0  \\
   0 & 0 & 0  &  \tfrac{1}{2} & 0 & 0 \\ 
   0 & 0 & 1  & 0 & 0 & 0 
\end{psmallmatrix} .
\end{equation}
\flo{With matrix $D^\dagger$ of \cref{Ddagger}, we also define,}
\begin{equation}
\mathbb{A}^{\dagger}_I  =   \mathbb{A}_I \,  D^\dagger\,, \quad \text{where} \,\, I = x,y,z\,.
\end{equation}

%It is thus diagonalizable with eigenvalues associated with 
%wavespeeds in elasticity (up to three), cf. discussion in \cref{KCmat::subsubsec}. 
%For isotropic elasticity, the well-known speeds are the P- and S- wave speed, 
%$\cp$ and $\cs$, \cref{eq:elastic-isotropic-cp-cs}, 
%and $\boldsymbol{\Gamma}(\boldsymbol{\nu})$ has 
%two distinct eigenvalues: $\rho \, \cp^2$ and $\rho\,\cs^2$.
%
%
%In our work, with anisotropic elasticity, we restrict to the case of three 
% listed in descending order:
%$\rho \cqp^2$, $ \rho \cqsa^2$, and $\rho \cqsb^2$. 

% ------------------------------------------------
\subsection{Examples} 
% ------------------------------------------------

\paragraph{Isotropic elasticity} 
The stiffness tensor for an isotropic material has the following form in Voigt notation \cref{Voigt-s4tensor},
%\begin{subequations}
%\begin{align}
%\tensC_\mathrm{iso}  & = (c_{ijkl})\,, \quad  \text{with} \,\, 
%c_{ijkl} \, = \, \mu \left( \delta_{ik }\, \delta_{jl}\, + \, \delta_{il}\, \delta_{jk}\right) 
%\, + \, \lambda\, \delta_{ij} \, \delta_{kl} \,;\\[0.4em]
%%
%\tensCS_\mathrm{iso} &= (s_{ijkl})\,, \quad  \text{with} \,\, 
%s_{ijkl} \, = \, \dfrac{1}{4\mu} \, \left( \delta_{ik }\, \delta_{jl}\, + \, \delta_{il}\, \delta_{jk}\right) - \dfrac{ \lambda}{ 2 \mu ( 3 \lambda + 2 \mu)} \delta_{ij} \, \delta_{kl}   \,.
%\end{align}
%\end{subequations}
\begin{equation}
\label{iso_stiffness_voigt}
\overset{=}{\tensC}_{\mathrm{iso}} = \begin{psmallmatrix}  \lambda + 2\mu &  \lambda & \lambda & 0 & 0 & 0 \\ \lambda & \lambda + 2\mu & \lambda & 0 & 0 & 0 \\
\lambda & \lambda & \lambda + 2\mu & 0 & 0 & 0 \\ 0 & 0 & 0 & \mu & 0 &0 \\ 0 & 0 & 0 & 0 & \mu & 0 \\ 0 &0 & 0 & 0 & 0 & \mu \end{psmallmatrix}\,, 
\hspace*{0.5cm}
\overset{=}{\tensCS}_{\mathrm{iso}} = 
  \begin{psmallmatrix}  1/E &  - {\nu}/{E} & -{\nu}/{E} & 0 & 0 & 0 \\[0.2em]
- {\nu}/{E} & {1}/{E}   & - {\nu}/{E} & 0 & 0 & 0 \\[0.2em]
- {\nu}/{E} &-{\nu}/{E} &   {1}/{E  } & 0 & 0 & 0 \\ 0 & 0 & 0 & 
  {1}/{4\mu} & 0 &0 \\ 
  0 & 0 & 0 & 0 & {1}/{4\mu}  & 0 \\[0.2em]
  0 &0 & 0 & 0 & 0 & {1}/{4\mu}  
 \end{psmallmatrix} .
\end{equation}
Here the two degrees of freedom are given by \flo{the} Lam\'e parameters $\lambda$ and $\mu$, 
or by the Young modulus $E$ and Poisson coefficient (ratio),
\begin{equation}
E := \dfrac{ (3\lambda + 2\mu)\mu}{\lambda + \mu }, \hspace*{0.5cm}
\nu := \dfrac{\lambda}{2(\lambda + \mu)},  \hspace*{0.5cm}  \dfrac{\nu}{E} = \dfrac{ \lambda }{  2\mu(3\lambda + 2\mu)}.
\end{equation}
   \modif{With $\rho$ denoting the material density,} 
   the KC matrix for isotropic elasticity is
\begin{equation}\label{KC_iso}
    \boldsymbol{\Gamma}_\mathrm{iso}(\boldsymbol{\nu}) = \mu \identity + (\lambda + \mu) \boldsymbol{\nu} \otimes     
    \boldsymbol{\nu}= \rho \left( \cs^2 \mathbb{Id} + (\cp^2 - \cs^2)  
    \boldsymbol{\nu} \otimes \boldsymbol{\nu}  \right) .
    \end{equation}    
The second expression is written in terms of the P- and S-wave speed. 
The isotropic eigenvalues are independent of $\boldsymbol{\nu}$, 
\begin{equation}
   \mu = \rho  \cs^2  \, \text{  with multiplicity 2},  \hspace*{1cm}\text{and}  \hspace*{0.5cm}
   \lambda+\mu = \rho \cp^2 \,  \text{ with multiplicity 1}. 
\end{equation}

\paragraph{Vertical transverse isotropy (VTI)} 
Elastic material\flo{s} in this family are rotational symmetric around $\mathbf{e}_z$. 
Their stiffness tensor in Voigt representation \cref{Voigt-s4tensor} is given with five independent
coefficients, cf. \cite[Equation (1.39)]{carcione2007wave}, 
% \cite[Equation (1.39) in Subsection 1.2.1]{carcione2007wave},
    \begin{equation}\label{vti_stiffness_voigt}
  \overline{\overline{\tensC}}_{\mathrm{VTI}} = \begin{psmallmatrix}
                        C_{11}  & C_{11} - 2C_{66} & C_{13} & 0 & 0 & 0 \\
                        C_{11} - 2C_{66} & C_{11}  & C_{13} & 0 & 0 & 0 \\
                        C_{13} & C_{13}  & C_{33} & 0 & 0 & 0 \\
                        0 & 0 & 0 & C_{44} & 0 & 0 \\
                        0 & 0 & 0 & 0 & C_{44} & 0 \\
                        0 & 0 & 0 & 0 & 0 & C_{66} \\
                    \end{psmallmatrix} \,.
\end{equation}
%Positive definiteness of $\mathbf{C}$ requires that
%\begin{equation}
%c_{11} > \lvert c_{12} \rvert \,, \quad ( c_{11} + c_{12}) \, c_{33} \, > \, 2\, c_{13}^2 \,, \quad c_{55} > 0\,.
%\end{equation}\
In this case, the KC matrix takes the following form,
\begin{equation} \label{KC_vti}
 \boldsymbol{\Gamma}_\mathrm{VTI}(\n)
   =  \begin{psmallmatrix}
  C_{11} \,\nu_x^2 + C_{66} \,\nu_y^2 + C_{\modif{44}}\, \nu_z^2   \hspace*{0.3cm} &   \hspace*{0.3cm} (C_{13} + C_{\modif{44}}) \nu_y \nu_z  \hspace*{0.3cm}  &   \hspace*{0.3cm} (C_{13} + C_{\modif{44}})\, \nu_x \nu_z\\[0.3em]
      & C_{66} \nu_x^2 + C_{11}\, \nu_y^2 + C_{\modif{44}} \,\nu_z^2     &(C_{11}-C_{66})\,\nu_x \,\nu_y  \\[0.3em]
      &       &  C_{\modif{44}}(\nu_x^2 + \nu_y^2) + C_{33} \,\nu_z^2
   \end{psmallmatrix}.
\end{equation}

%\begin{remark}
%The components of the stiffness tensor for VTI can be also represented with Thomsen parameters sets, cf. \cite[Section 3.5.2]{Pham2023hdgRR}.
%\end{remark}

% ==============================================================================
\section{HDG formulations for time-harmonic linear elasticity}
\label{Formulation::sec}
% ==============================================================================

% =====================================================================
% \subsection{First-order systems of wave equations in linear elasticity}
% =====================================================================

We consider the propagation of time-harmonic waves in linear 
elasticity in terms of the displacement $\displacement$. 
For a domain $\Omega$ and an interior source $\mathbf{f}$, 
the wave equations are given by (e.g., \cite{carcione2007wave}),
\begin{subequations} \label{1stsystem_variants:u}
\begin{empheq}[left = {\text{Formulation } \formulationU \quad \empheqlbrace}\,]{align}
  & - \omega^2 \,\rho(\bx)\,\displacement(\bx,\omega) \, - \,  \nabla \cdot \tensS(\bx,\omega) \,=\, \mathbf{f}(\bx,\omega) \,, \\
  &   \tensCS(\bx) \tensS(\bx,\omega) =  \nabla^\mathfrak{s}\displacement(\bx,\omega)\,. \label{eq:constitutive-law-U}
\end{empheq}
\end{subequations}
%\begin{subequations}  \label{1stsystem_variants:v}
%\begin{empheq}[left = {\text{Formulation } \formulationV \quad \empheqlbrace}\,]{align}
%  & - \ii \omega \, \rho(\bx)\,\velocity(\bx,\omega) \, - \,  \nabla \cdot \tensS(\bx,\omega) \,=\, \mathbf{f}(\bx,\omega) \,, \\
%  &   \ii \omega \, \tensCS(\bx) \tensS(\bx,\omega) =  -\nabla^\mathfrak{s}\velocity(\bx,\omega)\,.
%\end{empheq} \end{subequations}
Here, $\tensS$ is the stress tensor, $\rho$ is the density, 
$\tensCS$ is the compliance tensor which is the inverse of the 
stiffness tensor $\tensCS:=\tensC^{-1}$.

\begin{remark}[Alternative formulations]
  Several variants to \cref{1stsystem_variants:u} have been introduced and studied in 
  the literature, \modif{containing $\omega$ or $\omega^2$,}
  % first or second-order in $\omega$, 
  and with possibly different choices 
  of unknowns such as the strain tensor $\tensE:=\nabla^\mathfrak{s}\displacement$. 
  For instance, in our report \cite{Pham2023hdgRR}, we compare the formulation 
  \cref{1stsystem_variants:u} 
  % (second-order in $\omega$ with displacement field $\displacement$) 
  with \flo{one in terms of $(\velocity, \tensS)$ with the velocity 
  $\velocity \,:=\, - \ii\omega\, \displacement$}, resulting in a system 
  % of first-order in $\omega$.
  \modif{containing $\omega$ instead of $\omega^2$.}
  We note that \cref{1stsystem_variants:u} gives the formulation for 
  elastostatics as $\omega\rightarrow 0$, and we have not observed \flo{any difference 
  between these two formulations}, cf.~\cite[Sections 6 and 7]{Pham2023hdgRR}.
 % Moreover, our choice of working with the compliance tensor $\tensCS$ rather
 % than the stiffness tensor $\tensC$ is motivated as it is more convenient to
 % allow for physical properties that vary within each element of the discretized
 % domain, see \cref{subsection:elastic-iso:sun}.
\end{remark}

\begin{remark}
We write \flo{the equation} in terms of the compliance 
tensor $\tensCS$ instead of the stiffness tensor 
$\tensC \,=\, \tensCS^{-1}$. 
This allows us to easily \flo{encode highly-varying backgrounds 
on meshes containing large cells, see \cref{subsection:elastic-iso:sun}. 
In this case, high-order polynomial discretization is to be employed},
% It allows us to represent the physical properties on 
% large-cell discretization, hence maintaining} high-order 
% polynomial discretization
which is critical for the efficiency of the HDG method 
compared to CG one, \cite{fabien2020gpu}.
Secondly, it allows us to have a flexible representation of the 
physical properties on the discretized domain \flo{as seen in 
our experiments \cref{section:numerics}}, and which is a key 
ingredient for inversion, \cite{Faucher2020adjoint}.
\end{remark}

% ================================================================
\subsection{Discretization domain}
\label{dommesh::subsec}
% ================================================================

For numerical resolution, \modif{we employ a conforming triangulation of the domain $\Omega$, 
denoted by $\mathcal{T}_h$, which consists of} non-overlapping elements $K^e$, 
$1\leq e\leq \lvert \mathcal{T}_h\rvert$,
\begin{equation}\Omega\, = \, \bigcup_{e=1}^{\lvert \mathcal{T}_h\rvert} K^e :=  \Omega_h, \hspace*{1cm}
\mathcal{T}_h = \{ K^e \big| 1\leq e\leq \lvert \mathcal{T}_h\rvert\}.
\end{equation} 
The set of faces in the mesh $\mathcal{T}_h$, denoted by $\Sigma_h$,
 consists of %the union of the 
 faces $\mathsf{F}$ which also form the 
 boundary of all elements $K$ in $\mathcal{T}_h$, 
\begin{equation}
\Sigma_h  = \bigcup_{e=1}^{\lvert \mathcal{T}_h\rvert} \partial K^e  
 = \{\mathsf{F}^k \big| 1\leq k\leq  \lvert \Sigma_h\rvert\}.
\end{equation}
We distinguish between the interior faces, also called \emph{interfaces}, $ \Sigma_{\mathrm{int}}$
and boundary ones $\Sigma_\partial$ which form the boundary of $\Omega_h$:
\begin{equation}\Sigma_h = \Sigma_{\mathrm{int}} \cup \Sigma_\partial, \hspace*{0.5cm}
  \text{with} \hspace*{0.2cm}\Sigma_\partial = \partial \Omega_h =
  \Sigma_{\mathrm{N}} \, \cup \, \Sigma_{\mathrm{D}} \, \cup \, \Sigma_{\infty}\,.
  \end{equation}
The latter \flo{set} is partitioned into non-overlapping
regions on which  Neumann, Dirichlet, and Robin 
boundary conditions can be imposed. In \flo{our experiments}, 
we only consider tetrahedal meshes in 3D and triangle meshes in 2D.

\paragraph{Local and global indexes of an edge} 
An edge $\mathsf{F} \in \Sigma_h$ can be referred to in two ways,
\begin{enumerate}
\item as $\mathsf{F}^k$, with $1\leq k\leq \lvert \Sigma_h\rvert$, here $k$ is the index of the face 
      $\mathsf{F}$ in the ordered set $\Sigma_h$, 
\item as $\mathsf{F}^{(e,\ell)}$  with $1\leq e\leq \lvert \mathcal{T}_h\rvert$, 
$1\leq \ell\leq \nfe $, when it is considered as part of the boundary of an element $K^e$, and $\ell$ is its index among the set of faces of $K^e$. Note that $\nfe$ 
% is the number of edges in element $K^e$, 
% for the type of mesh we work with
\modif{is the number of edges 
in 2D, and the number of faces in 3D.
Working with simplexes, we have,}
\begin{equation}
\nfe=4 , \hspace*{0.2cm}   \text{in 3D}, \hspace*{1cm}  \hspace*{0.2cm}\nfe=3
, \hspace*{0.2cm}   \text{in 2D}.
 \end{equation}
\end{enumerate}
%We denote by $\mathfrak{f}$ the mapping which relates these local and global labels of an edge:
%\begin{equation}\label{betamap::def}
%  \mathfrak{f} \quad : \quad    (e,\ell)\qquad \mapsto \quad  k\,  =\,  \mathfrak{f}(e,\ell)\,, \quad \text{where} \quad \mathsf{F}^{(e,\ell)}\,  =\,  \mathsf{F}^{k}\,. 
%\end{equation}

\paragraph{Jump operator} 
At \flo{an} interface $\mathsf{F}$, the jumps $\llbracket \,\cdot\, \rrbracket$ 
of a vector $\velocity$ and of the traction of a matrix $\boldsymbol{\chi}$ 
defined on $\overset{\circ}{K^+} \cup \overset{\circ}{K^-}$ 
with $\overset{\circ}{\phantom{K}}$ denoting the interior, are respectively defined as
\begin{equation}
\llbracket \velocity \rrbracket := \velocity|_{K^+} - \velocity \big|_{K^-}
\,, \qquad\qquad \llbracket \, \boldsymbol{\chi}\, \boldsymbol{\nu} \, \rrbracket
 := \boldsymbol{\chi}|_{K^+} \boldsymbol{\nu}^+ + \boldsymbol{\chi}|_{K^-} \boldsymbol{\nu}^- \,,
\end{equation}
with  $\boldsymbol{\nu}^{\pm}$ the outward-pointing normal vector defined on $\partial K^{\pm}$.

% ===================================================================
\subsection{Statement of HDG problem for formulation $\formulationU$}
% ===================================================================

% In order to reduce the computational cost due to the doubling of 
% degrees of freedom at an interior edge in a DG method, HDG 
% formulation employs two strategies:  
%HDG formulation employs two strategies:  
%the \emph{hybridization} of mixed method, and static condensation. 
\modif{In the HDG method, the original global problem is statically condensed into a problem 
in terms of a hybrid variable defined on the skeleton of the mesh \cite{hungria2017hdg}.}
\modif{Specifically, it is first} written as a union of local boundary-valued problems defined 
for each mesh cell, \modif{having as boundary data the hybrid variable}. 
\modif{Local problems are linked} by transmission conditions (or jump conditions)
which constrain weakly \modif{the continuity of solutions}.
\modif{In this way}, a global problem is \modif{obtained and solved in terms of the hybrid variable;
after this,} the primal and mixed unknowns are found
\modif{by solving the} local problems, which can be realized in a parallel manner.
\modif{For a more in-depth discussion, we refer the readers to \cite{giacomini2021hdglab}.}
%as they are independent on each element of the discretized domain. 

\paragraph{Strong form} 
% The HDG method comprises of local problems which are boundary value 
\flo{The} local problems are boundary value problems defined 
on each cell $K \in \mathcal{T}_h$ with prescribed Dirichlet condition 
$\lambdaU\in L^2(\partial K)$. From \cref{1stsystem_variants:u}, we have,
\begin{subequations}\label{strong_local}
 \begin{empheq}[left = {\empheqlbrace}\quad ]{align}
& -\omega^2 \, \rho \, \displacement \, - \,  \nabla \cdot \tensS \,=\, \mathbf{f} \,, 
       \qquad \quad \text{on} \,\, K\,, \label{strong_local_mot}\\[0.3em]
&  \mathbf{S}\tensS \,=\, \nabla^\mathfrak{s} \displacement \,, 
       \hspace*{3.05cm} \text{on} \,\, K \,, \label{strong_local_consti}\\[0.3em]
& \displacement= \lambdaU \,, \hspace*{3.70cm} \text{on}\,\, \partial K  \label{strong_localBC}\,.
\end{empheq}
\end{subequations}
The trace  $\lambdaU$ is determined by the 
\modif{transmission} condition,
\begin{equation}\label{strong_conser}
\llbracket \tensS\boldsymbol{\nu}\rrbracket = 0
\, \quad \text{on } \,\,  \mathsf{F} \in \Sigma_{\mathrm{int}},
\end{equation}
and \flo{the} boundary conditions on the boundary of the domain $\Sigma_\partial$ are given by, 
\begin{equation}\label{strong_BC}
\begin{matrix*}[l]
 \lambdaU = \boldsymbol{\mathsf{g}}^\mathrm{D} \,                       & \text{on }  \,\, \Sigma_\mathrm{D}\,, \qquad\qquad & \text{Dirichlet boundary, } \\[.40em]
  \tensS\boldsymbol{\nu} = 0  \,\,\,\,  & \text{on }  \,\,\Sigma_\mathrm{N}\,,\qquad\qquad & \text{Neumann boundary, } \\[.40em]
  \tensS \boldsymbol{\nu} = -\mathrm{i}\omega \mathcal{Z} \, \lambdaU \qquad % \boldsymbol{\nu} \qquad 
                                  & \text{ on } \Sigma_\infty\,, \qquad\qquad & \text{Robin boundary.}
\end{matrix*}
\end{equation}
The Robin condition employs impedance-like matrix
$\mathcal{Z}= (\mathcal{Z}_{IJ})_{I,J=x,y,z}$.

% -------------------------------------------------------------
\paragraph{Finite element spaces} Below, tensor-valued 
functions whose components taking value in a function space 
$\mathcal{V}$ are written as, 
\begin{equation} \label{eq:definition-Vfunc}
\text{vector-valued}\hspace*{0.3cm}  \mathcal{V}^n \,, \hspace*{0.5cm} 
\text{matrix-valued}\hspace*{0.3cm} \mathcal{V}^{n \times n } \,, \hspace*{0.5cm} \substack{\text{symmetric}\\\text{ matrix-valued}} \hspace*{0.3cm}\mathcal{V}^{n\times n}_\mathrm{sym}\,.
\end{equation}
%Below we write $ L^2(\Omega_h)^{3\times 3}_\mathrm{sym} $ 
%to denote the space of symmetric $3\times 3$ matrices 
%(thus $\boldsymbol{\chi} = (\chi_{IJ})_{I,J=x,y,z,}$ 
%with $\chi_{IJ} = \chi_{JI}$).  
We introduce the global finite element spaces 
\begin{equation}\label{global_fem_spaces}
\begin{aligned}
U_h  &= \{ \displacement = (u_I)_{I=x,y,z}\in L^2(\Omega_h)^3 \hspace*{1cm}: \hspace*{0.5cm}  u_I|_K   \in U_h(K)   , \,\, \forall K \in \mathcal{T}_h\} \,;\\[0.5em]
V_h &= \{ \boldsymbol{v} = (v_{IJ})_{I,J=x,y,z} \in L^2(\Omega_h)^{3\times 3}_\mathrm{sym} \hspace*{0.1cm} :\hspace*{0.5cm} v_{IJ}|_{K}   \in V_h(K)   , \,\, \forall K \in \mathcal{T}_h\}\,;\\[0.5em]
M_h &= \{  \boldsymbol{w}=(w_I)_{I=x,y,z} \in L^2(\Sigma_h)^3 \hspace*{0.95cm}:\hspace*{0.5cm}   \mu_I|_{\mathsf{F}}   \in M_h(\mathsf{F})   , \,\, \forall \mathsf{F} \in \Sigma_h\}\,.
\end{aligned}
\end{equation}
\modif{We employ polynomials of equal degree in \cref{global_fem_spaces}, as was done in \cite{soon2009hybridizable,bonnasse2018hybridizable,fernandez2018hybridized}; specifically,
with} $\mathcal{P}^\mathsf{k}(D)$ denoting the space of polynomials of degree $\mathsf{k}$ defined on domain $D \subset \mathbb{R}^3$, 
\begin{equation}
  U_h(K) \,=\, V_h(K) \,=\, 
  \mathcal{P}^\mathsf{k}(K) \,, 
  \qquad M_h(\mathsf{F}) = \mathcal{P}^\mathsf{k}(\mathsf{F})\,.
\end{equation}

\medskip

\paragraph{Approximate problem}
To obtain the discrete problem,
we integrate the local problem \cref{strong_local} against test functions $
(\boldsymbol{\phi} , \Psi ) \in U_h \times V_h$,
and \flo{the} problems on \flo{the} interfaces against test functions $\boldsymbol{\upxi} \in M_h$, 
and carry out integration by parts. 
In the weak form associated with \cref{strong_local_mot}, we also employ the following definition for 
the numerical trace of \flo{the} traction, 
\begin{equation}\label{numtraction_u}
  \widehat{\tensS \boldsymbol{\nu}} \, = \, \tensS_h \boldsymbol{\nu} \, \textcolor{black}{-} \,\tauU (\displacement_h -  \lambdaUh) \,,\hspace*{0.5cm} \text{ with stabilization matrix } \tauU = (\tau_{IJ})_{I,J=x,y,z},
\end{equation} 
and carry out an inverse integration by parts. 
For more details in obtaining \cref{app_local_usig_Hdg}, see \cite{Pham2023hdgRR}.
% \flo{The explicit values of $\tauU$ considered in this work are given in \cref{subsection:stabilization}}.

\medskip

The approximate problem reads as follows: Find $(\textcolor{black}{\displacement_h},\textcolor{black}{\tensS_h},\textcolor{black}{ \lambdaUh }) \in U_h \times V_h \times M_h$ that solves, 
\begin{enumerate}[leftmargin=*]
\item Local problems on element $K^e \in \mathcal{T}_h$, for all test functions $
(\boldsymbol{\phi} , \Psi ) \in U_h \times V_h$, 
\begin{subequations}\label{app_local_usig_Hdg}
\begin{empheq}[left = { \empheqlbrace}]{align}
&-\omega^2 \hspace*{-.30em} \int_{K} \hspace*{-.20em}  \rho\textcolor{black}{ \displacement_h }\cdot \overline{\boldsymbol{\phi}}\,\mathrm{d} \mathbf{x}
 - \int_{K}  (\nabla\cdot\textcolor{black}{\tensS_h}) \cdot \overline{\boldsymbol{\phi}} \,\mathrm{d} \mathbf{x}
  \textcolor{black}{+} \int_{\partial K} \overline{\boldsymbol{\phi}} \cdot \textcolor{black}{\tauU}(\textcolor{black}{ \displacement_h}- \textcolor{black}{\lambdaUh})\, \mathrm{d} s_{\mathbf{x}}
  =  \int_{K} \mathbf{f}\cdot \overline{\boldsymbol{\phi}}\,  \mathrm{d} \mathbf{x}, \label{app_local_usig_Hdg::eqnmot}\\[0.3em]
& 
\hspace*{0.5cm} \int_{K}  \textcolor{black}{\tensS_h} : \mathbf{S} : \overline{\Psi}\,  \mathrm{d}\mathbf{x}  =  - \int_{K}\textcolor{black}{\displacement_h }  \cdot  \nabla\cdot   \overline{\Psi} \mathrm{d}\mathbf{x}
+ \int_{\partial K} 
\boldsymbol{\nu} \cdot  \overline{\Psi}  \cdot \modif{\lambdaUh}\,  \mathrm{d} s_{\mathbf{x}}\,.\label{app_local_usig_Hdg::conlaw}
\end{empheq}
\end{subequations}

\item Interface and boundary problems: 
\modif{the transmission} and Neumann conditions give, 
\begin{equation}\label{intface_neu_HDGprob}
\begin{aligned}
&\sum_{ e=1}^{\lvert \mathcal{T}_h\rvert} \int_{\partial K^e \cap( \Sigma_{\mathrm{int}} \cup \Sigma_\mathrm{N}) } \left(
\boldsymbol{\nu} \cdot \textcolor{black}{\tensS_h} \cdot\overline{\boldsymbol{\xi}}  \hspace*{0.1cm} \textcolor{black}{-} \hspace*{0.1cm}  \overline{\boldsymbol{\xi}} \cdot\tauU \cdot \left(\textcolor{black}{\displacement_h }
-\textcolor{black}{\lambdaUh } \right) \right) \, \mathrm{d}s_{\mathbf{x}} \, = \, 0\,, 
 \hspace*{0.2cm} \,\, \forall \boldsymbol{\upxi} \in M_h(\Sigma_{\mathrm{int}} \cup \Sigma_\mathrm{N}) \,.
\end{aligned}
\end{equation}
The Robin boundary conditions~\cref{strong_BC} give, $ \forall \boldsymbol{\upxi} \in M_h(\Sigma_{\infty})$,
\begin{equation}\label{abc_HDGprob}
\begin{aligned}
&\sum_{ e=1}^{\lvert \mathcal{T}_h\rvert} \int_{\partial K^e \cap \Sigma_{\infty}} \left(
\boldsymbol{\nu} \cdot \textcolor{black}{\tensS_h} \cdot\overline{\boldsymbol{\xi}} \hspace*{0.2cm} \textcolor{black}{-} \hspace*{0.2cm}  \overline{\boldsymbol{\xi}} \cdot
\tauU \cdot\left(\textcolor{black}{\displacement_h}-\textcolor{black}{\lambdaUh} \right)
 \hspace*{0.1cm}\, \modif{+ \, \ii \omega} \,
\overline{\boldsymbol{\xi}} \cdot \mathcal{Z}\, \cdot \textcolor{black}{\lambdaUh } \right) \mathrm{d}s_{\mathbf{x}} \, = \, 0 \,.
\end{aligned}
\end{equation}
The Dirichlet condition is imposed weakly, $ \forall \boldsymbol{\upxi} \in M_h(\Sigma_{\mathrm{D}})$,
\begin{equation}\label{homDirc_cond}
\sum_{ e=1}^{\lvert \mathcal{T}_h\rvert} \int_{\partial K^e \cap \Sigma_{\mathrm{D}}}
\left( \textcolor{black}{\lambdaUh}  - \boldsymbol{\mathsf{g}}^\mathrm{D} \right) \cdot \overline{\boldsymbol{\xi}}\,\mathrm{d}s_{\mathbf{x}} =  0 \,, \quad \text{on }  \,\, \Sigma_\mathrm{D}\,.
\end{equation}
\end{enumerate}

%%---------------------------------------
\subsection{Discretization of HDG problem}
\label{discretization::sec}
%%---------------------------------------

The discussion in this section employs the Voigt notation introduced in  \cref{VI::subsec}.

% ==============================================================
\subsubsection{Discrete unknowns}
\label{notation_discre::subsec}
% ==============================================================

\paragraph{Local Basis functions} We have three groups of basis functions 
for the local finite element spaces introduced in \cref{global_fem_spaces},
\begin{equation}\label{localbasisfcn}
\begin{array}{ c| c |c} \text{Local finite element space} &  \text{Basis functions}  & \text{Dimension}  \\[0.2em]\hline
 U_h(K^e)  ,  \,\, 1\leq e\leq \lvert \mathcal{T}_h\rvert &   \phi_j^e \,, \,\, 1\leq j\leq \dofpe & \dofpe \\[0.3em] 
V_h(K^e) ,  \,\,  1\leq e\leq \lvert \mathcal{T}_h\rvert\quad & \psi_j^e \,, \,\, 1\leq j\leq \dofve  & \dofve\\[0.3em]
M_h(\mathsf{F}^k) =M_h(\mathsf{F}^{(e,\ell)})  ,  \,\,  1\leq k\leq \lvert \Sigma_h\rvert\hspace*{0.2cm} & \hspace*{0.2cm}  \xi^k_j=\xi^{(e,\ell)}_j \,, \,\, 1\leq j\leq \hat{\mathrm{n}}_k \hspace*{0.2cm} & \hspace*{0.2cm}  \dofedge_k = \dofedge_{(e,\ell)}\   \,.  \end{array}
 \end{equation}
\flo{We} also denote the total number of face degrees of freedom in each direction $I$ by, 
\begin{equation}
\dofglobal_I =\dofglobal
 = \sum_{k  =1}^{|\Sigma_h|} \dofedge_k \,, \hspace*{0.2cm} I=x,y,z\,.
\end{equation}

\paragraph{Volume discrete unknowns}
The displacement 
vector field $\displacement$ and  strain tensor 
$\tensS$ are approximated on each cell $K^e$, $1\leq e\leq \lvert \mathcal{T}_h\rvert$ as follows.
For the displacement we have,
\begin{equation}
{\displacement}^h = \sum_{I =x,y,z} \textcolor{black}{u^h_I}\hspace*{0.1cm} \hat{\mathbf{e}}_I\,, \quad \text{with} \quad 
\textcolor{black}{u^h_I} \big|_{K^e} = \sum_{j=1}^{\dofpe} \textcolor{black}{ \mathsf{u}^e_{Ij} }\, \phi^e_j
\,, \quad \text{for} \,\, I = x,y,z \,.
\end{equation}
For the strain tensor, we work with
its Voigt representation as a vector of length 6, cf. \cref{decomsigma_v},
%\begin{equation}
%\begin{aligned}
%& \tensS_h = \sum_{\frak{J} \in \mathcal{I}_{\mathrm{sm}}}   \textcolor{black}{\sigma_{\frak{J}}^h}  \,  \overset{=}{\mathbb{e}}_{\frak{J}}\hspace*{1cm} \longleftrightarrow \hspace*{1cm} \overrightarrow{\tensS}_h = \sum_{\frak{J} \in  \mathcal{I}_{\mathrm{sm}} }\textcolor{black}{\sigma^h_{\frak{J}} } \hspace*{0.1cm} \hat{\mathbb{e}}_{\frak{J}}\,,  \\
%%
%& \hspace*{2cm} \text{with} \hspace*{0.7cm} \textcolor{black}{ \sigma^h_{\frak{J}}}\big|_{K^e} = \sum_{j=1}^{\dofve} \textcolor{black}{\mathsf{v}^e_{\frak{J}j}} \hspace*{0.1cm} \psi^e_j\,, \quad \text{for} \,\, \frak{J} \in \mathcal{I}_{\mathrm{sm}} \,.
%\end{aligned}
%\end{equation}
\begin{equation}
\begin{aligned}
 \tensS_h = \sum_{\frak{J} \in \mathcal{I}_{\mathrm{sm}}}   \textcolor{black}{\sigma_{\frak{J}}^h}  \,  \overset{=}{\mathbb{e}}_{\frak{J}}\hspace*{0.5cm} \longleftrightarrow \hspace*{0.5cm} \overrightarrow{\tensS}_h = \sum_{\frak{J} \in  \mathcal{I}_{\mathrm{sm}} }\textcolor{black}{\sigma^h_{\frak{J}} } \hspace*{0.1cm} \hat{\mathbb{e}}_{\frak{J}}\,,  
 \hspace*{0.2cm} \text{with} \hspace*{0.2cm} \textcolor{black}{ \sigma^h_{\frak{J}}}\big|_{K^e} = \sum_{j=1}^{\dofve} \textcolor{black}{\mathsf{v}^e_{\frak{J}j}} \, \psi^e_j\,\,.
\end{aligned}
\end{equation}
We gather all volume unknowns into vector $\boldsymbol{\mathsf{W}}$ such that,
\begin{equation}\label{vectorW_e}
\begin{aligned}
\boldsymbol{\mathsf{W}} = \left(\boldsymbol{\mathsf{W}}^e \right)_{1\leq e\leq \lvert \mathcal{T}_h\rvert}
\,, \quad &\text{with} \quad 
\boldsymbol{\mathsf{W}}^e : = \begin{pmatrix} \boldsymbol{\mathsf{U}}^e \\ \boldsymbol{\mathsf{V}}^e \end{pmatrix} \hspace*{0.2cm}\text{and} \quad \boldsymbol{\mathsf{U}}^e = \left(\mathsf{U}^e_I\right)_{I=x,y,z}\,, \,\, \boldsymbol{\mathsf{V}}^e =  \left( \mathsf{V}^e_\frak{J}\right)_{\frak{J}\in \mathcal{I}_\mathrm{sm}}\,,\\
&\hspace*{0cm} \text{with sub-blocks} \quad\mathsf{U}^e_I = \left( \textcolor{black}{\mathsf{u}^e_{Ij}}\right)_{j=1,\ldots,\dofpe}\,, \quad  \mathsf{V}^e_\frak{J} = \left(\textcolor{black}{ \mathsf{v}^e_{\mathfrak{J}j}}\right)_{j=1,\ldots,\dofve} \,.
\end{aligned}
\end{equation}

\paragraph{Edge discrete unknowns} 
\flo{The} unknowns defined on the edges $\Sigma_h$ are numerical approximations 
of the trace of the displacement ${\displacement}$ on the edges $\Sigma_h$. 
Specifically, on each edge $\mathsf{F}^k$, $k=1, \ldots, \dofedge_k$,
\begin{equation}
\boldsymbol{\uplambda}_h = \sum_{I \in \{x,y,z\}} \textcolor{black}{ \uplambda^h_{I}} \, \hat{\mathbf{e}}_I
\,, \quad \text{with} \quad
\textcolor{black}{\uplambda^h_{I}}\big|_{\mathsf{F}^k} = \sum_{j=1}^{\dofedge_k} \textcolor{black}{\lambda^{k}_{Ij}}  \hspace*{0.1cm}  \xi^k_j \,.
\end{equation}
\flo{We} gather these coefficients into \flo{a} global vector $\boldsymbol{\Lambda}$ which has the following substructure,
\begin{equation}\label{vectorLambda}
\boldsymbol{\Lambda} =  \left(\boldsymbol{\Lambda}^k \right)_{k=1,\ldots,\lvert \Sigma_h\rvert} 
\,, \quad\text{with} \quad  \boldsymbol{\Lambda}^k = \left(  \Lambda^k_I \right)_{I=x,y,z}
 \quad \text{and} \quad  \Lambda^k_I= \left( \lambda^k_{Ij}\right)_{1\leq j\leq\dofedge_k }\,.
\end{equation}

\begin{remark}
The above expressions are written \flo{in terms of the} global edge indices. 
They can also be written in
terms of the local edge indices, specifically
for $(e,\ell)$  identified 
with $k$, then
\begin{equation}
\mathsf{F}^k = \mathsf{F}^{(e,\ell)}
\,, \quad \lambda^{(e,\ell)}_{Ij} = \lambda^{k}_{Ij} 
\,, \quad \xi^{(e,\ell)}_j = \xi^k_j\,.
\end{equation}
\end{remark}

%% --------------------------------------------------------
%\subsection{Discretization for elastic wave formulation $\formulationU$}
%\label{dis_u_sigma_S::subsec}
%% --------------------------------------------------------

\subsubsection{Discretized local problem}

The discretized problem defined on each element $K^e$ \flo{is written as},
  \begin{equation}\label{dis_localprob_v0}
 \mathbb{A}^e \, \boldsymbol{\mathsf{W}}^e \, + \, \mathbb{D}^e \, \boldsymbol{\Lambda}^e \, = \, \boldsymbol{\mathsf{S}}^e\,, \hspace*{1.5cm} \text{for} \,\, 1 \leq e\leq \lvert \mathcal{T}_h\rvert ,
 \end{equation}
 with coefficient matrices and sources given by,
 \begin{equation}
 \mathbb{A}^e = \begin{pmatrix} -\omega^2\,  \mathbb{M}_{{\displacement}}^e  \,  \textcolor{black}{+}\,  \mathbb{M}_{\partial}^e  \,\, &\,\,  -\mathbb{K}_{\tensS}^e \\[0.5em]
 -\mathbb{K}^e_{{\displacement}} & -\mathbb{M}^e_{\tensS}
 \end{pmatrix}, \hspace*{0.2cm} \mathbb{D}^e = \begin{pmatrix}  \textcolor{black}{-}\mathbb{D}^{(e,1)}_{\mathrm{m}} & \ldots & \textcolor{black}{-}\mathbb{D}_{\mathrm{m}}^{(e,\nfe)} \\[0.3em]
\mathbb{D}_\mathrm{c}^{(e,1)}  & \ldots & \mathbb{D}_\mathrm{c}^{(e,\nfe)}  \end{pmatrix} ,\hspace*{0.2cm} \boldsymbol{\mathsf{S}}^e  =  \begin{pmatrix}
 \boldsymbol{\mathsf{S}}^e_\mathrm{m} \\[0.1em]
 \boldsymbol{0}_{ (6  \mathsf{m}_e)   \times 1} \end{pmatrix}.
\end{equation}
 We list below the components of the sub-blocks.
 
 \paragraph{Description of components}
 The non-zero entries of \flo{the} source term come from integrating the 
 volume source $\mathbf{f}$ in the equation of motion \flo{with test functions}:
  \begin{equation}
 \boldsymbol{\mathsf{S}}^e_\mathrm{m} = \left(\hspace*{-0.3em} ~_\mathrm{m}\boldsymbol{\mathsf{S}}^e_I \right)_{I=x,y,z}
 ,  \hspace*{0.4cm} \text{with} \hspace*{0.3cm} 
  \left[\hspace*{-0.3em}~_\mathrm{m}\boldsymbol{\mathsf{S}}^e_I \right]_i = \left<  \mathbf{f} |_{K^e} ,     \phi^e_i\, \hat{\mathbf{e}}_I \right>\,, \,\, 1\leq i\leq \dofpe\,.\label{vol_src}  \end{equation}
The upper row blocks of $\mathbb{A}^e$ and $\mathbb{D}^e$ come from integrating the equation of motion
\cref{app_local_usig_Hdg::eqnmot}, % against test functions,
\begin{equation}
\boldsymbol{\phi} = \phi^e_i \, \hat{\mathbf{e}}_I, \qquad \,\, I=x,y,z\, \,\, 1\leq i\leq \dofpe \,,
\end{equation}
while those of the second row-block come from\footnote{To describe the components of the matrices in Voigt quantities, the following identities are employed: for $\boldsymbol{\sigma}\in S_2$ and vector $\mathbf{w}$, 
\begin{equation}
 (\nabla\cdot \boldsymbol{\sigma}) \cdot \mathbf{e}_I
 = \sum_{\frak{J} \in \mathcal{I}_{\mathrm{sm}}} \mathbf{e}_I    \cdot  \mathbb{A}^\dagger(\nabla \sigma_{\frak{J}})  \cdot \hat{\mathbb{e}}_{\frak{J}}, \hspace*{0.7cm}
  (\tensCS \boldsymbol{\sigma})  : {\boldsymbol{w}} = \boldsymbol{\sigma} \cdot  \tensCS  \cdot {\boldsymbol{w}}
 = \overrightarrow{\boldsymbol{w}} \cdot ^\dagger\!\overset{=}{\tensCS}\,\! ^\dagger  \cdot \overrightarrow{\boldsymbol{\sigma}}\,.
\end{equation}} 
integrating the constitutive equation \cref{app_local_usig_Hdg::conlaw}
against test functions,
\begin{equation}\Psi  = \psi^e_i \, \overset{=}{\mathbb{e}}_\frak{I}  \hspace*{0.2cm} \underset{\text{Voigt identification}}{\leftrightarrow} \hspace*{0.2cm} \overrightarrow{\Psi} = \psi^e_i \, \hat{\mathbb{e}}_{\mathfrak{I}},  \hspace*{0.7cm} \text{for} \,\,  \mathfrak{I} \in \mathcal{I}_{\mathrm{sm}}\, , \,\,1\leq i\leq \dofve.
\end{equation}
 For $1\leq i,j\leq \dofpe$, 
\begin{subequations}
\begin{align}
\mathbb{M}_{{\displacement}}^e  &= \left( \mathbb{M}_{IJ}^{{\displacement}e} \right)_{I,J=x,y,z},  \hspace*{0.2cm} \text{with} \hspace*{0.3cm} \left[\mathbb{M}^{{\displacement} e}_{II} \right]_{ij}  = 
\begin{dcases}   \int_{K^e} \rho \, \overline{\phi^e_i} \, \phi^e_j \, \mathrm{d} \mathbf{x}  , & I=J\\
  \hspace*{1cm} 0, &I\neq J \end{dcases},   \\
 \mathbb{M}_{\partial}^e  &= \left( \mathbb{M}_{IJ}^{\partial e} \right)_{I,J=x,y,z}, \quad \text{with} \,\, \left[\mathbb{M}^{\partial e}_{IJ}\right]_{ij} =  \sum_{\ell=1}^{\nfe }\int_{\mathsf{F}^{(e,\ell)} }
\,   \textcolor{black}{\tau_{IJ}} \,  \overline{  \phi^{e}_i}   \, \phi^e_j \,  \mathrm{d}s_{\mathbf{x}} .\label{matMpartial} 
\end{align}
  \end{subequations}
  For $1\leq i,j\leq \dofve$,
  \begin{equation}
\mathbb{M}^e_{\tensS} = \left( \mathbb{M}^{\tensS e}_{\mathfrak{J} \mathfrak{J}'}\right)_{\mathfrak{J},\mathfrak{J}'\in \mathcal{I}_{\mathrm{sm}}}  \hspace*{0.2cm} \text{with} \hspace*{0.3cm}\left[ \mathbb{M}^{\tensS e}_{\mathfrak{J}\mathfrak{J}'}\right]_{ij}
 =   \hat{\mathbb{e}}_\frak{J}  \cdot \left(\int_{K^e}  \overline{\psi^e_i} \, \psi^e_j\,  ^\dagger\!\overset{=}{\mathbf{S}}\,\! ^\dagger\, \mathrm{d}\mathbf{x} \right)  \cdot \hat{\mathbb{e}}_{\frak{J}'} .
\end{equation}
  For $1\leq i\leq \dofpe, 1\leq j\leq \dofve$, 
  \begin{subequations}
\begin{align}
  \mathbb{K}_{\tensS}^e &=  \left( \mathbb{K}^{\tensS e}_{I\mathfrak{J}}\right)_{
  \substack{I=x,y,z,\\ \mathfrak{J}\in \mathcal{I}_{\mathrm{sm}}}} \hspace*{0.2cm}\text{with} \hspace*{0.2cm} 
\left[\mathbb{K}^{\tensS e}_{I\mathfrak{J}}\right]_{ij} = \hat{\mathbf{e}}_I \cdot \left(   \int_{K^e}     \overline{\phi^e_i}  \hspace*{0.1cm} \mathbb{A}^\dagger(\nabla \psi^e_j )\,\mathrm{d} \mathbf{x}  \right)\cdot  \hat{\mathbb{e}}_\frak{J}, \\
 \mathbb{K}_{{\displacement}}^e &= \left( \mathbb{K}^{{\displacement}e}_{\frak{I}J}\right)_{\substack{\mathfrak{I} \in \mathcal{I}_{\mathrm{sm}},\\J=x,y,z}} \hspace*{0.2cm} \text{with} \hspace*{0.3cm}\left[ \mathbb{K}^{{\displacement} e}_{\frak{I}J}\right]_{ij}  = \hat{\mathbf{e}}_J \cdot  \left(\int_{K^e}  \phi^e_j  \, \, \overline{\mathbb{A}^\dagger(\nabla\psi^e_i)} 
\, \mathrm{d} \mathbf{x} \right)\cdot \hat{\mathbb{e}}_{\frak{J}} .
\end{align}
\end{subequations} 
For $ 1\leq i\leq \dofve , 1\leq j\leq \dofedge_{(e,\ell)} $
\begin{equation}
\mathbb{D}^{(e,\ell)}_\mathrm{m} = \left(\hspace*{-0.5em}  ~_{\mathrm{m}}\mathbb{D}^{ (e,\ell)}_{ IJ}\right)_{I,J=x,y,z} \hspace*{0.2cm} \text{with} \hspace*{0.1cm} \left[\hspace*{-0.4em}  ~_{\mathrm{m}}\mathbb{D}^{(e,\ell)}_{ IJ} \right]_{ij} = \int_{\mathsf{F}^{(e,\ell)} }
 \textcolor{black}{\tau_{IJ}} \,  \overline{  \phi^{e}_i}   \, \xi_j^{(e,\ell)} \,  \mathrm{d}s_{\mathbf{x}}.\label{Dm}
 \end{equation}
For $1\leq i\leq \dofpe ,   1\leq j\leq \dofedge_{(e,\ell)}$,
\begin{equation}
\mathbb{D}^{(e,\ell)}_\mathrm{c} = \left( \hspace*{-0.5em}  ~_{\mathrm{c}}\mathbb{D}^{(e,\ell) }_{\frak{I}J}\right)_{
\substack{\frak{I} \in \mathcal{I}_{\mathrm{sm}},\\  J=x,y,z}} \hspace*{0.1cm} \text{with} \hspace*{0.1cm}\left[\hspace*{-0.4em}  ~_{\mathrm{c}}\mathbb{D}^{(e,\ell)}_{ \frak{I}J} \right]_{ij}
 \! = \!  \int_{\mathsf{F}^{(e,\ell)}} \overline{\psi^e_i} \, \xi^{(e,\ell)}_j 
\mathbb{A}^\dagger(\boldsymbol{\nu}^{(e,\ell)})_{J \frak{J}}  \, \mathrm{d}s_{\mathbf{x}} .  \label{Dc} 
 \end{equation}

% =========================================================
\subsubsection{Discretized problem defined on edges}
% =========================================================

The discretization of 
\cref{intface_neu_HDGprob}--\cref{homDirc_cond} defined on $\Sigma$ takes the following form,
\begin{equation}
\sum_{e=1}^{\lvert \mathcal{T}\rvert} \mathcal{R}^t_e 
\left( \mathbb{B}^e\, \boldsymbol{\mathsf{W}}^e \, \textcolor{black}{+} \, \mathbb{L}^e \, \mathcal{R}_e \boldsymbol{\Lambda}\right)
\, = \, \sum_{e=1}^{\lvert \mathcal{T}\rvert} \mathcal{R}^t_e\,  \boldsymbol{\frak{s}}^e ,
\end{equation}
with coefficient matrices and sources,
\begin{equation}
 \mathbb{B}^e = \begin{pmatrix} \textcolor{black}{-}\mathbb{B}^{(e,1)}_{\displacement}  & \mathbb{B}^{(e,1)}_{\tensS} \\
\vdots  & \vdots\\[0.3em]
\textcolor{black}{-}\mathbb{B}^{(e,\nfe)}_{{\displacement}}  & \mathbb{B}^{(e,\nfe)}_{\tensS} \end{pmatrix}
, \hspace*{0.3cm} \mathbb{L}^e = \begin{pmatrix} \mathbb{L}^{(e,1)}  & \boldsymbol{0} & \boldsymbol{0}& \boldsymbol{0} \\
 \boldsymbol{0} & \mathbb{L}^{(e,2)} & \boldsymbol{0} & \boldsymbol{0}\\
  \boldsymbol{0} &   \boldsymbol{0} & \ddots & \boldsymbol{0}\\
  \boldsymbol{0} & \boldsymbol{0} & \boldsymbol{0} &\mathbb{L}^{(e,\nfe)} \end{pmatrix} ,\hspace*{0.3cm}
  \boldsymbol{\frak{s}}^e = \begin{pmatrix}\boldsymbol{\frak{s}}^{(e,1)} \\ \vdots \\ \boldsymbol{\frak{s}}^{(e,\nfe)} \end{pmatrix} ,
\end{equation}
and the local-to-global map  
$\mathcal{R}_e$ which gives the restriction of \flo{a} global vector to \flo{an} element $K^e$, 
\setlength\arraycolsep{5pt}
\begin{equation}\label{matred}
\mathcal{R}_e  = \begin{psmallmatrix}
\left[\mathcal{R}_e\right]_{(e,1),1} &  \left[\mathcal{R}_e\right]_{(e,1),2 } & \ldots &   \left[\mathcal{R}_e\right]_{(e,1), \lvert \Sigma_h\rvert}       \\[0.4em]
\left[\mathcal{R}_e\right]_{(e,2),1}    &  \left[\mathcal{R}_e\right]_{(e,2),2 }  &  \ldots   & \left[\mathcal{R}_e\right]_{(e,2),\lvert \Sigma_h\rvert}     \\[0.4em]
\vdots   &   \vdots  & \vdots   &   \vdots  \\[0.4em]
[\mathcal{R}_e]_{(e,\nfe),1}   & \left[\mathcal{R}_e\right]_{(e,\nfe),2 }& \ldots& \left[\mathcal{R}_e\right]_{(e,\nfe),\lvert \Sigma_h\rvert}
\end{psmallmatrix}.
\end{equation}
 The components of their sub-blocks are listed below.

\paragraph{Descriptions of components} The blocks of the local-to-global map
operator are given by,  for  $1\leq \ell \leq \nfe$, 
and $k$ with $1\leq k \leq \lvert \Sigma\rvert$,
\begin{equation}
[\mathcal{R}_e]_{(e,\ell), k}
\, = \, \begin{cases} \displaystyle
\mathbb{I}_{3\dofedge_k \times 3\dofedge_k} \,     , &       \text{if } \,\,\overset{\circ}{\mathsf{F}^k} \cap \partial K^e \neq \emptyset \,\, \text{ with }\,\,  \mathsf{F}^k = \mathsf{F}^{(e,\ell)} \,,  \\[0.4em]
 \displaystyle\boldsymbol{0}_{3\dofedge_{(e,\ell)} \times 3\dofedge_k }  \, , & \text{if }  \,\,\overset{\circ}{\mathsf{F}^k} \cap \,\partial K^e = \emptyset \,.
 \end{cases} 
 \end{equation}
\flo{The} vector $\boldsymbol{\frak{s}}^e$, and matrices $\mathbb{B}^e$ and $\mathbb{L}^e$ have $\nfe$ row-blocks, 
labeled by $(e,\ell)$ and corresponding to faces of $K^e$. \flo{They have the following block structures}:
\begin{align*}
\mathbb{B}_{{\displacement}}^{(e,\ell)} = \left(~_{\displacement}\mathbb{B}^{(e,\ell)}_{IJ}\right)_{I,J=x,y,z},  \hspace*{0.5cm} &\mathbb{B}_{\tensS}^{(e,\ell)} = \left(\hspace*{-0.2em}~_{\tensS}\mathbb{B}^{(e,\ell)}_{I\frak{J}} \right)_{I=x,y,z\,, \,\frak{J}\in\mathcal{I}_{\mathrm{sm}}} ,\\
 \mathbb{L}^{(e,\ell)} = \left(\mathbb{L}^{(e,\ell)}_{IJ}\right) _{I,J=x,y,z} , \hspace*{0.5cm} &\boldsymbol{\frak{s}}^{(e,\ell)}
  =  (\boldsymbol{\frak{s}}^{(e,\ell)}_I )_{I=x,y,z}.
\end{align*}
The $j$-th row of block-row $(e,\ell)$
is obtained\footnote{\flo{We have employed the identity:} 
for two vectors $\boldsymbol{\nu}, \boldsymbol{w}$ and a symmetric matrix $\boldsymbol{\sigma}\in S_2$, 
\begin{equation}\label{nu_v_tau_voigt}
\boldsymbol{\sigma}  \boldsymbol{\nu} 
 = \mathbb{A}^\dagger(\boldsymbol{\nu}) \overrightarrow{\boldsymbol{\sigma} }
\hspace*{0.2cm}  ; \hspace*{0.2cm} \boldsymbol{\sigma}  \boldsymbol{\nu}  \cdot \boldsymbol{w}
  = \boldsymbol{w} \cdot \boldsymbol{\sigma}  \cdot \boldsymbol{v} 
   =  \boldsymbol{w} \cdot \mathbb{A}^\dagger(\boldsymbol{v})\cdot  \overrightarrow{\boldsymbol{\sigma} } \,.
\end{equation}} 
by integrating corresponding boundary conditions 
\cref{intface_neu_HDGprob,abc_HDGprob} with respect
to test functions \cref{testfcnedge},
\begin{equation}\label{testfcnedge}
\boldsymbol{\xi}\big|_{F^{(e,\ell)}}  = \xi^{(e,\ell)}_j\,  \hat{\mathbf{e}}_I  \,, \quad \text{with} \,\, I = x,y,z\,. 
\end{equation}
For $  1\leq i\leq \dofedge_k$, $1\leq j\leq \dofpe$,
\begin{equation}\label{matBu} 
\left[\hspace*{-0.2em} ~_{\displacement}\mathbb{B}^{(e,\ell)}_{IJ} \right]_{ij} = \begin{dcases}
  \int_{\mathsf{F}^{(e,\ell)}}  \textcolor{black}{\tau_{IJ}} \,  \overline{\xi^{(e,\ell)}_i}\,  \phi^e_j    \,  \mathrm{d}s_\mathbf{x}  \,, \hspace*{0.1cm} &\hspace*{0.1cm} \mathsf{F}^{(e,\ell)} \in \Sigma_\mathrm{int} \cup \Sigma_\mathrm{N} \cup \Sigma_\infty,\\[-0.4em]
 \hspace*{1.5cm} 0 \,,\hspace*{0.2cm} & \hspace*{0.2cm} \mathsf{F}^{(e,\ell)} \in \Sigma_\mathrm{D} . 
  \end{dcases}  
  \end{equation}
For $  1\leq i,j\leq \dofedge_k$,
\begin{equation}\label{matbbL_u}
\left[\mathbb{L}^{(e,\ell)}_{IJ} \right]_{ij} \hspace*{-0.1em}= \hspace*{-0.1em}
\begin{dcases} \int_{\mathsf{F}^{(e,\ell)}} \textcolor{black}{\tau_{IJ}}\,  \overline{\xi^{(e,\ell)}_i}\,  \xi^{(e,\ell)}_j    \,  \mathrm{d}s_\mathbf{x}, &\hspace*{-0.2cm} \mathsf{F}^{(e,\ell)} \in \Sigma_\mathrm{int} \cup \Sigma_\mathrm{N},\\[-0.1em]
%%%
  \int_{\mathsf{F}^{(e,\ell)}} \overline{\xi^{(e,\ell)}_i} \,  \xi_j^{(e,\ell)} \left(\textcolor{black}{\tau_{IJ}} + \mathcal{Z}^\mathrm{abc}_{IJ} \right) \mathrm{d} s_\mathbf{x} 
 \,,  &\hspace*{-0.2cm} \mathsf{F}^{(e,\ell)} \in \Sigma_\infty, \\[-0.3em]
%%
%%%
\int_{\mathsf{F}^{(e,\ell)}} \overline{\xi^{(e,\ell)}_i} \,\xi^{(e,\ell)}_j  \,  \mathrm{d}s_\mathbf{x}    \,, & \hspace*{0.cm} \mathsf{F}^{(e,\ell)} \in \Sigma_\mathrm{D} .
\end{dcases}
 \end{equation}
For  $1\leq i\leq \dofedge_k$,  $1\leq j\leq \dofve$, 
\begin{equation}\label{sigmatB}
\left[\hspace*{-0.4em}~_{\tensS}\mathbb{B}^{(e,\ell)}_{I\mathfrak{J}} \right]_{ij}= \begin{dcases}
\int_{\mathsf{F}^{(e,\ell)}}   \overline{\xi^{(e,\ell)}_i} \, \textcolor{black}{\psi^e_j} \,
 \, \textcolor{black}{\mathbb{A}^\dagger(\boldsymbol{\nu}^{(e,\ell)})_{I\frak{J}}}\,  \mathrm{d}s_{\mathbf{x}},  &\hspace*{-0.1cm} \mathsf{F}^{(e,\ell)} \in \Sigma_\mathrm{int} \cup \Sigma_\mathrm{N} \cup \Sigma_\infty,\\[-0.4em]
\hspace*{2cm} 0  \,, & \hspace*{0.1cm} \mathsf{F}^{(e,\ell)} \in \Sigma_\mathrm{D}.
 \end{dcases} 
\end{equation}
We recall the Dirichlet boundary vector source $\boldsymbol{\mathsf{g}}^D = ( \mathsf{g}^\mathrm{D}_I)_{I=x,y,z}$. The components of the row block for the source $\boldsymbol{\frak{s}}^e$ are, for $1\leq i \leq \dofedge_k$, 
\begin{equation}
[\boldsymbol{\frak{s}}^{(e,\ell)}_I]_i =  \begin{dcases}
  \hspace*{1cm} 0  &,\hspace*{0.1cm} \mathsf{F}^{(e,\ell)} \in \Sigma_\mathrm{int} \cup \Sigma_\mathrm{N} \cup \Sigma_\infty,\\[-0.4em]
 \int_{\mathsf{F}^{(e,\ell)}} \overline{\mathsf{g}^{\mathrm{D}}_I} \,\xi^{(e,\ell)}_j  \,  \mathrm{d}s_\mathbf{x} \,,\hspace*{0.2cm} &, \hspace*{0.2cm} \mathsf{F}^{(e,\ell)} \in \Sigma_\mathrm{D} . 
  \end{dcases}  
\end{equation}
% =============================================
\subsubsection{Summary of discrete problem in HDG method}
% =============================================

The discretization of \cref{app_local_usig_Hdg}--\cref{homDirc_cond} 
takes the following form 
with discrete unknowns $(\boldsymbol{\mathsf{W}},\boldsymbol{\Lambda})$, 
\begin{subequations}\label{dis_forwardHDG_v0}
\begin{empheq}[left={ \empheqlbrace\,}]{align}
&\hspace*{0.7cm}  \mathbb{A}^e \, \boldsymbol{\mathsf{W}}^e 
  \,\, + \,\,\mathbb{D}^e\, \mathcal{R}_e \, \boldsymbol{\Lambda}   = \, \boldsymbol{\mathsf{S}}^e \,\,, 
% \hspace*{1.5cm}   \substack{\text{Local problem} \\ \text{ on } K^e}  
\hspace*{1.cm} \forall \,\,e =1\, \ldots, \lvert \mathcal{T}_h\rvert \,, \label{forward_hdg_local_prob} \\[0.em]
&  \displaystyle\sum_{e=1}^{\lvert \mathcal{T}\rvert}\,\, \mathcal{R}^t_e\,   \Big(\,    \mathbb{B}^e \, \boldsymbol{\mathsf{W}}^e
\, \, + \,\,  \mathbb{L}^e \, \mathcal{R}_e \, \boldsymbol{\Lambda} \, \Big) \,  = \,  \displaystyle\sum_{e=1}^{\lvert \mathcal{T}\rvert}\,\, \mathcal{R}^t_e\, \boldsymbol{\frak{s}}^e. \label{forward_hdg_globalprob}
\end{empheq}
\end{subequations}
Problem \cref{dis_forwardHDG_v0} can be reduced to 
one in terms of $\boldsymbol{\Lambda}$ only, called 
the \emph{global problem},
\begin{subequations}\label{globalprob_coefficients}
\begin{align}
 \mathbb{K} \,  \boldsymbol{\Lambda}\, = \, \mathcal{S}, \hspace*{0.5cm}
 & \text{with}\hspace*{1.5cm} \mathcal{S}:= \sum_{e=1}^{\lvert \mathcal{T}\rvert} \mathcal{R}^t_e 
\left(   \boldsymbol{\frak{s}}^e -    \mathbb{B}_e \, \mathbb{A}_e^{-1} \, \boldsymbol{\mathsf{S}}^e\right)   ,\label{global_rhs}\\[-0.5em]
&\hspace*{0cm}\text{and } \hspace*{0.2cm} \mathbb{K}:= \sum_{e=1}^{\lvert \mathcal{T}_h\rvert}\,\mathcal{R}^t_e\,\mathbb{K}^e\,  \mathcal{R}_e
\,, \hspace*{0.3cm} \text{where} \hspace*{0.2cm}\mathbb{K}^e:=  -   \mathbb{B}_e 
\left( \mathbb{A}^e\right)^{-1}\mathbb{D}^e +   \mathbb{L}_e . \label{global_coeffmat}
\end{align}
\end{subequations}

This means that in the HDG method the problem is solved in two stages.
Firstly, one solves the global problem \cref{forward_global_prob} 
in terms of the trace $\boldsymbol{\Lambda}$; secondly, the values of
the volume unknowns $\boldsymbol{\mathsf{W}}^e$ are retrieved 
element-by-element with the right-hand side of \cref{forward_local_prob},
\begin{subequations}\label{dis_forwardHDG}
\begin{empheq}[left={ \cref{dis_forwardHDG_v0} \quad \Leftrightarrow \quad \empheqlbrace\,}]{align}
& \hspace*{1cm} \mathbb{K} \,  \boldsymbol{\Lambda}\, = \, \mathcal{S} , \label{forward_global_prob}\\
& \boldsymbol{\mathsf{W}}^e \, =\,   \left(\mathbb{A}^e\right)^{-1} \left( \modif{-} \mathbb{D}^e\, \mathcal{R}_e \, \boldsymbol{\Lambda}   \, + \,      \boldsymbol{\mathsf{S}}^e \right) \,, \quad \forall e =1, \ldots, \lvert \mathcal{T}_h\rvert. \label{forward_local_prob}
\end{empheq}
\end{subequations}

% ==============================================================================
\subsection{Stabilization matrix $\tauU$}
\label{subsection:stabilization}
% ==============================================================================

The stabilization matrix $\tauU$ appears in the definition of the 
numerical trace of the traction \cref{numtraction_u} and its 
subsequent discretization, e.g., \cref{matBu,matbbL_u}.
A judicial choice of $\tauU$ is necessary to provide accurate results. 
In \flo{its} most general form, $\tauU$ is a 
symmetric \modif{positive definite} matrix. The most common 
choice is \flo{an} \emph{identity-based} stabilization, 
in which $\tauU$ is a scalar multiple of the identity matrix.
Also proposed in the literature for elastostatics, cf., 
\cite{soon2009hybridizable,fu2015analysis,cockburn2013superconvergent} is
the \emph{Kelvin--Christoffel (KC) stabilization} which corresponds 
to a scalar scaling multiple of the $3\times3$ Kelvin--Christoffel matrix 
$\boldsymbol{\Gamma}$ \cref{KCmat_def}. 

In addition to these two families of stabilization, we construct
in \cref{godunov_flux::sec} the \emph{hybridized Godunov stabilization} 
defined in terms of the Godunov matrix $M_\mathrm{Godunov}$. \modif{Its definition
for a mesh cell $K$ with outward-pointing unit normal vector $\n$ is},
\begin{equation}\label{HDGGodflux}
 M_\mathrm{Godunov} (\boldsymbol{\nu})\,=\, 
 \begin{dcases} M_\mathrm{Giso}(\boldsymbol{\nu})\,, 
 & \text{isotropy}; \\
 M_\mathrm{Gani}(\boldsymbol{\nu})\,, &
 \text{anisotropy with 3 distinct speeds } \rho \cqp^2 > \rho \cqsa^2 > \rho \cqsb^2,
 \end{dcases}
\end{equation}
\begin{equation}
\text{with} \hspace*{0.6cm} M_\mathrm{Giso}(\boldsymbol{\nu})
 = \dfrac{\rho}{\cp + \cs}\left( \cp \cs \identity + \dfrac{\boldsymbol{\Gamma}_\mathrm{iso}(\boldsymbol{\nu})}{\rho}\right)
 = \rho \left( \cs \, \identity + (\cp - \cs) \boldsymbol{\nu}\otimes \boldsymbol{\nu} \right), \label{eq:MGodunov-iso}
\end{equation}
\begin{equation}
\begin{aligned}
 & \text{and} \hspace*{0.3cm}M_\mathrm{Gani} (\boldsymbol{\nu}) := \rho \left( \cqsa + \cqsb + \cqp \right)  
           \left( \identity + \gamma \left(\dfrac{\boldsymbol{\Gamma}(\boldsymbol{\nu})}{\rho} 
           +\frak{p}_2\,  \identity\right)^{-1}  \right),\\
& \text{where} \hspace*{0.5cm} \frak{p}_2 :=  \cqsa\, \cqsb
+ \cqsa \, \cqp + \cqsb\, \cqp\,, 
  \hspace*{0.2cm} \text{and} \hspace*{0.2cm} \gamma := \dfrac{\mathrm{c}_{q\mathrm{S}1}\, \mathrm{c}_{q\mathrm{S}2} \, \mathrm{c}_\mathrm{qP}}{\mathrm{c}_{qS1} + \mathrm{c}_{qS2} + \mathrm{c}_{q\mathrm{P}} } - \frak{p}_2.
\end{aligned}
\end{equation} 
\modif{Here, the density $\rho$ and wavespeeds are associated with mesh cell $K$}. 
Hybridization for isotropy was carried out in \cite{terrana2018spectral}, 
we also refer to \cite{Pham2023hdgRR} \flo{which shows the derivation}
of $M_\mathrm{Giso}$ in terms of $\boldsymbol{\Gamma}$. 
The derivation for $ M_\mathrm{Gani}$ is given in \cref{godunov_flux::sec}.
In \cref{section:numerics}, we  will investigate the three families
(identity-based, KC-based and Godunov-based) and show the optimality of the 
Godunov stabilization.

\section{Construction of hybridized Godunov stabilization operators}\label{godunov_flux::sec}

In this section, we  extend the hybridization method in \cite{terrana2018spectral} 
for isotropic elasticity to anisotropy, and construct the hybridized Godunov 
stabilization operator $ M_\mathrm{Godunov} (\boldsymbol{\nu})$ \cref{HDGGodflux}.
We recall that \flo{a} stabilization matrix is employed to define the numerical trace of traction $\widehat{\boldsymbol{\sigma}\boldsymbol{\nu}}$
on each mesh element.
We consider \flo{an} elastic material whose KC matrix $\boldsymbol{\Gamma}(\boldsymbol{\nu})$ \cref{KCmat_def} has three distinct eigenvalues,
\begin{equation}\label{3s::asump}
  \rho \cqp^2 \,\,>\,\, \rho \cqsa^2 \,\,>\,\, \rho \cqsb^2 .
\end{equation}
The discussion in this section employs the Voigt notation introduced in  \cref{VI::subsec}--
\cref{subsection:KVmatrix}.

In \cref{flux::subsec}, we introduce a first-order formulation of the time-dependent elastic equation \cref{1sys_esp_v} 
with unknown $\mathbf{q}=(\boldsymbol{\epsilon},\rho\mathbf{v})$, and make appear the 
flux term $\frak{B}(\boldsymbol{\nu})\mathbf{q}$ \cref{fluxop::def} whose last three components give  
$\boldsymbol{\sigma}\boldsymbol{\nu}$. 
This means that the numerical trace of traction $\widehat{\boldsymbol{\sigma}\boldsymbol{\nu}}$ at a face $\mathsf{F}$
 can be obtained from the numerical trace $\widehat{\frak{B}(\boldsymbol{\nu})\mathbf{q}}$.
Additionally, we consider $\mathbf{q}$ as a collection of values on each mesh element, i.e,
$\mathbf{q} = \cup_{K\in \mathcal{T}_h} \mathbf{q}^K$. 
At an interface $\mathsf{F}$ shared by two elements $K^+$ and $K^-$, 
% \modifspace{we denote the trace of $\mathbf{q}$ from the left and right 
% of $\mathsf{F}$ as $ \mathbf{q}^{\pm}$, i.e.}
\flo{we denote by $\mathbf{q}^{\pm}$ the trace of $\mathbf{q}$ 
from each side of $\mathsf{F}$ associated with the neighbor 
cells $K^{\pm}$, i.e.}
\begin{equation}\label{interFlr}
\modifspace{\mathsf{F} = \partial K^-\cap \partial K^+,  \hspace*{0.2cm}   \mathsf{F}^{\pm} := \partial K^{\pm} \cap \mathsf{F}.}
\end{equation}
For the construction of $\widehat{\frak{B}(\boldsymbol{\nu})\mathbf{q}}$, we
work with \flo{the} exact solutions of the Riemann problem in a neighborhood of $\mathsf{F}$,  
\flo{having as} initial data $\mathbf{q}^{\pm}$, denoted by $\mathrm{RP}(\mathbf{q}^-, \mathbf{q}^+)$.
The solution of $\mathrm{RP}(\mathbf{q}^-, \mathbf{q}^+)$ to the immediate left and right of $\mathsf{F}$
is called intermediate states $\mathbf{q}^\star$.

\medskip
The main construction is given in \cref{Godunov_ani::subsec} and consists of 3 main steps. 
\begin{itemize}[leftmargin=*]
\item The numerical flux $\widehat{\frak{B}(\boldsymbol{\nu})\mathbf{q}}$ is given by the Godunov flux, which by definition is the value of $\frak{B}(\boldsymbol{\nu})$ at \flo{the}
intermediate state $\mathbf{q}^\star$, i.e. $ \frak{B}(\boldsymbol{\nu})\mathbf{q}^\star$.

\item The relations between the flux of the intermediate states $\mathbf{q}^\star$ and the left 
and right data $\mathbf{q}^{\pm}$ are given by a system of Rankine--Hugoniot (RH) jump 
conditions, cf. \cref{ran_hug_jump_ani_strain}. 
 From this, 
 we derive the usual transmission conditions for elasticity \cref{inte_eq_ani}, and an 
 expression of the Godunov flux, cf.
\cref{rel_Aqstar_qm_0}, and thus of the traction, cf. \cref{star_minus_traction::ani}, 
\flo{which is} in terms of one-sided data (i.e., either $\mathbf{q}^+$ or $\mathbf{q}^-)$.

\item Finally, in a step called `hybridization', these above ingredients are combined
to obtain the hybridized HDG numerical flux \cref{hybrid_ani}.
\end{itemize}
%
%%The first two ideas are employed in the construction of 
%%upwind fluxes for DG implementation in time domain, 
%%cf. \cite{wilcox2010high,tie2018unified,tie2020systematic,zhan2018exact}, 
%%
%%For applications with DG, the constructed numerical flux depends 
%%on information from both left and right of the interface. 
%%For implementation with HDG, also starting from these 
%%ingredients, \cite{terrana2018spectral} however proceeds differently 
%%to arrive at an expression of the flux which contains only one-sided data. 
%%
%
%
%Although constructed from the same ingredients, the DG numerical flux in \cite{wilcox2010high,tie2018unified,tie2020systematic,zhan2018exact} depends 
%on information both from left and right of the interface, 
%while the hybridized HDG numerical fluxes in \cite{terrana2018spectral} 
%and in our work are defined from one-sided data. 

For the hybridization in isotropy, we refer to \cite{terrana2018spectral}
and \cite{Pham2023hdgRR}; in the latter, the derivation is written with Voigt notation
which brings out the connection \flo{between} the stabilization operator 
\flo{and the} Kelvin--Christoffel matrix.
Details for \flo{the} identities employed in this section with Voigt notation 
as well as \flo{a} review of the Riemann problem can be found in \cite{Pham2023hdgRR}.

\begin{remark}
We adopt the name `Godunov-upwind' from \cite{bui2015godunov,bui2015rankine,bui2016construction}, 
in which the hybridization was carried out following different approaches; however, employing Rankine--Hugoniot jump condition is recognized in sequel work \cite{bui2015rankine,bui2016construction} 
(where the hybridized flux is also called `Rankine--Hugoniot' flux) to be more natural and direct.
\end{remark}

% ===================================================================
\subsection{First order system and flux operator}
\label{flux::subsec}
% ===================================================================

We rewrite the time-dependent elastic equation as a first order system with 
unknowns $(\boldsymbol{\epsilon},\rho \mathbf{v})$, 
\modif{in recalling $\tensS = \tensC \boldsymbol{\epsilon}$, $\boldsymbol{\epsilon} =\nabla^{\mathrm{s}}\displacement$, 
and $ \velocity =  \partial_t\displacement$}:
\begin{equation}\label{1sys_esp_v}
\begin{dcases}
\partial_t \boldsymbol{\epsilon}  - \nabla^\mathfrak{s} \mathbf{v} \, =\, 0 \,, \\[0.3em]
\partial_t  \rho \mathbf{v}   - \nabla \cdot 
\modif{(} \tensC \boldsymbol{\epsilon} \modif{)}\, = \, \boldsymbol{0}
\end{dcases}\hspace*{0.3cm} 
\Leftrightarrow \hspace*{0.3cm}  \partial_t \begin{pmatrix} \boldsymbol{\epsilon} \\
\rho\mathbf{v} \end{pmatrix} = \begin{pmatrix} 0 & \nabla^\mathfrak{s}\, \rho^{-1}\\  \nabla \cdot \,\tensC & 0 \end{pmatrix}\begin{pmatrix} \boldsymbol{\epsilon} \\
\rho\mathbf{v} \end{pmatrix} \,.
\end{equation} 
With matrix $\mathbb{A}_I$ \cref{bbAmat::def}, define matrix-valued differential operator $ \mathbb{A}(\partial_\mathbf{x})$ and its transpose,
\begin{equation}
 \mathbb{A}(\partial_\mathbf{x}) := \sum_{I=x,y,z} \partial_I\, \mathbb{A}_I \,, \hspace*{0.7cm}
% ------------------------------------------------------------------- % %
\mathbb{A}(\partial_\mathbf{x})^t = \sum_I \partial_I \mathbb{A}_I^t \,.
\end{equation}
Employing \flo{the following} identity, cf. \cite{Pham2023hdgRR},
\begin{equation}
\overrightarrow{\nabla^\mathfrak{s} \mathbf{v}} =  \mathbb{A}(\partial_\mathbf{x})^t \mathbf{v}, \hspace*{0.5cm}
\nabla\cdot ( \tensC \boldsymbol{\epsilon}) =  \mathbb{A}(\partial_\mathbf{x}) \, ^\dagger\!\overset{=}{\tensC}\,\! ^\dagger   \overrightarrow{\boldsymbol{\epsilon}},
\end{equation}
equation \cref{1sys_esp_v} can be written as
\begin{equation}
\partial_t \mathbf{q} + \frak{B}(\partial\mathbf{x})\, 
\mathbf{q} =  \begin{pmatrix} 0 \\ \mathbf{f} \end{pmatrix}, \end{equation}
\begin{equation}\text{with} \hspace*{0.5cm}
\mathfrak{B}(\partial_\mathbf{x}) :=- \begin{pmatrix}  0 &  \textcolor{black}{\mathbb{A}^t(\partial_\mathbf{x})} \,\rho^{-1} \\
\textcolor{black}{\mathbb{A}(\partial_\mathbf{x})  ^\dagger\!\overset{=}{\tensC}\,\! ^\dagger } & 0 \end{pmatrix} \hspace*{0.1cm}
\text{and }\hspace*{0.1cm} \mathbf{q} := \begin{pmatrix} \overrightarrow{\boldsymbol{\epsilon}} \\ \rho \mathbf{v}\end{pmatrix}\,.
\end{equation}

\paragraph{Flux term} 
Along $\partial K$ which has outward pointing normal vector $\boldsymbol{\nu}$, we define the matrix
\begin{equation}\label{fluxop::def}
\mathfrak{B}(\boldsymbol{\nu}) :=- \begin{pmatrix}  0 &  \textcolor{black}{\mathbb{A}^t(\boldsymbol{\nu})} \rho^{-1} \\
\textcolor{black}{\mathbb{A}(\boldsymbol{\nu}) \, ^\dagger\!\overset{=}{\tensC}\,\! ^\dagger } & 0 \end{pmatrix}.
\end{equation}
This is called the flux matrix at $\partial K$ and arises from integrating both sides of \cref{1sys_esp_v} 
against test function $\begin{pmatrix} \boldsymbol{\chi} \\ \mathbf{w}\end{pmatrix}$ 
with $\boldsymbol{\chi}\in \mathcal{S}_2$ and vector $\mathbf{w} \in \mathbb{C}^3$,
and carrying out an integration by parts,
 \begin{equation}\label{IP_fluxop}
 \begin{aligned}
&\int_K \partial_t \mathbf{q} \cdot  \begin{pmatrix}\overrightarrow{\boldsymbol{\chi}}^\dagger \\ \mathbf{w}\end{pmatrix} \, \mathrm{d}\mathbf{x} \hspace*{0.1cm} + \hspace*{0.1cm}
 \int_K \mathfrak{B}(\partial_\mathbf{x}) \begin{pmatrix}\overrightarrow{\boldsymbol{\epsilon}} \\ \rho \mathbf{v}\end{pmatrix}
 \cdot \begin{pmatrix}\overrightarrow{\boldsymbol{\chi}}^\dagger \\ \mathbf{w}\end{pmatrix}\, \mathrm{d}\mathbf{x} \\[0.3em]
&\hspace*{0.1cm}  = \int_K \partial_t \mathbf{q} \cdot  \begin{pmatrix}\overrightarrow{\boldsymbol{\chi}}^\dagger \\ \mathbf{w}\end{pmatrix} \, \mathrm{d}\mathbf{x} +\int_K  \begin{pmatrix}\overrightarrow{\boldsymbol{\epsilon}} \\ \rho \mathbf{v}\end{pmatrix} \cdot \frak{B}^\mathrm{ad}(\partial_\mathbf{x}) \begin{pmatrix}\overrightarrow{\boldsymbol{\chi}}^\dagger \\ \mathbf{w}\end{pmatrix} \, \mathrm{d}\mathbf{x}
   \hspace*{0.1cm} + \hspace*{0.1cm} \int_{\partial K} \mathfrak{B}(\boldsymbol{\nu}) \begin{pmatrix}\overrightarrow{\boldsymbol{\epsilon}} \\ \rho \mathbf{v}\end{pmatrix}
 \cdot \begin{pmatrix}\overrightarrow{\boldsymbol{\chi}}^\dagger \\ \mathbf{w}\end{pmatrix} \, \mathrm{d}s_\mathbf{x}.
\end{aligned} 
 \end{equation}
Here `adjoint' operator associated with $\frak{B}$ is
\begin{equation}\label{adfrakB::def}
\frak{B}^\mathrm{ad}(\partial_\mathbf{x}) :=  \begin{pmatrix}  0 &  \textcolor{black}{^\dagger\!\overset{=}{\tensC}\,\! ^\dagger  \, \mathbb{A}(\partial_\mathbf{x})^t} \,   \\
\rho^{-1}\,\textcolor{black}{\mathbb{A}(\partial_\mathbf{x})}
 & 0 \end{pmatrix}.
\end{equation}

\paragraph{Eigenvalues of $\frak{B}(\boldsymbol{\nu})$} 
Consider the eigenvalue problem 
of the flux operator $\mathfrak{B}(\boldsymbol{\nu})$, 
\begin{equation}
\begin{pmatrix}  0 &  \textcolor{black}{\mathbb{A}^t(\boldsymbol{\nu})} \rho^{-1} \\
\textcolor{black}{\mathbb{A}(\boldsymbol{\nu}) \, ^\dagger\!\overset{=}{\tensC}\,\! ^\dagger } & 0 \end{pmatrix} \begin{pmatrix} 
\overset{\rightarrow}{\boldsymbol{\tau} }\\ \mathbf{w} \end{pmatrix}
 = \alpha \begin{pmatrix} 
\overset{\rightarrow}{\boldsymbol{\tau} }\\ \mathbf{w} \end{pmatrix}
\hspace*{0.5cm} \Leftrightarrow \hspace*{0.5cm} \begin{cases}
\textcolor{black}{\mathbb{A}^t(\boldsymbol{\nu})} \rho^{-1}\mathbf{w} \,\,=\,\, \textcolor{black}{\alpha}\, \overset{\rightarrow}{\boldsymbol{\tau} } \\
\textcolor{black}{\mathbb{A}(\boldsymbol{\nu}) \, ^\dagger\!\overset{=}{\tensC}\,\! ^\dagger }  \overset{\rightarrow}{\boldsymbol{\tau} } \,\,= \,\,\textcolor{black}{\alpha}\, \mathbf{w}
\end{cases} \,.
\end{equation}
Applying $\textcolor{black}{\mathbb{A}(\boldsymbol{\nu}) \, ^\dagger\!\overset{=}{\tensC}\,\! ^\dagger }$
to both sides of the first equation, we obtain the eigenproblem for \flo{the} Kelvin--Christoffel matrix 
$\boldsymbol{\Gamma} (\boldsymbol{\nu})$, which is also called the Christoffel equation,
\begin{equation}
\textcolor{black}{\mathbb{A}(\boldsymbol{\nu}) \, ^\dagger\!\overset{=}{\tensC}\,\! ^\dagger } \textcolor{black}{\mathbb{A}^t(\boldsymbol{\nu})} \, \mathbf{w}\, =\, \rho\, \textcolor{black}{\alpha}^2 \,  \mathbf{w} 
\quad \Leftrightarrow \quad
\boldsymbol{\Gamma} (\boldsymbol{\nu})\, \mathbf{w}\, =\, \rho \,\textcolor{black}{\alpha}^2 \, \mathbf{w} \,.
\end{equation}
We work under the assumption that $\boldsymbol{\Gamma}(\boldsymbol{\nu})$ has 
the three distinct eigenvalues,
\begin{equation}\label{3s::asump_rep}
\rho \,\cqp^2 \hspace*{0.2cm} <\hspace*{0.2cm} \rho\, \cqsa^2
\hspace*{0.2cm} <\hspace*{0.2cm} \rho\, \cqsb^2 \,.
\end{equation}
This implies that \flo{the} eigenvalues of $\frak{B}(\boldsymbol{\nu})$ are
\begin{equation}\label{EVfrakB_ani}
-\cqp \hspace*{0.2cm}\,, \hspace*{0.2cm} -\cqsa
\,, \hspace*{0.2cm} - \cqsb \hspace*{0.2cm}\,, \hspace*{0.2cm} \underset{\substack{\\[0.2em]\text{ multiplicity 3}}}{0}\hspace*{0.2cm}\,, \hspace*{0.2cm} \cqsb \hspace*{0.2cm}\,, \hspace*{0.2cm} \cqsa \hspace*{0.2cm}\,, \hspace*{0.2cm}\cqp \,.
\end{equation}

\begin{remark}
For elastic isotropy, \cref{iso_stiffness_voigt}, the 
eigenvalues of $\frak{B}(\boldsymbol{\nu})$ are
\begin{equation}\label{EVfrakB_iso}
-\cp \hspace*{0.2cm}\,, \hspace*{0.2cm} \underset{\substack{\\[0.2em]\text{multiplicity 2}}}{-\cs}
 \hspace*{0.2cm}\,, \hspace*{0.2cm} \underset{\substack{\\[0.2em]\text{ multiplicity 3}}}{0}\hspace*{0.2cm}\,, \hspace*{0.2cm} \underset{\substack{\\[0.2em]\text{multiplicity 2}}}{\cs}\hspace*{0.2cm}\,, \hspace*{0.2cm}\cp \,.
\end{equation} 
\end{remark}
\subsection{Derivation for anisotropic elasticity with distinct waves speeds}
\label{Godunov_ani::subsec}
% ========================================================================================
\modif{Along an interface $\mathsf{F}$ shared by mesh cells $K^+$ and $K^-$, cf. \cref{interFlr}}, 
we denote by $\boldsymbol{\nu}^{\pm}$ the outward-pointing normal vectors
 of $K^{\pm}$, and by $\boldsymbol{\nu} $ the normal vector 
\modif{pointing from} $K^-$ to $K^+$, i.e.
% \begin{equation}\label{left_right_interface_nu}
$ \boldsymbol{\nu} = \boldsymbol{\nu}^- = - \boldsymbol{\nu}^+$.
%\end{equation}
\modif{To distinguish the background parameters on each side of $\mathsf{F}$}, we write
\begin{equation}\label{left_right_interface}
\rho^{\pm},   \quad\tensC^{\pm}, \quad \modif{\boldsymbol{\Gamma}(\boldsymbol{\nu})^{\pm}},\quad \frak{B}(\boldsymbol{\nu})^{\pm}  \text{ with eigenvalues }  \mathrm{c}_{\alpha}^{\pm} \, .
\end{equation}
Here, $\modif{\frak{B}(\boldsymbol{\nu})^-}$ is defined in \cref{fluxop::def} 
with \flo{the} physical parameters \modif{of} $K^-$, i.e.,
$(\rho^-,\lambda^-,\mu^-)$ and \modif{its} normal vector $\boldsymbol{\nu}^-$; similarly for $\frak{B}(\boldsymbol{\nu})^+$.
\modif{The KC matrix associated with each mesh cell, $\boldsymbol{\Gamma}(\boldsymbol{\nu})^{\pm} $, is defined in the same manner with definition \cref{KCmat_def}.}

% ------------------------------------------------------------
\begin{figure}[h!]\centering
  \includegraphics[width=8.5cm]{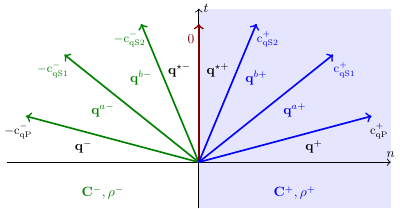}
  \caption{\flo{The eigenvalues,} spectral structure and the Godunov states appearing in the Rankine-Hugoniot jump condition \cref{ran_hug_jump_ani_strain} for anisotropic elasticity with 3 distinct speeds.}
  \label{rk_ani::fig}
\end{figure}
% ------------------------------------------------------------

Under the assumption of spectral structure \cref{EVfrakB_ani} of $\modif{\frak{B}(\boldsymbol{\nu})^{\pm}}$, we have 7 possible discontinuities propagating at the speed of corresponding eigenvalues indicated in  \cref{EVfrakB_ani}.
Denote the states in between the discontinuities \modif{by}, 
cf. \cref{rk_ani::fig},
\begin{equation}
\mathbf{q}^-  \hspace*{0.2cm}\,, \hspace*{0.2cm} \mathbf{q}^{a-} 
 \hspace*{0.2cm}\,, \hspace*{0.2cm} \mathbf{q}^{b-}   \hspace*{0.2cm}\,, \hspace*{0.2cm}\mathbf{q}^{\star-}
  \hspace*{0.2cm}\,, \hspace*{0.2cm}\mathbf{q}^{\star+}  \hspace*{0.2cm}\,, \hspace*{0.2cm}\mathbf{q}^{b+}  \hspace*{0.2cm}\,, \hspace*{0.2cm}\mathbf{q}^{a+}  \hspace*{0.2cm}\,, \hspace*{0.2cm} \mathbf{q}^+\,.
\end{equation}
They satisfy jump condition, cf. \cite[Equation (19)]{tie2020systematic}, 
\begin{subequations}\label{ran_hug_jump_ani_strain}
 \begin{empheq}[left = {  \empheqlbrace}\,]{align}
\modif{\frak{B}(\boldsymbol{\nu})^-} \left(  \mathbf{q}^{a-} - \mathbf{q}^-  \right) &= -\mathrm{c}_{q\mathrm{P}}^- \left(  \mathbf{q}^{a-} - \mathbf{q}^-  \right) \,;\label{ran_hug_jump_ani_strain_mcp}\\
\modif{\frak{B}(\boldsymbol{\nu})^-}\left(  \mathbf{q}^{b-} - \mathbf{q}^{a-}  \right) &= -\mathrm{c}_{q\mathrm{S}1}^- \left(  \mathbf{q}^{b-} - \mathbf{q}^{a-}  \right) \,;\label{ran_hug_jump_ani_strain_mcs1}\\
\modif{\frak{B}(\boldsymbol{\nu})^-}\left(  \mathbf{q}^{\star-} - \mathbf{q}^{b-}  \right) &= -\mathrm{c}_{q\mathrm{S}2}^- \left(  \mathbf{q}^{\star-} - \mathbf{q}^{b-}  \right) \,;\label{ran_hug_jump_ani_strain_mcs2}\\
\modif{\frak{B}(\boldsymbol{\nu})^-}   \mathbf{q}^{\star-} &=\modif{\frak{B}(\boldsymbol{\nu})^+}   \mathbf{q}^{\star+} \,;\label{ran_hug_jump_ani_strain_s0}\\
\modif{\frak{B}(\boldsymbol{\nu})^+} \left(  \mathbf{q}^{b+} - \mathbf{q}^{\star+}  \right) &= \mathrm{c}_{q\mathrm{S}1}^+ \left(  \mathbf{q}^{b+} - \mathbf{q}^{\star+}  \right) \,;\\
\modif{\frak{B}(\boldsymbol{\nu})^+} \left(  \mathbf{q}^{b+} - \mathbf{q}^{a+}  \right) &= \mathrm{c}_{q\mathrm{S}1}^+ \left(  \mathbf{q}^{b+} - \mathbf{q}^{a+}  \right) \,;\\
\modif{\frak{B}(\boldsymbol{\nu})^+}  \left(  \mathbf{q}^{a+} - \mathbf{q}^+  \right) &= \mathrm{c}_{q\mathrm{S}}^+ \left(  \mathbf{q}^{a+} - \mathbf{q}^+ \right) \,.
\end{empheq}
\end{subequations}

\paragraph{Transmission conditions for the interface states}
We obtain the 
 transmission conditions from the first 3 and last 3 components 
 of equation \cref{ran_hug_jump_ani_strain_s0} associated with the 
 non-propagative state $\mathbf{q}^{\star\pm}$,
\begin{equation}\label{inte_eq_ani}
  \begin{cases} \modif{\nu_I^-}\,  v_I^{\star-}\,  =\,  \modif{\nu_I^+} \, v_I^{\star+} \,, \quad I = x,y,z \\[0.3em]
  \hspace*{0.2cm} \left(\boldsymbol{\sigma}\boldsymbol{\nu} \right)^{\star-} \, = \,
  \left(\boldsymbol{\sigma} \boldsymbol{\nu}\right)^{\star+}\end{cases}
  \quad \Rightarrow \quad
  \begin{cases}  \mathbf{v}^{\star-}\, &= \, \mathbf{v}^{\star+}\,;\\
  \left(\boldsymbol{\sigma}\boldsymbol{\nu} \right)^{\star-} \, &= \,
  \left(\boldsymbol{\sigma} \boldsymbol{\nu}\right)^{\star+} \end{cases}\,.
\end{equation}
\modif{Here $\nu_I^{\pm}$ are the components of the normal vectors $\boldsymbol{\nu}^{\pm}$, i.e., $\boldsymbol{\nu}^{\pm}=\modif{(\nu_x^{\pm},\nu_y^{\pm},\nu_z^{\pm})^t}$}. 

\paragraph{\modif{Relation between the} intermediate states with the interface and boundary states} 
Working on each side of the interface $\mathsf{F}$, we will express 
the state $\mathbf{q}^{a\pm}$ and $\mathbf{q}^{b\pm}$ 
in terms of the interface states $\mathbf{q}^{\star\pm}$ 
and \flo{the} boundary states $\mathbf{q}^{\pm}$. 
\modif{To alleviate the notation in the derivation, we write}
\begin{equation}
\modif{\frak{B}_{\pm} :=\frak{B}(\boldsymbol{\nu})^{\pm} . }
\end{equation}

\noindent \textbf{Step 1a} 
We start by working with the first three equations 
of \cref{ran_hug_jump_ani_strain}. \modif{Summing them yields}, 
% three equations $ \cref{ran_hug_jump_ani_strain_mcp} +    \cref{ran_hug_jump_ani_strain_mcs1} +   \cref{ran_hug_jump_ani_strain_mcs2}$, 
\begin{equation}
\modif{\frak{B}_-} \left(  \mathbf{q}^{\star-} - \mathbf{q}^-  \right)
 = \left( -\mathrm{c}_{q\mathrm{P}}^-  + \mathrm{c}_{q\mathrm{S}1}^- \right)  \mathbf{q}^{a-}
+  \left( -\mathrm{c}_{q\mathrm{S}1}^-  + \mathrm{c}_{q\mathrm{S}2}^- \right)  \mathbf{q}^{b-} 
  -\mathrm{c}_{q\mathrm{S}2}^-  \mathbf{q}^{\star-}
   +   \mathrm{c}_{q\mathrm{P}}^- \mathbf{q}^-\,.
\end{equation}
After rearrangement, this leads to,
\begin{equation}\label{sum_ran_hug_jump_ani_strain_minus}
 \textcolor{black}{\left(\mathrm{c}_{q\mathrm{P}}^- -\mathrm{c}_{q\mathrm{S}1}^-   \right)  \mathbf{q}^{a-}
+  \left(\mathrm{c}_{q\mathrm{S}1}^-  -\mathrm{c}_{q\mathrm{S}2}^-   \right)  \mathbf{q}^{b-} }
 = -\modif{\frak{B}_-}  \left(  \mathbf{q}^{\star-} - \mathbf{q}^-  \right) -\mathrm{c}_{q\mathrm{S}2}^-  \mathbf{q}^{\star-}+   \mathrm{c}_{q\mathrm{P}}^- \mathbf{q}^-\,.
\end{equation}
Next, \modif{taking} $c_{qP}^{\modif{-}} \times \cref{ran_hug_jump_ani_strain_mcp} +  c_{qS1}^{\modif{-}}   \times\cref{ran_hug_jump_ani_strain_mcs1} + c_{qS2}^{\modif{-}}  \times \cref{ran_hug_jump_ani_strain_mcs2}$, 
we obtain,
\begin{equation}\label{wsum_ran_hug_jump_ani_strain_minus}
\begin{aligned}
&\modif{\frak{B}_-} 
\left[\textcolor{black}{\left( \mathrm{c}_{q\mathrm{P}}^-  - \mathrm{c}_{q\mathrm{S}1}^- \right)  \mathbf{q}^{a-}                   
+ \left( \mathrm{c}_{q\mathrm{S}1}^-  - \mathrm{c}_{q\mathrm{S}2}^- \right)  \mathbf{q}^{b-} }
 \right]
 = \mathrm{c}_{q\mathrm{P}}^- \modif{\frak{B}_-}  \mathbf{q}^-
 - \mathrm{c}_{q\mathrm{S2}}^- \modif{\frak{B}_-}  \mathbf{q}^{\star-} \\
&\hspace*{1cm} +   \left[ (\mathrm{c}_{q\mathrm{S}1}^-)^2-(\mathrm{c}_{q\mathrm{P}}^-)^2   \right]  \mathbf{q}^{a-} 
+  \left[ (\mathrm{c}^-_{q\mathrm{S}2})^2-(\mathrm{c}_{q\mathrm{S}1}^-)^2    \right]  \mathbf{q}^{b-} 
  -(\mathrm{c}_{q\mathrm{S}2}^-)^2  \mathbf{q}^{\star-}  +   (\mathrm{c}_{q\mathrm{P}}^-)^2 \mathbf{q}^- \,.
\end{aligned}
\end{equation}

\medskip

\noindent \textbf{Step 1b}
Using \cref{sum_ran_hug_jump_ani_strain_minus} to rewrite the left-hand side of equation \cref{wsum_ran_hug_jump_ani_strain_minus}, we obtain,
\begin{equation}
\begin{aligned}
& \modif{\frak{B}_-} \left(- \modif{\frak{B}_-}  \left(  \mathbf{q}^{\star-} - \mathbf{q}^-  \right)  -\mathrm{c}_{q\mathrm{S}2}^-  \mathbf{q}^{\star-}
   +   \mathrm{c}_{q\mathrm{P}}^- \mathbf{q}^- \right)
 \hspace*{0.1cm}  =\hspace*{0.1cm} \mathrm{c}_{q\mathrm{P}}^- \modif{\frak{B}_-}  \mathbf{q}^-
 - \mathrm{c}_{q\mathrm{S2}}^- \modif{\frak{B}_-} \mathbf{q}^{\star-} \\
&\hspace*{0.7cm} +  \left( -(\mathrm{c}_{q\mathrm{P}}^-)^2  + (\mathrm{c}_{q\mathrm{S}1}^-)^2 \right)  \mathbf{q}^{a-} 
+  \left( -(\mathrm{c}_{q\mathrm{S}1}^-)^2  + (\mathrm{c}^-_{q\mathrm{S}2})^2 \right)  \mathbf{q}^{b-} 
  -(\mathrm{c}_{q\mathrm{S}2}^-)^2  \mathbf{q}^{\star-}
   +   (\mathrm{c}_{q\mathrm{P}}^-)^2 \mathbf{q}^-\,.
\end{aligned}
\end{equation}
After simplification, we arrive at,
\begin{equation}\label{seeqn_ani_intstats}
\begin{aligned}
 - \modif{\frak{B}_-^2} \left(  \mathbf{q}^{\star-} - \mathbf{q}^-  \right)  
\hspace*{0.1cm}&=\hspace*{0.1cm} \left( -(\mathrm{c}_{q\mathrm{P}}^-)^2  + (\mathrm{c}_{q\mathrm{S}1}^-)^2 \right)  \mathbf{q}^{a-} 
+  \left( -(\mathrm{c}_{q\mathrm{S}1}^-)^2  + (\mathrm{c}^-_{q\mathrm{S}2})^2 \right)  \mathbf{q}^{b-} \\
 &\hspace*{5.5cm} -(\mathrm{c}_{q\mathrm{S}2}^-)^2  \mathbf{q}^{\star-}
   +   (\mathrm{c}_{q\mathrm{P}}^-)^2 \mathbf{q}^-\,.
\end{aligned}
\end{equation}
%\begin{equation}\label{seeqn_ani_intstats}
%\begin{aligned}
% - \modif{\frak{B}_-^2} \left(  \mathbf{q}^{\star-} - \mathbf{q}^-  \right)  
%&=\hspace*{0.cm} \left( -(\mathrm{c}_{q\mathrm{P}}^-)^2  + (\mathrm{c}_{q\mathrm{S}1}^-)^2 \right)  \mathbf{q}^{a-} 
% %%
%+  \left( -(\mathrm{c}_{q\mathrm{S}1}^-)^2  + (\mathrm{c}^-_{q\mathrm{S}2})^2 \right)  \mathbf{q}^{b-} -(\mathrm{c}_{q\mathrm{S}2}^-)^2  \mathbf{q}^{\star-}
%   +   (\mathrm{c}_{q\mathrm{P}}^-)^2 \mathbf{q}^-\,.
%\end{aligned}
%\end{equation}
\flo{Equations \cref{sum_ran_hug_jump_ani_strain_minus,seeqn_ani_intstats}}
give a linear system which determines uniquely $(\mathbf{q}^{a-},\mathbf{q}^{b-})$ \modif{in terms of} $\mathbf{q}^{\star-}$ and $\mathbf{q}^-$,
\begin{equation}\label{system_qa_qb_ani}
\begin{pmatrix}
 \alpha   &  \beta   \\
  \alpha \tilde{\alpha}  &
\beta \tilde{\beta}  
\end{pmatrix} \begin{pmatrix}
\mathbf{q}^{a-}  \\ \mathbf{q}^{b-} 
\end{pmatrix} = 
\begin{pmatrix}
-\modif{\frak{B}_-}  (\mathbf{q}^{\star-} - \mathbf{q}^-)  -\mathrm{c}_{q\mathrm{S}2}^-  \mathbf{q}^{\star-}
   +   \mathrm{c}_{q\mathrm{P}}^- \mathbf{q}^- \\[0.3em]
\modif{\frak{B}_-^2}(\mathbf{q}^{\star-} - \mathbf{q}^-) -(\mathrm{c}_{q\mathrm{S}2}^-)^2  \mathbf{q}^{\star-}
   +   (\mathrm{c}_{q\mathrm{P}}^-)^2 \mathbf{q}^-
\end{pmatrix} ,
\end{equation}
\begin{equation}
\text{where} \hspace*{0.3cm} \alpha = \mathrm{c}_{q\mathrm{P}}^- -\mathrm{c}_{q\mathrm{S}1}^-
\,, \quad \tilde{\alpha} =  \mathrm{c}_{q\mathrm{S}1}^-+\mathrm{c}_{q\mathrm{P}}^- 
\,, \quad \beta = \mathrm{c}_{q\mathrm{S1}}^- -\mathrm{c}_{q\mathrm{S}2}^-
\, \quad \tilde{\beta} = \mathrm{c}_{q\mathrm{S}1}^-+\mathrm{c}_{q\mathrm{S}2}^- \,.
\end{equation}
Note that the determinant of the coefficient matrix in \cref{system_qa_qb_ani} is
\begin{equation}
\alpha \beta ( \tilde{\beta} - \tilde{\alpha}) = 
\left(\mathrm{c}_{q\mathrm{P}}^- - \mathrm{c}_{q\mathrm{S}1}^-\right) \left(\mathrm{c}_{q\mathrm{S}1}^- - \mathrm{c}_{q\mathrm{S}2}^- \right)
\left(\mathrm{c}_{q\mathrm{S}2}^- - \mathrm{c}_{q\mathrm{P}}^- \right)\,.
 \end{equation}
 \modif{This matrix is thus} invertible under \modif{the} assumption of distinct wave speeds \cref{3s::asump_rep}. In this case, we obtain the expression of $\mathbf{q}^{a-}$ and $\mathbf{q}^{b-}$ in terms of $\mathbf{q}^{\star-}$
 and $\mathbf{q}^-$,
 \begin{equation}\label{qaqb}
 \begin{pmatrix}
 \mathbf{q}^{a-}\\[0.3em] \mathbf{q}^{b-}
 \end{pmatrix}
 = \dfrac{1}{\tilde{\beta} - \tilde{\alpha} } \begin{pmatrix}
\tfrac{ \tilde{\beta}}{\alpha}     &   - \tfrac{1}{\alpha}   \\[0.5em]
- \tfrac{\tilde{\alpha}}{\beta}  & \tfrac{1}{\beta}
\end{pmatrix}\begin{pmatrix}
-\modif{\frak{B}_-}  (\mathbf{q}^{\star-} - \mathbf{q}^-)  -\mathrm{c}_{q\mathrm{S}2}^-  \mathbf{q}^{\star-}
   +   \mathrm{c}_{q\mathrm{P}}^- \mathbf{q}^- \\[0.5em]
\modif{\frak{B}_-^2} (\mathbf{q}^{\star-} - \mathbf{q}^-) -(\mathrm{c}_{q\mathrm{S}2}^-)^2  \mathbf{q}^{\star-}
   +   (\mathrm{c}_{q\mathrm{P}}^-)^2 \mathbf{q}^-
\end{pmatrix} \,.
 \end{equation}
 From here, we can proceed by working with either expression of $ \mathbf{q}^{a-}$ and \cref{ran_hug_jump_ani_strain_mcp}, or with $ \mathbf{q}^{b-}$ and \cref{ran_hug_jump_ani_strain_mcs2}. 
\modif{The first option is chosen in the following step}, 
%for which it is useful to further 
%simplify the expression of $ \mathbf{q}^{a-}$ given in \cref{qaqb},
\begin{align*}
\alpha(\tilde{\beta} - \tilde{\alpha}) \mathbf{q}^{a-}
 &=   \tilde{\beta}\left( -\modif{\frak{B}_-}  (\mathbf{q}^{\star-} - \mathbf{q}^-)  -\mathrm{c}_{q\mathrm{S}2}^-  \mathbf{q}^{\star-}
   +   \mathrm{c}_{q\mathrm{P}}^- \mathbf{q}^-\right) \\
&\hspace*{1cm} -   \left(\modif{\frak{B}_-^2}  (\mathbf{q}^{\star-} - \mathbf{q}^-) -(\mathrm{c}_{q\mathrm{S}2}^-)^2  \mathbf{q}^{\star-}
   +   (\mathrm{c}_{q\mathrm{P}}^-)^2 \mathbf{q}^-\right) \\
 & = - \left( \tilde{\beta}\modif{\frak{B}_-}  + \modif{\frak{B}_-^2} \right)(\mathbf{q}^{\star-} - \mathbf{q}^-) 
 + \mathrm{c}^-_{q\mathrm{S}2}\left(-\tilde{\beta} + \mathrm{c}_{q\mathrm{S}2}^- \right) \mathbf{q}^{\star-} 
 + \mathrm{c}^-_{q\mathrm{P}}\left(\tilde{\beta} - \mathrm{c}_{q\mathrm{P}}^- \right) \mathbf{q}^- .
 %%
%  & = - \left( \tilde{\beta}\modif{\frak{B}_-} +\modif{\frak{B}_-^2} \right)(\mathbf{q}^{\star-} - \mathbf{q}^-) 
% - \mathrm{c}^-_{q\mathrm{S}2}\mathrm{c}^-_{q\mathrm{S}1} \mathbf{q}^{\star-} 
% + \mathrm{c}^-_{q\mathrm{P}}\left(\tilde{\beta} - \mathrm{c}_{q\mathrm{P}}^- \right) \mathbf{q}^-  \,.
 \end{align*}
\begin{equation}\label{diff_qa_ql_pre}
\Rightarrow \hspace*{0.2cm}\alpha(\tilde{\beta} - \tilde{\alpha}) \mathbf{q}^{a-}= - \left( \tilde{\beta}\modif{\frak{B}_-} +\modif{\frak{B}_-^2} \right)(\mathbf{q}^{\star-} - \mathbf{q}^-) 
 - \mathrm{c}^-_{q\mathrm{S}2}\mathrm{c}^-_{q\mathrm{S}1} \mathbf{q}^{\star-} 
 + \mathrm{c}^-_{q\mathrm{P}}\left(\tilde{\beta} - \mathrm{c}_{q\mathrm{P}}^- \right) \mathbf{q}^-  .
\end{equation}
With some algebraic manipulations\footnote{This is seen as, 
\begin{equation}
\begin{aligned}
&\left(\mathrm{c}_{q\mathrm{S}1}^--\mathrm{c}_{q\mathrm{P}}^- \right)
 \left(\mathrm{c}_{q\mathrm{S}2}^--\mathrm{c}_{q\mathrm{P}}^- \right)
  = \mathrm{c}_{q\mathrm{S}1}^-\mathrm{c}_{q\mathrm{S}2}^-
 \, - \, \mathrm{c}_{q\mathrm{P}}^-\left(\mathrm{c}_{q\mathrm{S}1}^-+\mathrm{c}_{q\mathrm{S}2}^--\mathrm{c}_{q\mathrm{P}}^-\right) \\
 \Rightarrow & \hspace*{0.5cm}\mathrm{c}^-_{q\mathrm{P}}\left(\tilde{\beta} - \mathrm{c}_{q\mathrm{P}}^- \right)  - \alpha(\tilde{\beta} - \tilde{\alpha}) 
 \, = \, \mathrm{c}^-_{q\mathrm{P}}\left(\mathrm{c}_{q\mathrm{S}1}^- + \mathrm{c}_{q\mathrm{S}2}^- - \mathrm{c}_{q\mathrm{P}}^- \right)
\, -\,  \left(\mathrm{c}_{q\mathrm{P}}^- -\mathrm{c}_{q\mathrm{S}1}^- \right)
 \left(\mathrm{c}_{q\mathrm{S}2}^--\mathrm{c}_{q\mathrm{P}}^- \right)
 \,=\, \mathrm{c}_{q\mathrm{S}1}^-\mathrm{c}_{q\mathrm{S}2}^- \,.
 \end{aligned}
 \end{equation}
 }, we obtain \modif{the following identity},
\begin{equation}
\mathrm{c}^-_{q\mathrm{P}}\left(\tilde{\beta} - \mathrm{c}_{q\mathrm{P}}^- \right)  - \alpha(\tilde{\beta} - \tilde{\alpha})  = \mathrm{c}_{q\mathrm{S}1}^-\mathrm{c}_{q\mathrm{S}2}^-  \,.
\end{equation}
\modif{We employ this to further rewrite} \cref{diff_qa_ql_pre} as, 
\begin{equation}\label{diff_qa_ql}
\alpha(\tilde{\beta} - \tilde{\alpha}) \left( \mathbf{q}^{a-} \, - \, \mathbf{q}^- \right)
 \, = \,  - \left( \tilde{\beta}\modif{\frak{B}_-} + \modif{\frak{B}_-^2} + \mathrm{c}_{q\mathrm{S}1}^-\mathrm{c}_{q\mathrm{S}2}^-\right)(\mathbf{q}^{\star-} - \mathbf{q}^-) \,.
\end{equation}

\medskip

%--------------------------
\noindent \textbf{Step 2a} Using \cref{ran_hug_jump_ani_strain_mcp}, we have \modif{the following equalities},
\begin{equation}
\left(\modif{\frak{B}_-}  +\mathrm{c}_{q\mathrm{P}}^- \right) 
\left(\textcolor{black}{\mathbf{q}^{a-}} - \mathbf{q}^-  \right) =0
\quad \Rightarrow \quad \left(\modif{\frak{B}_-}  +\mathrm{c}_{q\mathrm{P}}^- \right) 
\alpha(\tilde{\beta} - \tilde{\alpha})  \left(\textcolor{black}{\mathbf{q}^{a-}} - \mathbf{q}^-  \right) =0\,.
\end{equation}
Substitute the difference given by \cref{diff_qa_ql} into the above expression \modif{to obtain},
\begin{equation}\left(\modif{\frak{B}_-}  +\mathrm{c}_{q\mathrm{P}}^- \right) \left( \tilde{\beta}\modif{\frak{B}_-} + \modif{\frak{B}_-^2} + \mathrm{c}_{q\mathrm{S}1}^-\mathrm{c}_{q\mathrm{S}2}^-\right)(\mathbf{q}^{\star-} - \mathbf{q}^-) =0\,.
\end{equation}
In rewriting left-hand side and \modif{introducting} $\frak{p}_2^{\modif{-}}$ and $\frak{p}_3^{\modif{-}}$, 
\begin{equation}
  \frak{p}_3^{\modif{-}} :=\mathrm{c}_{q\mathrm{S}1}^-\mathrm{c}_{q\mathrm{S}2}^-\mathrm{c}_{q\mathrm{P}}^-
  \,, \quad \frak{p}_2 :=  \mathrm{c}_{q\mathrm{S}1}^-\mathrm{c}_{q\mathrm{S}2}^-
 + \mathrm{c}_{q\mathrm{S}1}^-\mathrm{c}_{q\mathrm{P}}^- 
 + \mathrm{c}_{q\mathrm{S}2}^-\mathrm{c}_{q\mathrm{P}}^- = \mathrm{c}_{q\mathrm{S}1}^-\mathrm{c}_{q\mathrm{S}2}^-+ \mathrm{c}_{q\mathrm{P}}^-\tilde{\beta},
 \end{equation}
we rewrite the above equation as,
\begin{equation}\label{rel_Aqstar_qm_0}
\left( \textcolor{black}{\frak{B}_{-}^3} \,+\,  \frak{p}_2^-\,\textcolor{black}{\frak{B}_-}
\, +\, \left( \mathrm{c}_{q\mathrm{S}1}^-+\mathrm{c}_{q\mathrm{S}2}^- + \mathrm{c}_{q\mathrm{P}}^-\right)
\textcolor{black}{\frak{B}_-^2}  \, +\, \frak{p}_3^- \mathbb{Id} \right)
 \left( \mathbf{q}^{\star-} - \mathbf{q}^-\right) = 0 \,.
\end{equation}

\medskip

\noindent \textbf{Step 2b} By algebraic computation, cf. \cite{Pham2023hdgRR}, we have \modif{the following identities for the action of powers of $\frak{B}$ on vector $\mathbf{q} = \begin{pmatrix} \overrightarrow{\boldsymbol{\epsilon}} \\  \rho\mathbf{v}\end{pmatrix}$},
\begin{equation}
\begin{aligned}
 \frak{B}(\boldsymbol{\nu}) \mathbf{q}
 = -\begin{pmatrix} 
\textcolor{black}{\mathbb{A}^t(\boldsymbol{\nu})} \mathbf{v} \\ \boldsymbol{\sigma \nu}\end{pmatrix}, \hspace*{0.3cm}
 \frak{B}(\boldsymbol{\nu})^2 \mathbf{q}
  =  \begin{pmatrix}    \rho^{-1} \mathbb{A}^t(\boldsymbol{\nu}) \boldsymbol{\sigma \nu}  \\
   \boldsymbol{\Gamma (\nu)} \mathbf{v}
  \end{pmatrix}
  , \hspace*{0.3cm} \frak{B}(\boldsymbol{\nu})^3 \mathbf{q}
 = - \dfrac{1}{\rho}\begin{pmatrix}  \mathbb{A}^t(\boldsymbol{\nu})\,  \boldsymbol{\Gamma}(\boldsymbol{\nu}) \mathbf{v} \\
 \boldsymbol{\Gamma}(\boldsymbol{\nu}) \boldsymbol{\sigma \nu}\end{pmatrix} \,.
\end{aligned}
\end{equation}
We \modif{employ these identities with $\frak{B}_-$ in} \cref{rel_Aqstar_qm_0}, 
\modif{whose last three components give},
\begin{equation}
\begin{aligned}
- \left(\dfrac{\boldsymbol{\Gamma}(\boldsymbol{\nu})^{\modif{-}}}{\rho^{\modif{-}}} +\frak{p}_2^-\right)  \left( (\boldsymbol{\sigma} \nu)^{\star-}  -(\boldsymbol{\sigma} \nu)^{-} \right)
+ \rho^-\left( \left( \modif{\mathrm{c}_{q\mathrm{S}1}^-+\mathrm{c}_{q\mathrm{S}2}^-} + \mathrm{c}_{q\mathrm{P}}^-\right)\dfrac{\boldsymbol{\Gamma}(\boldsymbol{\nu})^{\modif{-}}}{\rho^{\modif{-}}}+ \frak{p}_3^-\right)
\left(\mathbf{v}^{\star-} - \mathbf{v}^{-} \right) = 0\,\\
\Rightarrow (\boldsymbol{\sigma} \nu)^{\star-}  -(\boldsymbol{\sigma} \nu)^{-}
 = \rho^-\left(\dfrac{\boldsymbol{\Gamma}(\boldsymbol{\nu})^{\modif{-}}}{\rho^{\modif{-}}} +\frak{p}_2^-\right)^{-1} \left( \left( \modif{\mathrm{c}_{q\mathrm{S}1}^-+\mathrm{c}_{q\mathrm{S}2}^- }+ \mathrm{c}_{q\mathrm{P}}^-\right) \dfrac{\boldsymbol{\Gamma}(\boldsymbol{\nu})^{\modif{-}}}{\rho^{\modif{-}}} + \frak{p}_3^- \right)
\left(\mathbf{v}^{\star-} - \mathbf{v}^{-} \right) \,.
\end{aligned}
\end{equation}
With some algebraic manipulations, we arrive at
\begin{equation}\label{star_minus_traction::ani}
\begin{aligned}
(\boldsymbol{\sigma} \nu)^{\star-}  -(\boldsymbol{\sigma} \nu)^{-}
 &= \rho^- \left(\mathrm{c}_{q\mathrm{S}1}^-+\mathrm{c}_{q\mathrm{S}2}^- + \mathrm{c}_{q\mathrm{P}}^{\modif{-}}\right)  \left( 1+ \gamma^{\modif{-}} \left(\dfrac{\boldsymbol{\Gamma}(\boldsymbol{\nu})^{\modif{-}}}{\rho^{\modif{-}}} \,+\, \frak{p}_2^{\modif{-}} \right)^{-1}  \right)
\left(\mathbf{v}^{\star-} - \mathbf{v}^{-} \right) ,\\
& \hspace*{2cm} \text{where} \hspace*{0.3cm} \gamma^{\modif{-}}:= \dfrac{\frak{p}_3^{\modif{-}}}{\mathrm{c}_{q\mathrm{S}1}^-+\mathrm{c}_{q\mathrm{S}2}^- + \mathrm{c}_{q\mathrm{P}}^{\modif{-}} } -  \frak{p}_2^{\modif{-}}\,.
\end{aligned}
\end{equation} 
\modif{A similar relation is obtained for mesh cell $K^+$, with `$-$' replaced by `$+$' in the above expression.}

%
%\begin{equation}
%(\boldsymbol{\sigma} \nu)^{\star-}
%=(\boldsymbol{\sigma} \nu)^{-}
%+ M^{\mathrm{G-ani}} \left( \mathbf{v}^{\star-} - \mathbf{v}^-\right)
%\end{equation}

\paragraph{Derivation of HDG traces} 
From this point \modif{on, to arrive at the hybridized traces along $\mathsf{F}^{\pm}$ \cref{interFlr}, we follow the hybridization procedure employed in \cite{terrana2018spectral} for isotropy}. \modif{First, we let the numerical trace along $\mathsf{F}^{\pm}$ be given by the intermediate states,} 
\begin{equation}\label{numtracepre}
\widehat{\mathbf{v}}^{\modif{\pm}}  := \mathbf{v}^{\star \modif{\pm}}    \,, \hspace*{1cm}
\widehat{\boldsymbol{\sigma}\boldsymbol{\nu} }^{\modif{\pm}}  :=  (\boldsymbol{\sigma}\boldsymbol{\nu})^{\star\modif{\pm}} \,.
\end{equation} 
Next, employing the transmission condition  \cref{inte_eq_ani}
and relation \cref{star_minus_traction::ani}, \modif{and introducing the quantity $\lambdaV$ defined on $\mathsf{F}$ to represent $\mathbf{v}^{\star \modif{\pm}}$},
\modif{the numerical traces \cref{numtracepre} along $\mathsf{F}^-$ are rewritten as,}
\begin{equation}\label{hybrid_ani}
\widehat{\mathbf{v}}^{\modif{-}}   =\lambdaV, \hspace*{0.5cm}\text{and} \hspace*{0.5cm}
\widehat{\boldsymbol{\sigma}\boldsymbol{\nu} }^{\modif{-}} 
 =  (\boldsymbol{\sigma}\boldsymbol{\nu})^{-}
 +   \modif{M_{\mathrm{Gani}}(\boldsymbol{\nu})^-}\left( \lambdaV   - \mathbf{v}^-\right) \, ,
\end{equation}
\begin{equation} \label{eq:MGodunov-aniso}
\text{with} \hspace*{0.5cm} \modif{M_{\mathrm{Gani}}(\boldsymbol{\nu})^-}  := \rho^- \left( \cqsa^- + \cqsb^- + \mathrm{c}^-_{q\mathrm{P}}\right)  \left( 1+ \gamma^{\modif{-}} \left(\dfrac{\boldsymbol{\Gamma}(\boldsymbol{\nu})^{\modif{-}}}{\rho^-} + \frak{p}_2^-\right)^{-1}  \right) \,.
\end{equation}
\modif{The numerical traces $ \widehat{\mathbf{v}}^+$,  $ \widehat{\boldsymbol{\sigma}\boldsymbol{\nu} }^+ $ along $\mathsf{F}^+$,  are defined similarly with `$-$' replaced by `$+$' above.}

% ------------------------------------------------------------------------------- %
\section{Numerical experiments} 
\label{section:numerics}
% ------------------------------------------------------------------------------- %

The HDG method is implemented in the open-source parallel 
software \texttt{hawen}\footnote{\url{https://ffaucher.gitlab.io/hawen-website/}},
\cite{Hawen2021}. 
We carry out numerical experiments to evaluate the accuracy of 
numerical solutions depending on the choice of stabilization.
Three families of stabilization are investigated: based on the identity matrix, 
the Kelvin--Christoffel (KC) matrix $\Gambold$ (\cref{KC_iso,KC_vti}), and the 
Godunov matrix $M_{\mathrm{Godunov}}$ (\cref{eq:MGodunov-iso,eq:MGodunov-aniso}). 
Each of them is first considered with a scaling factor $\tau$:
% whose optimal value \flo{needs to be} determined. We have,
\vspace*{-.50em}
\begin{equation} \label{eq:numerics:stabilization-u:3Dpw}
\begin{aligned}
   \begin{matrix*}[l]
   & \text{Choices of $\tauU$ for} \\[.25em]
   & \widehat{\boldsymbol{\sigma \nu}} = \tensS_h - \tauU (\displacement_h - \lambdaUh),
   \end{matrix*} \quad
     \tauU \,=\, \left\lbrace \begin{matrix*}[l]
     \text{Identity-based:} & \tauSuIdRe(\tau) \,=\, \pm\,     \omega \,\,\tau \, \identity, \\
                            & \tauSuIdIm(\tau) \,\,\,=\, \pm\, \ii \omega \,  \tau \, \identity \,;  \\[.30em]
      \text{KC based:} 
    & \tauSuKCa(\tau)  \,\,=\, -\ii\omega\,\tau\, \Gambold \,; \\[.30em]
      \text{Godunov based:} 
    & \tauSuMGa(\tau)  \,\,=\, -\ii\omega\,\tau\, M_{\mathrm{Godunov}} \,.
   \end{matrix*} \right.
\end{aligned}
\end{equation}

\begin{remark}
  The scaling factor $\ii \omega$ 
  is motivated \flo{by} the relation between \flo{the} velocity 
  and \flo{the} displacement, and we refer to our extended 
  report \cite{Pham2023hdgRR} for more details. 
  % By comparing the identity-based stabilizations, we 
  % further highlight below that the imaginary unit is required.
  This result is also illustrated in the comparison within the 
  identity-based family, with $\tau$ purely real or imaginary.
\end{remark}

We consider isotropic and anisotropic materials. Our investigation consists of two steps:
\begin{enumerate}[leftmargin=*] \setlength{\itemsep}{-2pt}
  \item In \cref{subsection:numerics-3d}, we investigate the optimal choice 
        for \flo{the} scaling factor $\tau$ within each of the three families 
        of stabilization in \cref{eq:numerics:stabilization-u:3Dpw}. 
        The investigation is carried out with planewaves, 
        and thus concerns single-typed waves propagating in \flo{a} homogeneous material.
        Here, analytical solutions are used as references to evaluate
        the accuracy of the numerical simulations. 
        % We consider isotropic and VTI medium.  

  \item In \cref{subsection:numerics-2d}, we compare the best 
        candidates from each family in \flo{a} highly heterogeneous 
        medium with a Dirac point source. 
        Here, the wavefield contains all types of waves and,
        as there are no analytical solutions, 
       \flo{a reference solution is constructed by using a highly}
        % we employ a  
        % reference solution obtained with a 
        refined mesh and high-order polynomial basis functions. 
        This experiment is carried out for multiple frequencies 
        and orders of polynomial, in isotropic and TTI medium. 
        % we consider highly heterogeneous
        % medium and Dirac point-source for the propagation. 
        % Here, the wavefield contain all types of waves, and 
        % there is no analytical solution such that we use numerical 
        % reference solutions obtained on refined mesh and high-order 
        % polynomial to evaluate the accuracy.
        % We consider (2a) elastic isotropic and (2b) elastic TTI media.

\end{enumerate}

To evaluate the difference between a reference solution 
(either an analytic or a numerical one computed with
a refined mesh) and simulations, we introduce 
relative errors $\err$ and $\errall$. 
With $w$ representing a component of the displacement 
field $\displacement$ or the stress tensor $\tensS$,
and $\bx_k$ the $k^{\text{th}}$ position where the 
solutions are evaluated, we define,
\begin{equation}\label{eq:relative-error_field}
  \err[w]\,:=\, 
  \dfrac{1}{N} \,\, \sum_{k=1}^{N} \, 
  \dfrac{\vert \,w^\mathrm{ref}(\bx_k) \,-\, w(\bx_k) \,\vert}{\Vert w^\mathrm{ref} \Vert}\,,
  \qquad \text{with } \qquad 
  \Vert w^\mathrm{ref} \Vert \,=\, \sqrt{\sum_{k=1}^{N} \vert w^\mathrm{ref} (\bx_k)\vert^2 }
  \,\,\, .
\end{equation}
The error for the total field is given by,
\begin{equation}\label{eq:relative-error_sum}
  \errall[\displacement] \,:=\, 
  \dfrac{1}{\mathfrak{n}_d} \,\, \sum_{j\in I} \,\, \err[u_j] \,; \qquad\qquad
  \errall[\tensS] \,:=\,   \dfrac{1}{\mathfrak{N}_d} \,\, \sum_{j\in \mathcal{I}_{\mathrm{sm}}} \,\,
                                          \err[\sigma_j] \,,
\end{equation}
with $\mathfrak{n}_d=2$ or $3$ for two and three dimensions respectively, 
and $\mathfrak{N}_d=3$ or $6$.
%\begin{subequations}\label{eq:relative-error_all}
%\begin{align}
%& \err[w]\,:=\, 
%  \dfrac{1}{N} \,\, \sum_{k=1}^{N} \, 
%  \dfrac{\vert \,w^\mathrm{ref}(\bx) \,-\, w(\bx) \,\vert}{\Vert w^\mathrm{ref} \Vert}\,; 
%  \qquad \Vert w^\mathrm{ref} \Vert \,=\, \ldots \\
%%  \widehat{\err}(w; \bx_k) \,,
%%  \qquad
%%  \text{with} \qquad
%%  \widehat{\err}_w(\bx) \,=\, \dfrac{\vert \,w^\mathrm{ref}(\bx) \,-\, w(\bx) \,\vert}
%%                                 {\Vert w^\mathrm{ref} \Vert}\,; \\
%& \errall[\displacement] \,:=\, 
%  \dfrac{1}{3} \,\, \sum_{j=\{x,y,z\}} \,\, \err[u_j] \,; \qquad
%  \errall[\tensS] \,:=\,   \dfrac{1}{6} \,\, \sum_{j=\{xx,yy,zz,xy,xz,yz\}} \,\,
%                                          \err[\sigma_j] \,.
%\end{align}
%\end{subequations}
In the experiments, we use a Cartesian grid 
for the positions $\bx_k$ where the solutions are evaluated 
% in order to balance the space contributions, 
even though our computational mesh can be unstructured. 
We however exclude the positions near the boundaries, and near 
% the source when we have a 
Dirac point source (i.e., for experiments of \cref{subsection:numerics-2d}).

% ------------------------------------------------------------------------------- %
\subsection{3D experiments with planewaves}
\label{subsection:numerics-3d}
% isotropy
%\renewcommand{\plottauSuIda}{tau6}
%\renewcommand{\plottauSuIdb}{tau4}
%\renewcommand{\plottauSuKCa}{tau16}
%\renewcommand{\plottauSuKCb}{tau13}
%\renewcommand{\plottauSuMGa}{tau23}
%\renewcommand{\plottauSvIda}{tau0}
%\renewcommand{\plottauSvIdb}{tau1}
%\renewcommand{\plottauSvKCb}{tau11}
%\renewcommand{\plottauSvMGa}{tau20}
% ------------------------------------------------------------------------------- %

We consider % the propagation of elastic planewaves in 
a three-dimensional domain $\Omega$ corresponding to the 
cube $(-1,1)^3$, with boundary denoted by $\partial\Omega$.
% Reminding \cref{1stsystem_variants:u}, we have,
\modif{We consider equation \cref{1stsystem_variants:u} on $\Omega$
with boundary condition
$\displacement \,=\, \displacement_{\mathrm{pw}}$  on $\partial\Omega$.}
%\begin{subequations} \label{eq:numerics:eq-elastic-planewave}
%\begin{empheq}[left = \empheqlbrace\,]{align}
%  & - \omega^2 \,\rho \,\displacement  \, - \,  \nabla \cdot \tensS  \,=\, 0 \,,\qquad \text{in $\Omega$,} \\
%  &   \tensCS \tensS \,=\,  \nabla^\mathfrak{s}\displacement \,,\hspace*{6.9em} \text{in $\Omega$,} \\
%  &   \displacement \,=\, \displacement_{\mathrm{pw}} \,, \hspace*{7.9em}\text{on $\partial\Omega$} \,.
%\end{empheq}
%\end{subequations}
We consider different types of planewaves 
% depending on the definition of the function $\displacement_{\mathrm{pw}}$ 
% for the Dirichlet boundary condition. We 
and refer to~\cite[Appendix~A]{Pham2023hdgRR} 
for more details on the derivation of planewaves in linear elasticity.
 
% ------------------------------------------------------------------------------- %
% \renewcommand{\experiment}{\texttt{benchmark3Diso-Pw}}
% \subsubsection{Elastic isotropy (\experiment)}
\renewcommand{\experiment}{\texttt{Configuration~PwIso}}
\subsubsection{Elastic isotropic medium}
\label{subsection:iso-3D-planewaves}
% ------------------------------------------------------------------------------- %

We consider % an isotropic elastic material and 
the propagation of P- and S-planewaves, 
% given by the solutions to \cref{eq:numerics:eq-elastic-planewave}
\flo{with corresponding} Dirichlet condition $\displacement_{\mathrm{pw}}^\mathrm{P}$
and $\displacement_{\mathrm{pw}}^\mathrm{S}$ respectively. 
We select the planewave direction $\mathbf{d}=\big(1/\sqrt{2},0,1/\sqrt{2}\big)^t$,
and further impose a strong contrast between the P- and S-wavespeeds, 
with the following configuration:
\begin{equation}\label{benchmark02}
%% \experiment\quad 
\left\lbrace
\begin{aligned}
%% & \text{Cube $(-1,1)^3$ with planewave} \\  % \nonumber \\
  & \cp \,=\, \num{2.5e-3}\si{\meter\per\second}   \,, \quad
    \cs \,=\, \num{e-4}   \si{\meter\per\second}   \,, \quad
    \rho\,=\, \num{1}     \si{\kg\per\meter\cubed} \,, \\
 &  \qquad\qquad \Leftrightarrow \,\, 
    \lambda\,=\, \num{6.24e-06}\si{\pascal}\quad\text{and}\quad \mu\,=\,\num{e-8}\si{\pascal} \,; \\[.30em]
 &  \text{P-planewave propagation:} \qquad 
    \displacement_{\mathrm{pw}}^\mathrm{P}(\bx)\,=\,(1, 0, 1)^t \,\, \, e^{\ii\frac{\omega}{\cp}\, (\mathbf{d} \, \cdot \, \bx)}\,,\\[.20em]
 &  \text{S-planewave propagation:} \qquad 
    \displacement_{\mathrm{pw}}^\mathrm{S}(\bx)\,=\,(0, 1, 0)^t \,\, \, e^{\ii\frac{\omega}{\cs}\, (\mathbf{d} \, \cdot \, \bx)}\,.
\end{aligned} \right. \end{equation}
The simulations use a mesh consisting of \num{48000} tetrahedra 
and polynomials of order 4. % for the discretization. 
Due to the strong contrast between the P- and S-wavelength, 
we use frequency $\omega/(2\pi)=\num{8}$\si{\milli\Hz}
for the P-planewave, and $\omega/(2\pi)=\num{0.4}$\si{\milli\Hz} for the 
S-planewave, resulting in wavelength of size $\num{0.3}\si{\meter}$ and $\num{0.25}\si{\meter}$,
respectively. 

\begin{remark}\label{remark:scaling}

Our experiment \cref{benchmark02} on a unitary cube 
is equivalent to working with a cube 
$(-\num{1},\num{1})\si{\km\cubed}$ with 
$\cp =\num{2500}\si{\meter\per\second}$ and 
$\cs =\num{100} \si{\meter\per\second}$,
using P-planewave frequency \num{8}\si{\Hz} and 
S-planewave frequency \num{0.4}\si{\Hz}.
This is obtained by maintaining the number of wavelengths 
propagating in \flo{a} domain of size $L^3$, which is given by
$(L \,\times\, \text{frequency})/wave speed$.
%\begin{equation}
%  \text{Number of wavelengths} =
%  \dfrac{L \,\times\, \text{frequency}}{wave speed} \,.
%\end{equation}

\end{remark}

In \cref{figure:3D-tau-amplitude:isotropic}, 
we investigate the accuracy of the solution % for varying scaling factor $\tau$
\flo{as the scaling factor $\tau$ in \cref{eq:numerics:stabilization-u:3Dpw} varies}.
% in the stabilization coefficients of \cref{eq:numerics:stabilization-u:3Dpw}.
For the identity-based stabilization, 
we investigate the scaling factor on the purely real and complex lines 
and determine the optimal sign.
For the sake of conciseness, the relative error $\err$ 
is shown only for some of the non-zero wave fields: 
$u_x$ and $\sigma_{zz}$ for the P-planewave, and 
$u_y$ for the S-planewave.
% identity-based stabilization \cref{eq:numerics:stabilization-u},
% with $\tau$ varying from \num{e-9} to \num{e0}.
% In vertical lines, we also indicate the error obtained with 
% the Kelvin--Christoffel stabilization $\tauSuMGa$ without scaling
% (i.e., with $\tau_{\mathrm{kc}}=1$ in \cref{eq:numerics:stabilization-u}),
% and with the Godunov stabilization $\tauSuMGa$~\cref{eq:numerics:stabilization-u},
% with the Godunov matrix for isotropy given in \cref{eq:MGodunov-iso}.
% In \cref{figure:3D-tau-amplitude:iso02_qP-Su,figure:3D-tau-amplitude:iso02_sH-Su},
% we show the relative error for each of the non-zero wave fields, for the 
% P- and S-planewave propagation respectively.
% In these plots, we also indicate the error level 
% obtained with the  Godunov and Kelvin--Christoffel 
% stabilization $\tauSuMGa$ and $\tauSuKCa$ (horizontal lines).
% We also show the values of the P- and S-impedance multiplied 
% by the frequency (vertical lines), the latter provides the 
% error using stabilization $\tauSuIdb$ in {eq:numerics:stabilization-u}.

% ---------------------------------------------------------------------- %
\setlength{\plotwidth} {4.50cm}
\setlength{\plotheight}{3.50cm}
% ---------------------------------------------------------------------- %
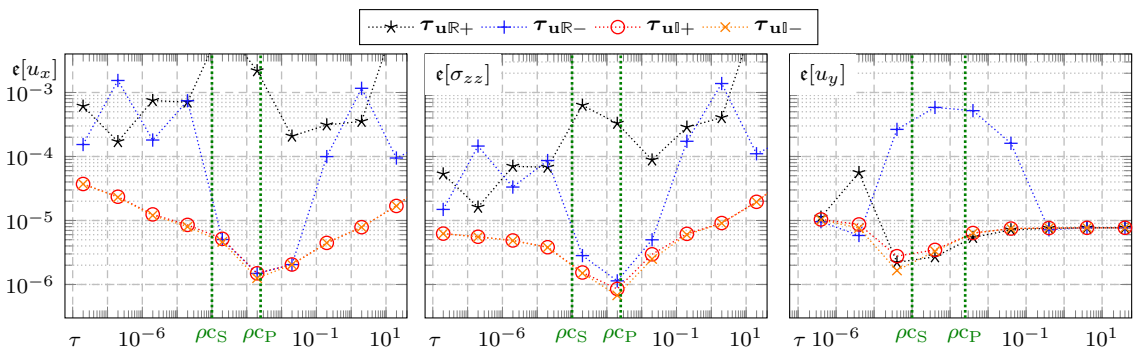
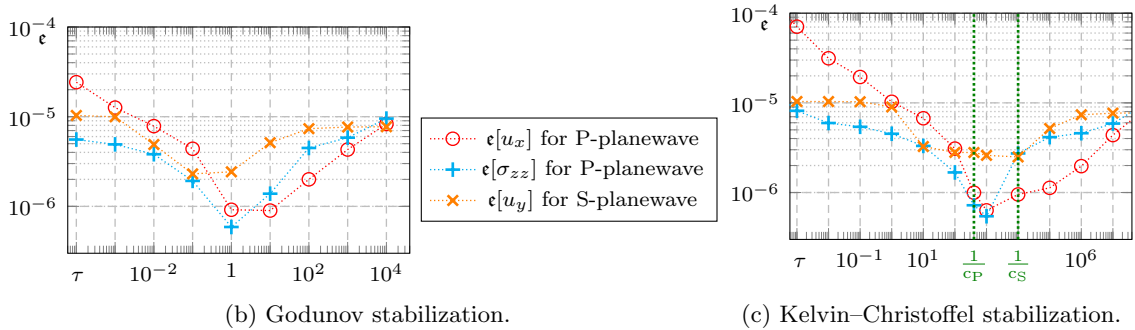
\begin{figure}[ht!] \centering
  \renewcommand{\datafile}{figures/numerics/3D-planewave_optimal-stabilization/data/data_elastic-iso_TauId.txt}
  \renewcommand{\myxlabel}{$\tau$}
  \renewcommand{\myylabel}{$\err$}

  \pgfkeys{/pgf/fpu=true}
  \renewcommand{\qPomega}{0.00250} % Rho Cp only
  \renewcommand{\qSomega}{0.00010} % Rho Cs only
  \pgfmathsetmacro{\scaleomega}{0.0503}
  \pgfmathsetmacro{\xmin}{0.60e-8} \pgfmathsetmacro{\xmax}{0.40e2} 
  \pgfmathsetmacro{\ymin}{3e-7}    \pgfmathsetmacro{\ymax}{4e-3} 
  \pgfkeys{/pgf/fpu=false}
  \renewcommand{\dataA}{qP-double_ux_freq8.00mHz} 
  \subfloat[][{Identity-based stabilization using real and imaginary scaling factors. 
               Relative error for P-planewave $\err[u_x]$ (left) and $\err[\sigma_{zz}]$ (middle);
               for S-planewave $\err[u_y]$ (right).}]{
              \renewcommand{\myylabel}{${\err[u_x]}$}
              %% --------------------------------------------------------------------------------------
\begin{tikzpicture}
%% --------------------------------------------------------------------------------------
\begin{axis}[
             enlargelimits=false, 
             xlabel={\scriptsize \myxlabel },
             ylabel={\scriptsize \myylabel },
             enlarge y limits=false,
             enlarge x limits=false,
             xmin=\xmin,xmax=\xmax, 
             ymin=\ymin,ymax=\ymax,             
             yminorticks=true,
             height=\plotheight,width=\plotwidth,
             scale only axis,
             ylabel style={rotate=-90},
             ylabel style = {yshift=4.0em, xshift=2.50em},
             xlabel style = {yshift = 1em, xshift=-5.50em},
             clip mode=individual,
             ytick={1e-6,1e-5,1e-4,1e-3},
             xtick={1e-8,1e-7,1e-6,1e-5,1e-4,1e-3,1e-2,1e-1,1e0,1e1,1e2,1e3},
             xticklabels={,,\num{e-6},,,,,\num{e-1},,\num{e1},,},
             minor grid style={line width=.5pt, draw=gray!50, densely dotted},
             major grid style={line width=.5pt, draw=gray!50, densely dashed},
             xmajorgrids=true,
           % xminorgrids=true,
             yminorgrids=true,
             ymajorgrids=true,             ymode=log,
             xmode=log,
             %y  tick label style={/pgf/number format/fixed},
             %y  tick label style={/pgf/number format/fixed zerofill},
             %y  tick label style={/pgf/number format/precision=1},
             %y  tick label style={/pgf/number format/precision=2},
             % y  tick label style={/pgf/number format/sci},
             %y tick label style={/pgf/number format/.cd},
             label style={font=\scriptsize},
             tick label style={font=\scriptsize},
             legend style={font=\scriptsize\selectfont},
             legend style={at={(0.99,1.17)}},
             legend columns=4,
             ]  %% foreground marks

%     \pgfmathsetmacro{\scale}{1.0}
     \addplot[color=black,densely dotted, line width=0.5, mark=star,% only marks,
              mark size=2.50,mark options={solid,fill=black, line width=.50}]
              table[x expr = \thisrow{tau}/\scaleomega, 
                    y expr = \thisrow{\dataA_plusRe}]
              {\datafile}; % \addlegendentry{$\tau \identity\,\,\,$}
     \addplot[color=\myblue,densely dotted, line width=0.5, mark=+, % only marks,
              mark size=2.50,mark options={solid,fill=\myblue, line width=.50}]
              table[x expr = \thisrow{tau}/\scaleomega, 
                    y expr = \thisrow{\dataA_minusRe}]
              {\datafile}; % \addlegendentry{$-\tau \identity \,\,\,$}

     \addplot[color=red,densely dotted, line width=0.5, mark=o, % only marks,
              mark size=2.50,mark options={solid,fill=red, line width=.50}]
              table[x expr = \thisrow{tau}/\scaleomega, 
                    y expr = \thisrow{\dataA_plusIm}]
              {\datafile}; % \addlegendentry{$\ii\tau \identity\,\,\,$}
     \addplot[color=orange,densely dotted, line width=0.5, mark=x, % only marks,
              mark size=2.50,mark options={solid,fill=orange, line width=.50}]
              table[x expr = \thisrow{tau}/\scaleomega, 
                    y expr = \thisrow{\dataA_minusIm}]
              {\datafile}; % \addlegendentry{$-\ii\tau \identity\,\,\,$}

     \draw[densely dotted, \mygreen, line width=1] (\qPomega,\ymin) -- (\qPomega,\ymax);
     \draw[densely dotted, \mygreen, line width=1] (\qSomega,\ymin) -- (\qSomega,\ymax);
     \node[\mygreen, line width=1,yshift=0em,anchor=north,align=center] at (\qPomega,\ymin)  
          {\scriptsize $\rho \cp$};
     \node[\mygreen, line width=1,yshift=0em,anchor=north,align=center] at (\qSomega,\ymin)  
          {\scriptsize $\rho \cs$};
\end{axis}
\end{tikzpicture}
              \renewcommand{\dataA}{qP-double_szz_freq8.00mHz} \hspace*{-3.50em}
              \renewcommand{\myylabel}{${\err[\sigma_{zz}]}$}
              %% --------------------------------------------------------------------------------------
\begin{tikzpicture}
%% --------------------------------------------------------------------------------------
\begin{axis}[
             enlargelimits=false, 
             xlabel={\scriptsize \myxlabel },
             enlarge y limits=false,
             enlarge x limits=false,
             xmin=\xmin,xmax=\xmax, 
             ymin=\ymin,ymax=\ymax,             
             yminorticks=true,
             height=\plotheight,width=\plotwidth,
             scale only axis,
             ylabel style={rotate=-90},
             ylabel style = {yshift=4.0em, xshift=1.00em},
             xlabel style = {yshift = 1em, xshift=-5.50em},
             clip mode=individual,
             ytick={1e-6,1e-5,1e-4,1e-3,1e-2},
             xtick={1e-8,1e-7,1e-6,1e-5,1e-4,1e-3,1e-2,1e-1,1e0,1e1,1e2,1e3},
             xticklabels={,,\num{e-6},,,,,\num{e-1},,\num{e1},,},
             yticklabels={,,},
             minor grid style={line width=.5pt, draw=gray!50, densely dotted},
             major grid style={line width=.5pt, draw=gray!50, densely dashed},
             xmajorgrids=true,
           % xminorgrids=true,
             yminorgrids=true,
             ymajorgrids=true,             ymode=log,
             xmode=log,
             %y  tick label style={/pgf/number format/fixed},
             %y  tick label style={/pgf/number format/fixed zerofill},
             %y  tick label style={/pgf/number format/precision=1},
             %y  tick label style={/pgf/number format/precision=2},
             % y  tick label style={/pgf/number format/sci},
             %y tick label style={/pgf/number format/.cd},
             label style={font=\scriptsize},
             tick label style={font=\scriptsize},
             legend style={font=\scriptsize\selectfont},
             legend style={at={(1.15,1.17)}},
             legend columns=4,
             ]  %% foreground marks

%     \pgfmathsetmacro{\scale}{1.0}
     \addplot[color=black,densely dotted, line width=0.5, mark=star, % only marks,
              mark size=2.50,mark options={solid,fill=black, line width=.50}]
              table[x expr = \thisrow{tau}/\scaleomega, 
                    y expr = \thisrow{\dataA_plusRe}]
              {\datafile}; \addlegendentry{$\tauSuIdRep$}
     \addplot[color=\myblue,densely dotted, line width=0.5, mark=+, % only marks,
              mark size=2.50,mark options={solid,fill=\myblue, line width=.50}]
              table[x expr = \thisrow{tau}/\scaleomega, 
                    y expr = \thisrow{\dataA_minusRe}]
              {\datafile}; \addlegendentry{$\tauSuIdRem$}

     \addplot[color=red,densely dotted, line width=0.5, mark=o, % only marks,
              mark size=2.50,mark options={solid,fill=red, line width=.50}]
              table[x expr = \thisrow{tau}/\scaleomega, 
                    y expr = \thisrow{\dataA_plusIm}]
              {\datafile}; \addlegendentry{$\tauSuIdImp$}
     \addplot[color=orange,densely dotted, line width=0.5, mark=x, % only marks,
              mark size=2.50,mark options={solid,fill=orange, line width=.50}]
              table[x expr = \thisrow{tau}/\scaleomega, 
                    y expr = \thisrow{\dataA_minusIm}]
              {\datafile}; \addlegendentry{$\tauSuIdImm$}

     \node[yshift=0em,anchor=north west,align=center,fill=white] at (\xmin,\ymax)  
          {\scriptsize \myylabel};

     \draw[densely dotted, \mygreen, line width=1] (\qPomega,\ymin) -- (\qPomega,\ymax);
     \draw[densely dotted, \mygreen, line width=1] (\qSomega,\ymin) -- (\qSomega,\ymax);
     \node[\mygreen, line width=1,yshift=0em,anchor=north,align=center] at (\qPomega,\ymin)  
          {\scriptsize $\rho \cp$};
     \node[\mygreen, line width=1,yshift=0em,anchor=north,align=center] at (\qSomega,\ymin)  
          {\scriptsize $\rho \cs$};
\end{axis}
\end{tikzpicture}
              \hspace*{-5.0em}
              \pgfkeys{/pgf/fpu=true}
              \pgfmathsetmacro{\xmin}{0.60e-7} \pgfmathsetmacro{\xmax}{0.60e2} 
              \pgfmathsetmacro{\scaleomega}{0.0025}
              \pgfkeys{/pgf/fpu=false}
              \renewcommand{\dataA}{sH-double_uy_freq0.40mHz} 
              \renewcommand{\myylabel}{${\err[u_y]}$}
              %% --------------------------------------------------------------------------------------
\begin{tikzpicture}
%% --------------------------------------------------------------------------------------
\begin{axis}[
             enlargelimits=false, 
             xlabel={\scriptsize \myxlabel },
             enlarge y limits=false,
             enlarge x limits=false,
             xmin=\xmin,xmax=\xmax, 
             ymin=\ymin,ymax=\ymax,             
             yminorticks=true,
             height=\plotheight,width=\plotwidth,
             scale only axis,
             ylabel style={rotate=-90},
             ylabel style = {yshift=4.0em, xshift=1.00em},
             xlabel style = {yshift =1em, xshift=-5.50em},
             clip mode=individual,
             ytick={1e-6,1e-5,1e-4,1e-3,1e-2},
             xtick={1e-8,1e-7,1e-6,1e-5,1e-4,1e-3,1e-2,1e-1,1e0,1e1,1e2,1e3,1e4,1e5,1e6,1e7},
             xticklabels={,,\num{e-6},,,,,\num{e-1},,\num{e1},},
             yticklabels={,,,,},
             minor grid style={line width=.5pt, draw=gray!50, densely dotted},
             major grid style={line width=.5pt, draw=gray!50, densely dashed},
             xmajorgrids=true,
           % xminorgrids=true,
             yminorgrids=true,
             ymajorgrids=true,             ymode=log,
             xmode=log,
             %y  tick label style={/pgf/number format/fixed},
             %y  tick label style={/pgf/number format/fixed zerofill},
             %y  tick label style={/pgf/number format/precision=1},
             %y  tick label style={/pgf/number format/precision=2},
             % y  tick label style={/pgf/number format/sci},
             %y tick label style={/pgf/number format/.cd},
             label style={font=\scriptsize},
             tick label style={font=\scriptsize},
             legend style={font=\scriptsize\selectfont},
             legend style={at={(0.99,1.17)}},
             legend columns=4,
             ]  %% foreground marks

%     \pgfmathsetmacro{\scale}{1.0}
     \addplot[color=black,densely dotted, line width=0.5, mark=star, %only marks,
              mark size=2.50,mark options={solid,fill=black, line width=.50}]
              table[x expr = \thisrow{tau}/\scaleomega, 
                    y expr = \thisrow{\dataA_plusRe}]
              {\datafile}; %\addlegendentry{$\tau \identity\,\,\,$}
     \addplot[color=\myblue,densely dotted, line width=0.5, mark=+, %only marks,
              mark size=2.50,mark options={solid,fill=\myblue, line width=.50}]
              table[x expr = \thisrow{tau}/\scaleomega, 
                    y expr = \thisrow{\dataA_minusRe}]
              {\datafile}; %\addlegendentry{$-\tau \identity \,\,\,$}

     \addplot[color=red,densely dotted, line width=0.5, mark=o, %only marks,
              mark size=2.50,mark options={solid,fill=red, line width=.50}]
              table[x expr = \thisrow{tau}/\scaleomega, 
                    y expr = \thisrow{\dataA_plusIm}]
              {\datafile}; %\addlegendentry{$\ii\tau \identity\,\,\,$}
     \addplot[color=orange,densely dotted, line width=0.5, mark=x, %only marks,
              mark size=2.50,mark options={solid,fill=orange, line width=.50}]
              table[x expr = \thisrow{tau}/\scaleomega, 
                    y expr = \thisrow{\dataA_minusIm}]
              {\datafile}; %\addlegendentry{$-\ii\tau \identity\,\,\,$}

     \node[yshift=0em,anchor=north west,align=center,fill=white] at (\xmin,\ymax)  
          {\scriptsize \myylabel};

     \draw[densely dotted, \mygreen, line width=1] (\qPomega,\ymin) -- (\qPomega,\ymax);
     \draw[densely dotted, \mygreen, line width=1] (\qSomega,\ymin) -- (\qSomega,\ymax);
     \node[\mygreen, line width=1,yshift=0em,anchor=north,align=center] at (\qPomega,\ymin)  
          {\scriptsize $\rho \cp$};
     \node[\mygreen, line width=1,yshift=0em,anchor=north,align=center] at (\qSomega,\ymin)  
          {\scriptsize $\rho \cs$};
\end{axis}
\end{tikzpicture}
              \label{figure:3D-tau-amplitude:isotropic:id}}

  %% =======================================================================
  \setlength{\plotheight}{3.00cm}
  \pgfkeys{/pgf/fpu=true}
  \pgfmathsetmacro{\xmin}{0.60e-4} \pgfmathsetmacro{\xmax}{0.40e5} 
  \pgfmathsetmacro{\ymax}{1e-4} 
  \pgfkeys{/pgf/fpu=false}
  \renewcommand{\dataA}{qP-double_ux_freq8.00mHz}
  \renewcommand{\dataB}{sH-double_uy_freq0.40mHz}
  \renewcommand{\dataC}{qP-double_szz_freq8.00mHz}
  \renewcommand{\datafile}
               {figures/numerics/3D-planewave_optimal-stabilization/data/data_elastic-iso_TauGodunov.txt}
  \renewcommand{\myxlabel}{$\tau$}
  \subfloat[][Godunov stabilization.]
             {%% --------------------------------------------------------------------------------------
\begin{tikzpicture}
%% --------------------------------------------------------------------------------------
\begin{axis}[
             enlargelimits=false, 
             xlabel={\scriptsize \myxlabel },
             ylabel={\scriptsize \myylabel },
             enlarge y limits=false,
             enlarge x limits=false,
             xmin=\xmin,xmax=\xmax, 
             ymin=\ymin,ymax=\ymax,             
             yminorticks=true,
             height=\plotheight,width=\plotwidth,
             scale only axis,
             ylabel style={rotate=-90},
             ylabel style = {yshift=3.5em, xshift=2em},
             xlabel style = {yshift=1em, xshift=-5.50em},
             clip mode=individual,
             ytick={1e-6,1e-5,1e-4,1e-3,1e-2,1e-1},
             xtick={1e-5,1e-4,1e-3,1e-2,1e-1,1e0,1e1,1e2,1e3,1e4,1e5},
             xticklabels={,,,\num{e-2},,\num{1},,\num{e2},,\num{e4},},
             minor grid style={line width=.5pt, draw=gray!50, densely dotted},
             major grid style={line width=.5pt, draw=gray!50, densely dashed},
             xmajorgrids=true,
           % xminorgrids=true,
             yminorgrids=true,
             ymajorgrids=true,             ymode=log,
             xmode=log,
             %y  tick label style={/pgf/number format/fixed},
             %y  tick label style={/pgf/number format/fixed zerofill},
             %y  tick label style={/pgf/number format/precision=1},
             %y  tick label style={/pgf/number format/precision=2},
             % y  tick label style={/pgf/number format/sci},
             %y tick label style={/pgf/number format/.cd},
             label style={font=\scriptsize},
             tick label style={font=\scriptsize},
             legend style={font=\scriptsize\selectfont},
             legend style={at={(1.900,0.60)}},
             legend columns=1,
             ]  %% foreground marks

%     \pgfmathsetmacro{\scale}{1.0}
%     \addplot[color=black,densely dashed, line width=0.2, mark=star, only marks,
%              mark size=2.50,mark options={solid,fill=black, line width=.50}]
%              table[x expr = \thisrow{tau}, 
%                    y expr = \thisrow{\dataA}]
%              {\datafile}; \addlegendentry{P-planewave}
%     \addplot[color=\myblue,densely dashed, line width=0.2, mark=+, only marks,
%              mark size=2.50,mark options={solid,fill=\myblue, line width=.50}]
%              table[x expr = \thisrow{tau}, 
%                    y expr = \thisrow{\dataB}]
%              {\datafile}; \addlegendentry{S-planewave}

     \addplot[color=red,densely dotted, line width=0.5, mark=o, % only marks,
              mark size=2.50,mark options={solid,fill=red, line width=.50}]
              table[x expr = \thisrow{tau}, 
                    y expr = \thisrow{\dataA}]
              {\datafile}; \addlegendentry{$\err[u_x]$ for P-planewave}
     \addplot[color=cyan,densely dotted, line width=0.5, mark=+,% only marks,
              mark size=2.75,mark options={solid,fill=cyan, line width=1}]
              table[x expr = \thisrow{tau}, 
                    y expr = \thisrow{\dataC}]
              {\datafile}; \addlegendentry{$\err[\sigma_{zz}]$ for P-planewave}
     \addplot[color=orange,densely dotted, line width=0.5, mark=x,% only marks,
              mark size=2.75,mark options={solid,fill=orange, line width=1}]
              table[x expr = \thisrow{tau}, 
                    y expr = \thisrow{\dataB}]
              {\datafile}; \addlegendentry{$\err[u_y]$ for S-planewave}
              
%     \draw[densely dotted, magenta, line width=1] (\qPomega,\ymin) -- (\qPomega,\ymax);
%     \draw[densely dotted, magenta, line width=1] (\qSomega,\ymin) -- (\qSomega,\ymax);
%     \node[magenta, line width=1,yshift=0em,anchor=north,align=center] at (\qPomega,\ymin)  
%          {\scriptsize $\rho \omega c_p$};
%     \node[magenta, line width=1,yshift=0em,anchor=north,align=center] at (\qSomega,\ymin)  
%          {\scriptsize $\rho \omega c_s$};
\end{axis}
\end{tikzpicture}
              \label{figure:3D-tau-amplitude:isotropic:g}} 
             \hspace*{-4.50em}
  \pgfkeys{/pgf/fpu=true}
  \pgfmathsetmacro{\xmin}{0.60e-3} \pgfmathsetmacro{\xmax}{0.40e8} 
  \renewcommand{\qPomega}{400} % 1/(Rho Cp)
  \renewcommand{\qSomega}{1e4} % 1/(Rho Cs)
  \pgfkeys{/pgf/fpu=false}
  \renewcommand{\datafile}
               {figures/numerics/3D-planewave_optimal-stabilization/data/data_elastic-iso_TauKC.txt}
  \subfloat[][Kelvin--Christoffel stabilization.]
             {\makebox[.40\linewidth][c]{{
              %% --------------------------------------------------------------------------------------
\begin{tikzpicture}
%% --------------------------------------------------------------------------------------
\begin{axis}[
             enlargelimits=false, 
             xlabel={\scriptsize \myxlabel },
             ylabel={\scriptsize \myylabel },
             enlarge y limits=false,
             enlarge x limits=false,
             xmin=\xmin,xmax=\xmax, 
             ymin=\ymin,ymax=\ymax,             
             yminorticks=true,
             height=\plotheight,width=\plotwidth,
             scale only axis,
             ylabel style={rotate=-90},
             ylabel style = {yshift=3.5em, xshift=2em},
             xlabel style = {yshift =1em, xshift=-5.50em},
             clip mode=individual,
             ytick={1e-6,1e-5,1e-4,1e-3,1e-2,1e-1},
             xtick={1e-5,1e-4,1e-3,1e-2,1e-1,1e0,1e1,1e2,1e3,1e4,1e5,1e6,1e7,1e8,1e9},
             xticklabels={,,,,\num{e-1},,\num{e1},,,,,\num{e6},,\num{e8},},
             minor grid style={line width=.5pt, draw=gray!50, densely dotted},
             major grid style={line width=.5pt, draw=gray!50, densely dashed},
             xmajorgrids=true,
           % xminorgrids=true,
             yminorgrids=true,
             ymajorgrids=true,             ymode=log,
             xmode=log,
             %y  tick label style={/pgf/number format/fixed},
             %y  tick label style={/pgf/number format/fixed zerofill},
             %y  tick label style={/pgf/number format/precision=1},
             %y  tick label style={/pgf/number format/precision=2},
             % y  tick label style={/pgf/number format/sci},
             %y tick label style={/pgf/number format/.cd},
             label style={font=\scriptsize},
             tick label style={font=\scriptsize},
             legend style={font=\scriptsize\selectfont},
             legend style={at={(0.99,1.17)}},
             legend columns=4,
             ]  %% foreground marks

%     \pgfmathsetmacro{\scale}{1.0}
%     \addplot[color=black,densely dashed, line width=0.2, mark=star, only marks,
%              mark size=2.50,mark options={solid,fill=black, line width=.50}]
%              table[x expr = \thisrow{tau}, 
%                    y expr = \thisrow{\dataA}]
%              {\datafile}; \addlegendentry{P-planewave}
%     \addplot[color=\myblue,densely dashed, line width=0.2, mark=+, only marks,
%              mark size=2.50,mark options={solid,fill=\myblue, line width=.50}]
%              table[x expr = \thisrow{tau}, 
%                    y expr = \thisrow{\dataB}]
%              {\datafile}; \addlegendentry{S-planewave}

     \addplot[color=red,densely dotted, line width=0.5, mark=o, % only marks,
              mark size=2.50,mark options={solid,fill=red, line width=.50}]
              table[x expr = \thisrow{tau}, 
                    y expr = \thisrow{\dataA}]
              {\datafile}; % \addlegendentry{}
     \addplot[color=cyan,densely dotted, line width=0.5, mark=+,% only marks,
              mark size=2.75,mark options={solid,fill=cyan, line width=1}]
              table[x expr = \thisrow{tau}, 
                    y expr = \thisrow{\dataC}]
              {\datafile}; 
     \addplot[color=orange,densely dotted, line width=0.5, mark=x,% only marks,
              mark size=2.75,mark options={solid,fill=orange, line width=1}]
              table[x expr = \thisrow{tau}, 
                    y expr = \thisrow{\dataB}]
              {\datafile}; % \addlegendentry{}

     \draw[densely dotted, \mygreen, line width=1] (\qPomega,\ymin) -- (\qPomega,\ymax);
     \draw[densely dotted, \mygreen, line width=1] (\qSomega,\ymin) -- (\qSomega,\ymax);
     \node[\mygreen, line width=1,yshift=0em,anchor=north,align=center] at (\qPomega,\ymin)  
          {\scriptsize $\tfrac{1}{\cp}$};
     \node[\mygreen, line width=1,yshift=0em,anchor=north,align=center] at (\qSomega,\ymin)  
          {\scriptsize $\tfrac{1}{\cs}$};
\end{axis}
\end{tikzpicture}}}
              \label{figure:3D-tau-amplitude:isotropic:kc}}

  \caption{Investigation of the optimal scaling $\tau$ within the three
           families of stabilization in \cref{eq:numerics:stabilization-u:3Dpw} 
           for planewave propagation in isotropic medium using 
           parameters~\cref{benchmark02}.
           Our criterion of accuracy is the relative error 
           $\err$ \cref{eq:relative-error_field}.
           }
  \label{figure:3D-tau-amplitude:isotropic}
\end{figure}
% ---------------------------------------------------------------------- %

We can draw the following observations.
\begin{itemize}\renewcommand{\itemsep}{-2pt} % [leftmargin=*]
 \item From the plots using the identity-based stabilization, 
       \cref{figure:3D-tau-amplitude:isotropic:id}, it is clear 
       that the real-valued stabilization 
       $\tauSuIdRe$ is inaccurate
       % while the imaginary-one $\tauSuIdIm$ can be. 
       \modif{compared to the imaginary one $\tauSuIdIm$.}
       On the other hand, the choice of sign ($\pm \ii$) 
       barely affects the accuracy.

 \item The optimal value of the coefficient $\tau$ 
       in the identity-based stabilization $\tauSuIdIm(\tau)$
       \cref{eq:numerics:stabilization-u:3Dpw}
       depends on the type of waves propagating 
       and appears to be given by the impedance. 
       Namely, for P-planewave, the optimal coefficient would be
       the P-impedance, while it is the S-impedance for the S-planewave,
       see~\cref{figure:3D-tau-amplitude:isotropic:id}.

 \item For the Kelvin--Christoffel stabilization, \cref{figure:3D-tau-amplitude:isotropic:kc},
       the optimal scaling coefficient $\tau$ 
       is proportional \flo{to} the inverse of the velocity, that is,
       the slowness. 
       It seems to depend only slightly on the type of waves.
       For P-planewave, the value of the optimal scaling factor 
       is hard to identify \flo{and lies} between $\tau=1/\cp$ and $\tau=1/\cs$. 
       The S-planewave shows more flexibility, with the optimal value 
       for the scaling given in range $\tau \in (1/\cp, 1/\cs)$.
       % There is more flexibility for S-planewave, where it 
       % seems slightly better to use $\tau=1/\cs$ even though 
       % the error level remains relatively similar in the interval
       % $\tau \in (1/\cp, 1/\cs)$.

 \item For the Godunov stabilization, \cref{figure:3D-tau-amplitude:isotropic:g},
       the optimal scaling coefficient does not depend on the type of wave\flo{s}, 
       and is simply given by $\tau=1$.
 
% \item For the P-planewave, the accuracy for the different fields is relatively 
%       similar (\cref{figure:3D-tau-amplitude:iso02_qP-Su}), except for the 
%       component $\sigma_{xz}$ which is slightly less accurate (\num{6e-6} 
%       relative error compared to \num{e-6} for the other components).
%       In this experiment, the S-planewave is slightly less accurate than 
%       the P-planewave because it has a smaller wavelength (\num{0.25} compared 
%       to \num{0.31}) 
%       at the selected frequencies.

% \item For the P-planewave (\cref{figure:3D-tau-amplitude:iso02_qP-Su}),
%       the Godunov stabilization $\tauSuMGa$ gives the most accurate results, 
%       accuracy that does not seem to be reached by identity-based stabilization, 
%       nor by the Kelvin--Christoffel coefficient. 
%       The best performance with identity-based stabilization is 
%       obtained with scale $\rho\omega \cp$, indicated in vertical line.

% \item For the S-planewave (\cref{figure:3D-tau-amplitude:iso02_sH-Su}), 
%       the accuracy obtained with the Godunov stabilization can be 
%       improved with an identity-based stabilization (from $\sim$\num{e-5} 
%       to $\sim$\num{5e-6}). 
%       Here, the optimal scale seems close to $\rho \omega  \cs$.
 \item Comparing the accuracy of the different families of stabilization, 
       the Godunov stabilization (with $\tau=1$) and the Kelvin-Christoffel  
       one (with appropriate slowness scaling) give similar accuracy. 
       The identity-based stabilization, although close, does 
       not match this accuracy.

\end{itemize}

Therefore, the Godunov stabilization appears as the most versatile choice 
as it accurately treats both the P-planewave and the S-planewave, 
without any specific scaling coefficient, with $\tau=1$. 
% On the other hand, the Kelvin--Christoffel and 
% identity-based stabilizations require a correct 
% scaling factor in order to reach the same level 
% of performance as the Godunov stabilization. 
For the identity-based, the optimal scaling factor 
clearly depends on the types of waves. 
For the Kelvin--Christoffel stabilization, the same value 
of \flo{the} scaling factor can be chosen for both types of waves. 
This can be expected since the Kelvin--Christoffel matrix 
contains the information of the material. 
% Nevertheless, the need for scaling factor tuning can be problematic 
% in practice when simulating waves in heterogeneous media where different 
% types of waves propagate, this need is avoided with the Godunov stabilization.

% ------------------------------------------------------------------------------- %
% \renewcommand{\experiment}{\texttt{benchmark3Dvti-Pw}}
% \subsubsection{Elastic VTI (\experiment)}
% \renewcommand{\experiment}{\texttt{Configuration~PwVTI}}
\subsubsection{Elastic vertical transverse isotropic medium}
\label{subsection:aniso-3D-planewaves_01}
% anisotropy
%\renewcommand{\plottauSuIda}{tau9}
%\renewcommand{\plottauSuIdb}{tau5}
%\renewcommand{\plottauSuKCa}{tau17}
%\renewcommand{\plottauSuKCb}{tau15}
%\renewcommand{\plottauSuMGa}{tau23}
%\renewcommand{\plottauSuMGb}{tau33}
% ------------------------------------------------------------------------------- %

\begin{figure}[ht!] \centering
\subfloat[][Field $u_x$ for the qP-planewave at 8 \si{\milli\Hz}.]
         {\includegraphics[scale=.80]{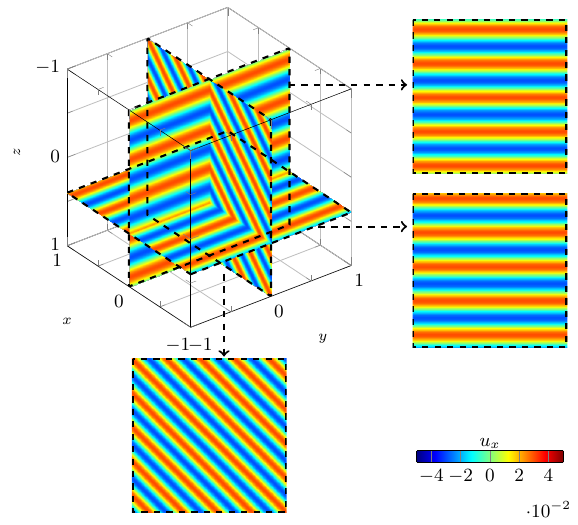}}
\subfloat[][Field $u_y$ for the sH-planewave at \num{0.4} \si{\milli\Hz}.]
         {\includegraphics[scale=.80]{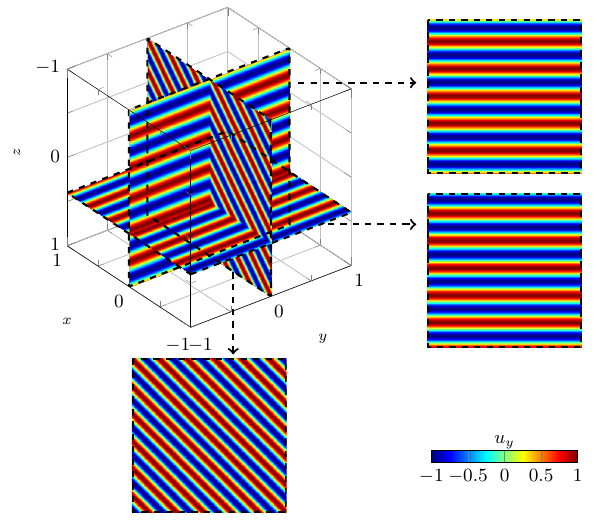}}
\caption{Reference solution for the 3D planewave propagation in VTI medium~\cref{benchmarkvti03}.
         %using parameters~\cref{benchmarkvti03}.
         }
\label{figure:3Dvti-pw:wavefield}
\end{figure}

We modify the previous experiment to include anisotropy in the medium.
We work with Thomsen's parameters, \cite{thomsen1986weak}, with notation
of \cite[Section 3.5.2]{Pham2023hdgRR} for VTI, and 
select\footnote{The values of the anisotropic 
coefficients $\epsilon$, $\delta$ and $\gamma$ correspond 
to muscovite crystal in \cite{thomsen1986weak}. The 
values of $\lambda^{\mathrm{TI}}$ and $\mu^{\mathrm{TI}}$ 
are selected such that the qP- and sH-wavespeeds 
in \cref{eq:numerics:3dvti:wavespeeds} correspond to the 
P- and S-wavespeed of the isotropic experiment~\cref{benchmark02}, i.e., 
\num{2.5e-3} and \num{e-4} \si{\meter\per\second} respectively.}
%\begin{equation}\label{benchmarkvti03}
%% \experiment\qquad 
%\left\lbrace
%\begin{aligned}
%  & \lambda^{\mathrm{TI}}  \,=\, \num{3.6107e-6}\si{\pascal} \,, \quad
%    \mu^{\mathrm{TI}}      \,=\, \num{3.0490e-9}\si{\pascal} \,, \quad
%    \rho     \,=\, \num{1}\si{\kg\per\meter\cubed} \,, \\
%  & \epsilon \,=\, \num{1.12}\,, \quad
%    \delta   \,=\, \num{-0.235} \, \quad
%    \gamma   \,=\, \num{2.28} \,.\\
%  &  \text{qP-planewave propagation:} \qquad 
% \displacement_{\mathrm{pw}}(\bx) \,:=\, \displacement_{\mathrm{pw}}^\mathrm{qP}(\bx)\,,\\
%  &  \text{sH-planewave propagation:} \qquad 
% \displacement_{\mathrm{pw}}(\bx) \,:=\, \displacement_{\mathrm{pw}}^\mathrm{sH}(\bx)\,.
%\end{aligned} \right. \end{equation}
\begin{equation}\label{benchmarkvti03}
\begin{aligned}
  & \lambda^{\mathrm{TI}}  \,=\, \num{3.6107e-6}\si{\pascal} \,, \quad
    \mu^{\mathrm{TI}}      \,=\, \num{3.0490e-9}\si{\pascal} \,, \quad
    \rho     \,=\, \num{1}\si{\kg\per\meter\cubed} \,, \\
  & \epsilon \,=\, \num{1.12}\,, \quad
    \delta   \,=\, \num{-0.235} \, \quad
    \gamma   \,=\, \num{2.28} \,.
\end{aligned} \end{equation}
This amounts to the following values of the stiffness tensor coefficients:
% given in \si{\pascal} (see \cref{remark:scaling} for discussion on the scales):
%\begin{subequations}\begin{align}
%& C_{11} \,=\, C_{22} \,=\, \num{1.1719e-5} \,, \quad
%  C_{33} \,=\,              \num{3.6168e-6}  \,, \\
%& C_{44} \,=\, C_{55} \,=\, \num{3.0490e-9}      \,, \quad
%  C_{66} \,=\,              \num{1.6953e-8}  \,, \\
%& C_{12} \,=\, C_{11} - 2 C_{66} \,=\,\num{1.1685e-05} \,, \quad
%  C_{23} \,=\, C_{23} \,=\, \num{2.6268e-06}\,.
%\end{align}\end{subequations}
\begin{subequations}\begin{align}
& C_{11} \,=\, C_{22} \,=\, \num{1.1719e-5}  \,, \quad
  C_{33} \,=\,              \num{3.6168e-6}  \,, \quad
  C_{44} \,=\, C_{55} \,=\, \num{3.0490e-9}  \,, \nonumber \\ 
& C_{66} \,=\,              \num{1.6953e-8}  \,, \quad
  C_{12} \,=\, C_{11} - 2 C_{66} \,=\,\num{1.1685e-05} \,, \quad 
  C_{23} \,=\, C_{23} \,=\, \num{2.6268e-06}\,. \nonumber
\end{align}\end{subequations}
% In vertical transverse isotropic medium, 
\flo{In direction $\mathbf{d}=(d_x,0,d_z)$, 
the qP- and sH-planewaves are given by, \cite[Appendix~A]{Pham2023hdgRR},}
\begin{subequations}\begin{align}
\displacement_{\mathrm{pw}}^\mathrm{qP}(\bx)&\,=\,
\dfrac{e^{\ii\frac{\omega}{\cqp}\, (\mathbf{d} \, \cdot \, \bx)}}
      {\sqrt{(C_{11}+C_{55}) d_x^2 + (C_{55}+C_{33}) d_z^2 - 2\rho \cqp}}
\begin{pmatrix}
 \sqrt{C_{55} d_x^2 \,+\, C_{33} d_z^2 - \rho \cqp^2} \\
 0 \\
 \sqrt{C_{11} d_x^2 \,+\, C_{55} d_z^2 - \rho \cqp^2} 
\end{pmatrix} \,, \\[.50em]
\displacement_{\mathrm{pw}}^\mathrm{sH}(\bx)&\,=\,
(0, 1, 0)^t \,\,
e^{\ii\frac{\omega}{\csh}\, (\mathbf{d} \, \cdot \, \bx)} \, ,
\end{align}\end{subequations}
with corresponding qP- and sH-wavespeeds, denoted by $\cqp$ 
and $\csh$,  %given by \cite[Appendix~A]{Pham2023hdgRR},
\begin{subequations}\label{eq:numerics:3dvti:wavespeeds}
\begin{align}
  \cqp^2(d_x,0,d_z) & \,=\, 
    \dfrac{1}{2\rho}\bigg( C_{11} d_x^2 + C_{33} d_z^2 + C_{55} \\
  & \,+\, 
    \sqrt{\big( (C_{11}-C_{55})d_x^2 + (C_{55}-C_{33})d_z^2\big)^2   
           + 4 (C_{13} + C_{55})^2 \, d_x^2\, d_z^2}
      \,\, \bigg) \nonumber
  \,, \\
  \csh^2(d_x,0,d_z) & \,=\, \dfrac{C_{66} d_x^2 \,+\, C_{55} d_z^2}{\rho}\,.
\end{align}\end{subequations}
In the experiment, we \flo{work with} direction 
$\mathbf{d}=\big(1/\sqrt{2},0,1/\sqrt{2}\big)^t$,
and illustrate the corresponding analytic solutions 
in \cref{figure:3Dvti-pw:wavefield}.

% ============================
% ============================

\begin{figure}[ht!] \centering
  \renewcommand{\datafile}{figures/numerics/3D-planewave_optimal-stabilization/data/data_elastic-vti-version04_TauId.txt}
  \renewcommand{\myxlabel}{$\tau$}
  \renewcommand{\myylabel}{$\err$}

  \pgfkeys{/pgf/fpu=true}
  \renewcommand{\qPomega}{0.00250} % Rho Cp only
  \renewcommand{\qSomega}{0.00010} % Rho Cs only
  \pgfmathsetmacro{\scaleomega}{0.0503}
  \pgfmathsetmacro{\xmin}{0.60e-8} \pgfmathsetmacro{\xmax}{0.40e2} 
  \pgfmathsetmacro{\ymin}{4e-7}    \pgfmathsetmacro{\ymax}{2e-2} 
  \pgfkeys{/pgf/fpu=false}
  \renewcommand{\dataA}{qP-double_ux_freq8.00mHz} 
  \subfloat[][{Identity-based stabilization using real and imaginary scaling factors. 
               Relative error for P-planewave $\err[u_x]$ (left) and $\err[\sigma_{zz}]$ (middle);
               for S-planewave $\err[u_y]$ (right).}]{
              \renewcommand{\myylabel}{${\err[u_x]}$}
              %% --------------------------------------------------------------------------------------
\begin{tikzpicture}
%% --------------------------------------------------------------------------------------
\begin{axis}[
             enlargelimits=false, 
             xlabel={\scriptsize \myxlabel },
             ylabel={\scriptsize \myylabel },
             enlarge y limits=false,
             enlarge x limits=false,
             xmin=\xmin,xmax=\xmax, 
             ymin=\ymin,ymax=\ymax,             
             yminorticks=true,
             height=\plotheight,width=\plotwidth,
             scale only axis,
             ylabel style={rotate=-90},
             ylabel style = {yshift=4.0em, xshift=2.50em},
             xlabel style = {yshift = 1em, xshift=-5.50em},
             clip mode=individual,
             ytick={1e-6,1e-5,1e-4,1e-3},
             xtick={1e-8,1e-7,1e-6,1e-5,1e-4,1e-3,1e-2,1e-1,1e0,1e1,1e2,1e3},
             xticklabels={,,\num{e-6},,,,,\num{e-1},,\num{e1},,},
             minor grid style={line width=.5pt, draw=gray!50, densely dotted},
             major grid style={line width=.5pt, draw=gray!50, densely dashed},
             xmajorgrids=true,
           % xminorgrids=true,
             yminorgrids=true,
             ymajorgrids=true,             ymode=log,
             xmode=log,
             %y  tick label style={/pgf/number format/fixed},
             %y  tick label style={/pgf/number format/fixed zerofill},
             %y  tick label style={/pgf/number format/precision=1},
             %y  tick label style={/pgf/number format/precision=2},
             % y  tick label style={/pgf/number format/sci},
             %y tick label style={/pgf/number format/.cd},
             label style={font=\scriptsize},
             tick label style={font=\scriptsize},
             legend style={font=\scriptsize\selectfont},
             legend style={at={(0.99,1.17)}},
             legend columns=4,
             ]  %% foreground marks

%     \pgfmathsetmacro{\scale}{1.0}
     \addplot[color=black,densely dotted, line width=0.5, mark=star,% only marks,
              mark size=2.50,mark options={solid,fill=black, line width=.50}]
              table[x expr = \thisrow{tau}/\scaleomega, 
                    y expr = \thisrow{\dataA_plusRe}]
              {\datafile}; % \addlegendentry{$\tau \identity\,\,\,$}
     \addplot[color=\myblue,densely dotted, line width=0.5, mark=+, % only marks,
              mark size=2.50,mark options={solid,fill=\myblue, line width=.50}]
              table[x expr = \thisrow{tau}/\scaleomega, 
                    y expr = \thisrow{\dataA_minusRe}]
              {\datafile}; % \addlegendentry{$-\tau \identity \,\,\,$}

     \addplot[color=red,densely dotted, line width=0.5, mark=o, % only marks,
              mark size=2.50,mark options={solid,fill=red, line width=.50}]
              table[x expr = \thisrow{tau}/\scaleomega, 
                    y expr = \thisrow{\dataA_plusIm}]
              {\datafile}; % \addlegendentry{$\ii\tau \identity\,\,\,$}
     \addplot[color=orange,densely dotted, line width=0.5, mark=x, % only marks,
              mark size=2.50,mark options={solid,fill=orange, line width=.50}]
              table[x expr = \thisrow{tau}/\scaleomega, 
                    y expr = \thisrow{\dataA_minusIm}]
              {\datafile}; % \addlegendentry{$-\ii\tau \identity\,\,\,$}

     \draw[densely dotted, \mygreen, line width=1] (\qPomega,\ymin) -- (\qPomega,\ymax);
     \draw[densely dotted, \mygreen, line width=1] (\qSomega,\ymin) -- (\qSomega,\ymax);
     \node[\mygreen, line width=1,yshift=0em,anchor=north,align=center] at (\qPomega,\ymin)  
          {\scriptsize $\rho \cqp$};
     \node[\mygreen, line width=1,yshift=0em,anchor=north,align=center] at (\qSomega,\ymin)  
          {\scriptsize $\rho \csh$};
\end{axis}
\end{tikzpicture}
              \renewcommand{\dataA}{qP-double_szz_freq8.00mHz} \hspace*{-3.50em}
              \renewcommand{\myylabel}{${\err[\sigma_{zz}]}$}
              %% --------------------------------------------------------------------------------------
\begin{tikzpicture}
%% --------------------------------------------------------------------------------------
\begin{axis}[
             enlargelimits=false, 
             xlabel={\scriptsize \myxlabel },
             enlarge y limits=false,
             enlarge x limits=false,
             xmin=\xmin,xmax=\xmax, 
             ymin=\ymin,ymax=\ymax,             
             yminorticks=true,
             height=\plotheight,width=\plotwidth,
             scale only axis,
             ylabel style={rotate=-90},
             ylabel style = {yshift=4.0em, xshift=1.00em},
             xlabel style = {yshift = 1em, xshift=-5.50em},
             clip mode=individual,
             ytick={1e-6,1e-5,1e-4,1e-3,1e-2},
             xtick={1e-8,1e-7,1e-6,1e-5,1e-4,1e-3,1e-2,1e-1,1e0,1e1,1e2,1e3},
             xticklabels={,,\num{e-6},,,,,\num{e-1},,\num{e1},,},
             yticklabels={,,},
             minor grid style={line width=.5pt, draw=gray!50, densely dotted},
             major grid style={line width=.5pt, draw=gray!50, densely dashed},
             xmajorgrids=true,
           % xminorgrids=true,
             yminorgrids=true,
             ymajorgrids=true,             ymode=log,
             xmode=log,
             %y  tick label style={/pgf/number format/fixed},
             %y  tick label style={/pgf/number format/fixed zerofill},
             %y  tick label style={/pgf/number format/precision=1},
             %y  tick label style={/pgf/number format/precision=2},
             % y  tick label style={/pgf/number format/sci},
             %y tick label style={/pgf/number format/.cd},
             label style={font=\scriptsize},
             tick label style={font=\scriptsize},
             legend style={font=\scriptsize\selectfont},
             legend style={at={(1.15,1.17)}},
             legend columns=4,
             ]  %% foreground marks

%     \pgfmathsetmacro{\scale}{1.0}
     \addplot[color=black,densely dotted, line width=0.5, mark=star, % only marks,
              mark size=2.50,mark options={solid,fill=black, line width=.50}]
              table[x expr = \thisrow{tau}/\scaleomega, 
                    y expr = \thisrow{\dataA_plusRe}]
              {\datafile}; \addlegendentry{$\tauSuIdRep$}
     \addplot[color=\myblue,densely dotted, line width=0.5, mark=+, % only marks,
              mark size=2.50,mark options={solid,fill=\myblue, line width=.50}]
              table[x expr = \thisrow{tau}/\scaleomega, 
                    y expr = \thisrow{\dataA_minusRe}]
              {\datafile}; \addlegendentry{$\tauSuIdRem$}

     \addplot[color=red,densely dotted, line width=0.5, mark=o, % only marks,
              mark size=2.50,mark options={solid,fill=red, line width=.50}]
              table[x expr = \thisrow{tau}/\scaleomega, 
                    y expr = \thisrow{\dataA_plusIm}]
              {\datafile}; \addlegendentry{$\tauSuIdImp$}
     \addplot[color=orange,densely dotted, line width=0.5, mark=x, % only marks,
              mark size=2.50,mark options={solid,fill=orange, line width=.50}]
              table[x expr = \thisrow{tau}/\scaleomega, 
                    y expr = \thisrow{\dataA_minusIm}]
              {\datafile}; \addlegendentry{$\tauSuIdImm$}

     \node[yshift=0em,anchor=north west,align=center,fill=white] at (\xmin,\ymax)  
          {\scriptsize \myylabel};

     \draw[densely dotted, \mygreen, line width=1] (\qPomega,\ymin) -- (\qPomega,\ymax);
     \draw[densely dotted, \mygreen, line width=1] (\qSomega,\ymin) -- (\qSomega,\ymax);
     \node[\mygreen, line width=1,yshift=0em,anchor=north,align=center] at (\qPomega,\ymin)  
          {\scriptsize $\rho \cqp$};
     \node[\mygreen, line width=1,yshift=0em,anchor=north,align=center] at (\qSomega,\ymin)  
          {\scriptsize $\rho \csh$};
\end{axis}
\end{tikzpicture}
              \hspace*{-5.0em}
              \pgfkeys{/pgf/fpu=true}
              \pgfmathsetmacro{\xmin}{0.60e-7} \pgfmathsetmacro{\xmax}{0.60e2} 
              \pgfmathsetmacro{\scaleomega}{0.0025}
              \pgfkeys{/pgf/fpu=false}
              \renewcommand{\dataA}{sH-double_uy_freq0.40mHz} 
              \renewcommand{\myylabel}{${\err[u_y]}$}
              %% --------------------------------------------------------------------------------------
\begin{tikzpicture}
%% --------------------------------------------------------------------------------------
\begin{axis}[
             enlargelimits=false, 
             xlabel={\scriptsize \myxlabel },
             enlarge y limits=false,
             enlarge x limits=false,
             xmin=\xmin,xmax=\xmax, 
             ymin=\ymin,ymax=\ymax,             
             yminorticks=true,
             height=\plotheight,width=\plotwidth,
             scale only axis,
             ylabel style={rotate=-90},
             ylabel style = {yshift=4.0em, xshift=1.00em},
             xlabel style = {yshift =1em, xshift=-5.50em},
             clip mode=individual,
             ytick={1e-6,1e-5,1e-4,1e-3,1e-2},
             xtick={1e-8,1e-7,1e-6,1e-5,1e-4,1e-3,1e-2,1e-1,1e0,1e1,1e2,1e3,1e4,1e5,1e6,1e7},
             xticklabels={,,\num{e-6},,,,,\num{e-1},,\num{e1},},
             yticklabels={,,,,},
             minor grid style={line width=.5pt, draw=gray!50, densely dotted},
             major grid style={line width=.5pt, draw=gray!50, densely dashed},
             xmajorgrids=true,
           % xminorgrids=true,
             yminorgrids=true,
             ymajorgrids=true,             ymode=log,
             xmode=log,
             %y  tick label style={/pgf/number format/fixed},
             %y  tick label style={/pgf/number format/fixed zerofill},
             %y  tick label style={/pgf/number format/precision=1},
             %y  tick label style={/pgf/number format/precision=2},
             % y  tick label style={/pgf/number format/sci},
             %y tick label style={/pgf/number format/.cd},
             label style={font=\scriptsize},
             tick label style={font=\scriptsize},
             legend style={font=\scriptsize\selectfont},
             legend style={at={(0.99,1.17)}},
             legend columns=4,
             ]  %% foreground marks

%     \pgfmathsetmacro{\scale}{1.0}
     \addplot[color=black,densely dotted, line width=0.5, mark=star, %only marks,
              mark size=2.50,mark options={solid,fill=black, line width=.50}]
              table[x expr = \thisrow{tau}/\scaleomega, 
                    y expr = \thisrow{\dataA_plusRe}]
              {\datafile}; %\addlegendentry{$\tau \identity\,\,\,$}
     \addplot[color=\myblue,densely dotted, line width=0.5, mark=+, %only marks,
              mark size=2.50,mark options={solid,fill=\myblue, line width=.50}]
              table[x expr = \thisrow{tau}/\scaleomega, 
                    y expr = \thisrow{\dataA_minusRe}]
              {\datafile}; %\addlegendentry{$-\tau \identity \,\,\,$}

     \addplot[color=red,densely dotted, line width=0.5, mark=o, %only marks,
              mark size=2.50,mark options={solid,fill=red, line width=.50}]
              table[x expr = \thisrow{tau}/\scaleomega, 
                    y expr = \thisrow{\dataA_plusIm}]
              {\datafile}; %\addlegendentry{$\ii\tau \identity\,\,\,$}
     \addplot[color=orange,densely dotted, line width=0.5, mark=x, %only marks,
              mark size=2.50,mark options={solid,fill=orange, line width=.50}]
              table[x expr = \thisrow{tau}/\scaleomega, 
                    y expr = \thisrow{\dataA_minusIm}]
              {\datafile}; %\addlegendentry{$-\ii\tau \identity\,\,\,$}

     \node[yshift=0em,anchor=north west,align=center,fill=white] at (\xmin,\ymax)  
          {\scriptsize \myylabel};

     \draw[densely dotted, \mygreen, line width=1] (\qPomega,\ymin) -- (\qPomega,\ymax);
     \draw[densely dotted, \mygreen, line width=1] (\qSomega,\ymin) -- (\qSomega,\ymax);
     \node[\mygreen, line width=1,yshift=0em,anchor=north,align=center] at (\qPomega,\ymin)  
          {\scriptsize $\rho \cqp$};
     \node[\mygreen, line width=1,yshift=0em,anchor=north,align=center] at (\qSomega,\ymin)  
          {\scriptsize $\rho \csh$};
\end{axis}
\end{tikzpicture}
              \label{figure:3D-tau-amplitude:vti:id}}

  %% =======================================================================
  \setlength{\plotheight}{3.00cm}
  \pgfkeys{/pgf/fpu=true}
  \pgfmathsetmacro{\xmin}{0.60e-4} \pgfmathsetmacro{\xmax}{0.40e5} 
  \pgfmathsetmacro{\ymin}{4e-7} \pgfmathsetmacro{\ymax}{2e-4} 
  \pgfkeys{/pgf/fpu=false}
  \renewcommand{\dataA}{qP-double_ux_freq8.00mHz}
  \renewcommand{\dataB}{sH-double_uy_freq0.40mHz}
  \renewcommand{\dataC}{qP-double_szz_freq8.00mHz}
  \renewcommand{\datafile}
               {figures/numerics/3D-planewave_optimal-stabilization/data/data_elastic-vti-version04_TauGodunov.txt}
  \renewcommand{\myxlabel}{$\tau$}
  \subfloat[][Using Godunov stabilization.]
             {%% --------------------------------------------------------------------------------------
\begin{tikzpicture}
%% --------------------------------------------------------------------------------------
\begin{axis}[
             enlargelimits=false, 
             xlabel={\scriptsize \myxlabel },
             ylabel={\scriptsize \myylabel },
             enlarge y limits=false,
             enlarge x limits=false,
             xmin=\xmin,xmax=\xmax, 
             ymin=\ymin,ymax=\ymax,             
             yminorticks=true,
             height=\plotheight,width=\plotwidth,
             scale only axis,
             ylabel style={rotate=-90},
             ylabel style = {yshift=3.5em, xshift=2em},
             xlabel style = {yshift=1em, xshift=-5.50em},
             clip mode=individual,
             ytick={1e-6,1e-5,1e-4,1e-3,1e-2,1e-1},
             xtick={1e-5,1e-4,1e-3,1e-2,1e-1,1e0,1e1,1e2,1e3,1e4,1e5},
             xticklabels={,,,\num{e-2},,\num{1},,\num{e2},,\num{e4},},
             minor grid style={line width=.5pt, draw=gray!50, densely dotted},
             major grid style={line width=.5pt, draw=gray!50, densely dashed},
             xmajorgrids=true,
           % xminorgrids=true,
             yminorgrids=true,
             ymajorgrids=true,             ymode=log,
             xmode=log,
             %y  tick label style={/pgf/number format/fixed},
             %y  tick label style={/pgf/number format/fixed zerofill},
             %y  tick label style={/pgf/number format/precision=1},
             %y  tick label style={/pgf/number format/precision=2},
             % y  tick label style={/pgf/number format/sci},
             %y tick label style={/pgf/number format/.cd},
             label style={font=\scriptsize},
             tick label style={font=\scriptsize},
             legend style={font=\scriptsize\selectfont},
             legend style={at={(1.900,0.60)}},
             legend columns=1,
             ]  %% foreground marks

%     \pgfmathsetmacro{\scale}{1.0}
%     \addplot[color=black,densely dashed, line width=0.2, mark=star, only marks,
%              mark size=2.50,mark options={solid,fill=black, line width=.50}]
%              table[x expr = \thisrow{tau}, 
%                    y expr = \thisrow{\dataA}]
%              {\datafile}; \addlegendentry{P-planewave}
%     \addplot[color=\myblue,densely dashed, line width=0.2, mark=+, only marks,
%              mark size=2.50,mark options={solid,fill=\myblue, line width=.50}]
%              table[x expr = \thisrow{tau}, 
%                    y expr = \thisrow{\dataB}]
%              {\datafile}; \addlegendentry{S-planewave}

     \addplot[color=red,densely dotted, line width=0.5, mark=o, % only marks,
              mark size=2.50,mark options={solid,fill=red, line width=.50}]
              table[x expr = \thisrow{tau}, 
                    y expr = \thisrow{\dataA}]
              {\datafile}; \addlegendentry{$\err[u_x]$ for P-planewave}
     \addplot[color=cyan,densely dotted, line width=0.5, mark=+,% only marks,
              mark size=2.75,mark options={solid,fill=cyan, line width=1}]
              table[x expr = \thisrow{tau}, 
                    y expr = \thisrow{\dataC}]
              {\datafile}; \addlegendentry{$\err[\sigma_{zz}]$ for P-planewave}
     \addplot[color=orange,densely dotted, line width=0.5, mark=x,% only marks,
              mark size=2.75,mark options={solid,fill=orange, line width=1}]
              table[x expr = \thisrow{tau}, 
                    y expr = \thisrow{\dataB}]
              {\datafile}; \addlegendentry{$\err[u_y]$ for S-planewave}
              
%     \draw[densely dotted, magenta, line width=1] (\qPomega,\ymin) -- (\qPomega,\ymax);
%     \draw[densely dotted, magenta, line width=1] (\qSomega,\ymin) -- (\qSomega,\ymax);
%     \node[magenta, line width=1,yshift=0em,anchor=north,align=center] at (\qPomega,\ymin)  
%          {\scriptsize $\rho \omega c_p$};
%     \node[magenta, line width=1,yshift=0em,anchor=north,align=center] at (\qSomega,\ymin)  
%          {\scriptsize $\rho \omega c_s$};
\end{axis}
\end{tikzpicture}
              \label{figure:3D-tau-amplitude:vti:g}} 
             \hspace*{-4.50em}
  \pgfkeys{/pgf/fpu=true}
  \pgfmathsetmacro{\xmin}{0.60e-3} \pgfmathsetmacro{\xmax}{0.40e8} 
  \renewcommand{\qPomega}{400} % 1/(Rho Cp)
  \renewcommand{\qSomega}{1e4} % 1/(Rho Cs)
  \pgfkeys{/pgf/fpu=false}
  \renewcommand{\datafile}
               {figures/numerics/3D-planewave_optimal-stabilization/data/data_elastic-vti-version04_TauKC.txt}
  \subfloat[][Using Kelvin--Christoffel stabilization.]
             {\makebox[.40\linewidth][c]{{
              %% --------------------------------------------------------------------------------------
\begin{tikzpicture}
%% --------------------------------------------------------------------------------------
\begin{axis}[
             enlargelimits=false, 
             xlabel={\scriptsize \myxlabel },
             ylabel={\scriptsize \myylabel },
             enlarge y limits=false,
             enlarge x limits=false,
             xmin=\xmin,xmax=\xmax, 
             ymin=\ymin,ymax=\ymax,             
             yminorticks=true,
             height=\plotheight,width=\plotwidth,
             scale only axis,
             ylabel style={rotate=-90},
             ylabel style = {yshift=3.5em, xshift=2em},
             xlabel style = {yshift =1em, xshift=-5.50em},
             clip mode=individual,
             ytick={1e-6,1e-5,1e-4,1e-3,1e-2,1e-1},
             xtick={1e-5,1e-4,1e-3,1e-2,1e-1,1e0,1e1,1e2,1e3,1e4,1e5,1e6,1e7,1e8,1e9},
             xticklabels={,,,,\num{e-1},,\num{e1},,,,,\num{e6},,\num{e8},},
             minor grid style={line width=.5pt, draw=gray!50, densely dotted},
             major grid style={line width=.5pt, draw=gray!50, densely dashed},
             xmajorgrids=true,
           % xminorgrids=true,
             yminorgrids=true,
             ymajorgrids=true,             ymode=log,
             xmode=log,
             %y  tick label style={/pgf/number format/fixed},
             %y  tick label style={/pgf/number format/fixed zerofill},
             %y  tick label style={/pgf/number format/precision=1},
             %y  tick label style={/pgf/number format/precision=2},
             % y  tick label style={/pgf/number format/sci},
             %y tick label style={/pgf/number format/.cd},
             label style={font=\scriptsize},
             tick label style={font=\scriptsize},
             legend style={font=\scriptsize\selectfont},
             legend style={at={(0.99,1.17)}},
             legend columns=4,
             ]  %% foreground marks

%     \pgfmathsetmacro{\scale}{1.0}
%     \addplot[color=black,densely dashed, line width=0.2, mark=star, only marks,
%              mark size=2.50,mark options={solid,fill=black, line width=.50}]
%              table[x expr = \thisrow{tau}, 
%                    y expr = \thisrow{\dataA}]
%              {\datafile}; \addlegendentry{P-planewave}
%     \addplot[color=\myblue,densely dashed, line width=0.2, mark=+, only marks,
%              mark size=2.50,mark options={solid,fill=\myblue, line width=.50}]
%              table[x expr = \thisrow{tau}, 
%                    y expr = \thisrow{\dataB}]
%              {\datafile}; \addlegendentry{S-planewave}

     \addplot[color=red,densely dotted, line width=0.5, mark=o, % only marks,
              mark size=2.50,mark options={solid,fill=red, line width=.50}]
              table[x expr = \thisrow{tau}, 
                    y expr = \thisrow{\dataA}]
              {\datafile}; % \addlegendentry{}
     \addplot[color=cyan,densely dotted, line width=0.5, mark=+,% only marks,
              mark size=2.75,mark options={solid,fill=cyan, line width=1}]
              table[x expr = \thisrow{tau}, 
                    y expr = \thisrow{\dataC}]
              {\datafile}; 
     \addplot[color=orange,densely dotted, line width=0.5, mark=x,% only marks,
              mark size=2.75,mark options={solid,fill=orange, line width=1}]
              table[x expr = \thisrow{tau}, 
                    y expr = \thisrow{\dataB}]
              {\datafile}; % \addlegendentry{}

     \draw[densely dotted, \mygreen, line width=1] (\qPomega,\ymin) -- (\qPomega,\ymax);
     \draw[densely dotted, \mygreen, line width=1] (\qSomega,\ymin) -- (\qSomega,\ymax);
     \node[\mygreen, line width=1,yshift=0em,anchor=north,align=center] at (\qPomega,\ymin)  
          {\scriptsize $\tfrac{1}{\cqp}$};
     \node[\mygreen, line width=1,yshift=0em,anchor=north,align=center] at (\qSomega,\ymin)  
          {\scriptsize $\tfrac{1}{\csh}$};
\end{axis}
\end{tikzpicture}}}
              \label{figure:3D-tau-amplitude:vti:kc}}

  \caption{Investigation of the optimal scaling $\tau$ within the three
           families of stabilization in \cref{eq:numerics:stabilization-u:3Dpw} 
           for planewave propagation in VTI medium with parameters~\cref{benchmarkvti03}.
           Our criterion of accuracy is the relative error 
           $\err$ \cref{eq:relative-error_field}.
           }
  \label{figure:3D-tau-amplitude:vti}
\end{figure}
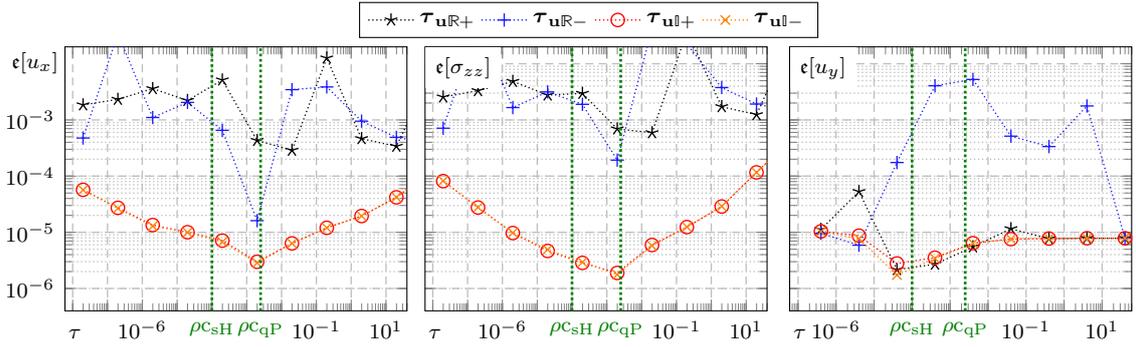
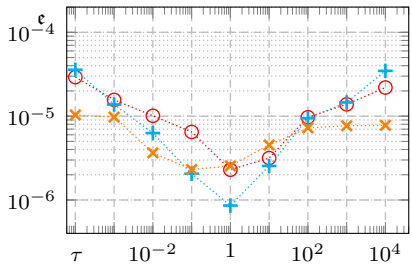
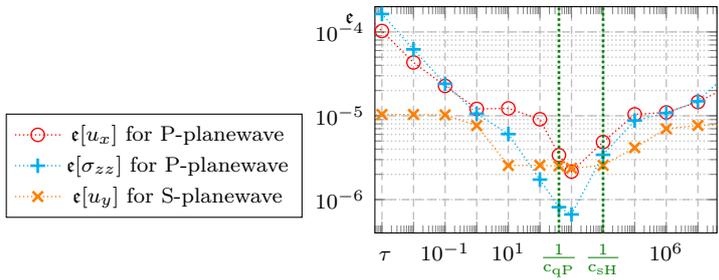

In \cref{figure:3D-tau-amplitude:vti}, we evaluate the accuracy of 
the solutions depending on the scaling of the coefficient $\tau$ for 
the three families of stabilization in \cref{eq:numerics:stabilization-u:3Dpw}. 
% for the HDG formulation $\formulationU$. 
The qP-planewave is computed at frequency \num{8} \si{\milli\Hz} and
the sH-planewave at frequency \num{0.4} \si{\milli\Hz}. 
The observations for the elastic VTI medium are close to those 
of the elastic isotropic case of \cref{subsection:iso-3D-planewaves}:
\begin{itemize}[leftmargin=*] \setlength{\itemsep}{-2pt}
  \item The optimal value of the scaling coefficient for the Godunov stabilization 
        is $\tau=1$, and it does not depend on the type of wave propagating.
  \item The optimal value of the scaling for the Kelvin--Christoffel stabilization 
        is in between $1/\cqp$ and $1/\csh$, with the value of the optimal scaling
        hard to find for the qP-planewave.
  \item For \flo{the} identity-based stabilization, the optimal coefficient 
        (which is imaginary) depends on the type of waves and is 
        given by the impedances. 
        Furthermore, the accuracy of the results, even with the optimal 
        scaling, is less than the accuracy obtained using other stabilizations.
\end{itemize}
In short, the Godunov stabilization is the most robust to treat both 
types of waves, with its default scaling factor, $\tau=1$, being optimal.

% ---------------------------------------------------------------------- %
% ---------------------------------------------------------------------- %

% ------------------------------------------------------------------------------- %
\subsection{2D experiments with highly varying properties and point-source}
\label{subsection:numerics-2d}
% ------------------------------------------------------------------------------- %

To \flo{study further} the efficiency of stabilizations,
we consider a heterogeneous medium on a two-dimensional 
disk domain $\Omega$.
The variation of the properties is too strong for \flo{a} 
piecewise-constant representation to be reliable and we 
instead use piecewise-polynomial representations, 
\modif{whose implementation is straightforward with the formulation based 
on the compliance tensor, \cref{1stsystem_variants:u}.
In this experiment, the piecewise polynomial representation is 
independent on each cell, and allows
discontinuities across the interfaces.}
We consider \flo{a problem with a Dirac point-source},
%\begin{subequations} \label{eq:numerics:eq-elastic-2d}
%\begin{empheq}[left = \empheqlbrace\,]{align}
%  & - \omega^2 \,\rho \,\displacement  \, - \,  \nabla \cdot \tensS  \,=\, \delta_\by \,,\qquad \text{in $\Omega$,} \\
%  &   \tensCS \tensS \,=\, \nabla^\mathfrak{s}\displacement \,,\hspace*{7.3em} \text{in $\Omega$,} \\
%  &   \tensS\,\n \,=\, -\ii \omega \mathcal{Z} \displacement \,, \hspace*{6em}\text{on $\partial\Omega$} \,,
%\end{empheq}
%\end{subequations}
\begin{subequations} \label{eq:numerics:eq-elastic-2d}
\begin{empheq}[left = \empheqlbrace\,]{align}
  & - \omega^2 \,\rho \,\displacement  \, - \,  \nabla \cdot \tensS  \,=\, \delta_\by \,,
  \qquad \tensCS \tensS \,=\, \nabla^\mathfrak{s}\displacement \,,\hspace*{1.3em} \text{in $\Omega$,} \\
  &   \tensS\,\n \,=\, -\ii \omega \mathcal{Z} \displacement \,, \hspace*{6em}\text{on $\partial\Omega$} .
\end{empheq}
\end{subequations}
\flo{Its solution is a superposition
of all types of waves supported by the medium}.
In \cref{eq:numerics:eq-elastic-2d}, $\by$ is the position of the source (the center of the disk
in our experiments), and the impedance coefficient $\mathcal{Z}$ corresponds 
to a low-order absorbing boundary condition (e.g., \cref{eq:numerics:abc}).

Following our previous experiments to select scaling, we compare the 
following stabilizations:
\begin{equation} \label{eq:numerics:stabilization-u:2Dhe}
\begin{aligned}
\begin{matrix*}[l]
    \text{Identity-based:} & \tauSuIdb   \,:=\, \tauSuIdImm(\rho\cqsb) 
                                       & \,=\,  -\, \ii \omega \,\rho \, \cqsb \, \identity \,;  \\[.40em]
                           &\tauSuIda    \,:=\, \tauSuIdImm(1) 
                                       & \,=\,  -\, \ii \omega \,\, \identity \,;  \\[.50em]
      \text{KC based:} 
    & \tauSuKCb \,:=\, \tauSuKCa(1/\cqp ) & \,=\, -\dfrac{\ii\omega}{\cqp}\,\, \Gambold \,; \\[.40em]
    & \tauSuKCunit \,:=\, \tauSuKCa(1)    & \,=\, - \ii\omega \,\, \Gambold \,; \\[.50em]
      \text{Godunov based:} \hspace*{2em}
    & \tauSuMGunit \,:=\, \tauSuMGa(1)    & \,=\, -\ii\omega\,\, M_{\mathrm{Godunov}} \,.
   \end{matrix*} 
\end{aligned}
\end{equation}

\begin{remark}
  We also \flo{investigated} the qP-wavespeed in \flo{the} identity-based
  and \flo{the} sH-wavespeed for \flo{the} Kelvin--Christoffel stabilizations, 
  \flo{denoted} respectively as $\tauSuIdImm$ and $\tauSuKCa$ in \cref{eq:numerics:stabilization-u:2Dhe}.
  They give (slightly) worse results, and are not shown 
  for the sake of conciseness.
  For the identity-based stabilization, the choice of the S-wavespeed 
  \flo{for $\tau$} is further discussed in \cref{subsection:numerics-summary}.
\end{remark}

% ------------------------------------------------------------------
\subsubsection{Numerical representation of the models of parameters} 
% ------------------------------------------------------------------

One advantage \flo{in} writing the elastic system in terms 
of the compliance matrix $\tensCS$ rather than the 
stiffness tensor $\tensC$ is that it facilitates 
\flo{a} non-constant description of background materials within a mesh cell.
For this experiment, we construct a background model inspired from 
model \texttt{S}, a standard solar model in helioseismology, 
\cite{Dalsgaard1996,Pham2020Siam}. 
Model \texttt{S} provides us with radially symmetric density and wave speed, 
which is employed as \flo{a} P-wave speed for our model.
Here, we consider a two-dimensional disk geometry of 
radius 1, on which the radial profiles are \flo{defined}.

%\begin{figure}[ht!]\centering
%\subfloat[][Wavespeed model.]{\makebox[.45\linewidth][c]
%             {\includegraphics[scale=1]{figures/numerics/2D-solar/models/standalone_plot_vp}}}
%\subfloat[][Density model on logarithmic scale.]{\makebox[.45\linewidth][c]
%             {\includegraphics[scale=1]{figures/numerics/2D-solar/models/standalone_plot_rho}}}
%\caption{Solar background wave speed and density extracted from model \texttt{S}
%         for scaled radius between 0 and 1. 
%         In our elastic experiment, we use $\cs = \num{0.70}\cp$.
%         }
%\label{figure:2D-solar:models-1D}
%\end{figure}

% The models of physical properties are 
\flo{This model is}
particularly challenging because of the high variation of
amplitude, 
with an exponential decrease of the density 
close to the boundary, and wave speed variation of 
about two orders of magnitude, \cite[Figure~3]{Pham2020Siam}. 
\flo{In order to capture correctly the variation of the models,
using a piecewise-constant model representation would require}
%Therefore, it is fundamental to correctly represent
%the models within the discretized domain, and using
%a piecewise-constant approximation would require 
extreme mesh refinement near the boundary. 
\flo{Exploiting the flexibility of the formulation with compliance tensor \cref{1stsystem_variants:u},}
% Using the flexibility of the discretization with the 
% compliance tensor, 
we instead represent the physical parameters (wave speeds and density) 
using a basis of Lagrange polynomials on each cell. 
In \cref{figure:2D-solar:models} we show the mesh 
of the 2D disk \flo{employed in our simulations.
In the same figure, we compare the two options of representing 
the density model near the surface: 
piecewise constant and piecewise polynomial (here with order 2 
Lagrange polynomials) representations.
The Lagrange polynomial representation preserves the 
spherical nature of the model, which is not the case
for the piecewise constant representation}.
% The piecewise-constant representation loses the spherical 
% nature of the medium. 
% The Lagrange polynomial representation preserves the 
% spherical nature of the model.
Note that we still have a finer discretization near \flo{the} surface,
to ensure accuracy as the wavelength drastically decreases
(as the wave speed drops), and to have an accurate circular geometry.

% ----------------------------
\setlength{\modelwidth} {5.0cm}
\setlength{\modelheight}{3.8cm}
% ----------------------------
\begin{figure}[ht!] \centering
  \pgfkeys{/pgf/fpu=true} \pgfmathsetmacro{\cmin}{2e-4} \pgfmathsetmacro{\cmax}{1} \pgfkeys{/pgf/fpu=false}

  \subfloat[][Mesh of the unit disk with about \num{50000} cells.]
             {\makebox[.30\linewidth][c]{
              \includegraphics[width=4cm]{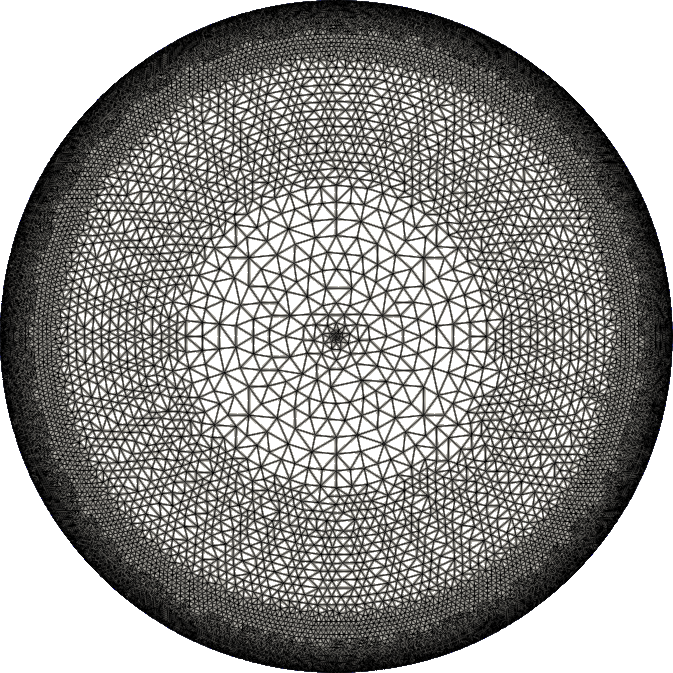}}
              \label{figure:2D-solar:models_a}} \hfill
  \renewcommand{\modelfile}{figures/numerics/2D-solar/models_v2/rho_logscale_pcste_scale2e-4-to-1_crop}
  \subfloat[][Zoom near surface for piecewise-constant representation.]
             {\makebox[.32\linewidth][c]{\begin{tikzpicture}

%\pgfmathsetmacro{\xmin}{0.00} 
%\pgfmathsetmacro{\xmax}{0.05} 
%\pgfmathsetmacro{\ymin}{0.98} 
%\pgfmathsetmacro{\ymax}{1.00} 
%\pgfmathsetmacro{\xminloc}{0.00} 
%\pgfmathsetmacro{\xmaxloc}{0.02}
%\pgfmathsetmacro{\yminloc}{0.98} 
%\pgfmathsetmacro{\ymaxloc}{1.00} 

\pgfmathsetmacro{\xmin}{-1} 
\pgfmathsetmacro{\xmax}{ 1} 
\pgfmathsetmacro{\ymin}{ 0.99} 
\pgfmathsetmacro{\ymax}{ 1.00} 
\pgfmathsetmacro{\xminloc}{\xmin} 
\pgfmathsetmacro{\xmaxloc}{\xmax}
\pgfmathsetmacro{\yminloc}{\ymin} 
\pgfmathsetmacro{\ymaxloc}{\ymax}

\begin{axis}[%
  width=\modelwidth, height=\modelheight,
  axis on top, separate axis lines,
  xmin=\xminloc, xmax=\xmaxloc, 
  ymin=\yminloc, ymax=\ymaxloc, 
  xlabel={$x$},
  ylabel={$z$}, 
% y dir=reverse,
  ytick={0.99,1},
  xtick={-.5,0,.5},
  xticklabels={\num{-0.005},0,\num{0.005}},
  yticklabels={\num{0.99},\num{1}},
ylabel style = {yshift=-2em, xshift=1.50em, rotate=90},
xlabel style = {yshift =0mm, xshift=0mm},
colormap/jet,colorbar,
colorbar style={title={$\rho$ (\si{\kg\per\meter\cubed})},
width=.5em,xshift=-.5em,ymode=log},
point meta min=\cmin,point meta max=\cmax,
label style={font=\scriptsize},
tick label style={font=\scriptsize},
legend style={font=\scriptsize\selectfont},
]
\addplot [forget plot] graphics [xmin=\xmin,xmax=\xmax,
                                 ymin=\ymin,ymax=\ymax] 
                                 {{\modelfile}.png};
\end{axis}
\end{tikzpicture}%
}} \hfill
  \renewcommand{\modelfile}{figures/numerics/2D-solar/models_v2/rho_logscale_dofp2_scale2e-4-to-1_crop}
  \subfloat[][Zoom near surface for representation with Lagrange basis of order 2.]
             {\makebox[.33\linewidth][c]{\begin{tikzpicture}

%\pgfmathsetmacro{\xmin}{0.00} 
%\pgfmathsetmacro{\xmax}{0.05} 
%\pgfmathsetmacro{\ymin}{0.98} 
%\pgfmathsetmacro{\ymax}{1.00} 
%\pgfmathsetmacro{\xminloc}{0.00} 
%\pgfmathsetmacro{\xmaxloc}{0.02}
%\pgfmathsetmacro{\yminloc}{0.98} 
%\pgfmathsetmacro{\ymaxloc}{1.00} 

\pgfmathsetmacro{\xmin}{-1} 
\pgfmathsetmacro{\xmax}{ 1} 
\pgfmathsetmacro{\ymin}{ 0.99} 
\pgfmathsetmacro{\ymax}{ 1.00} 
\pgfmathsetmacro{\xminloc}{\xmin} 
\pgfmathsetmacro{\xmaxloc}{\xmax}
\pgfmathsetmacro{\yminloc}{\ymin} 
\pgfmathsetmacro{\ymaxloc}{\ymax}

\begin{axis}[%
  width=\modelwidth, height=\modelheight,
  axis on top, separate axis lines,
  xmin=\xminloc, xmax=\xmaxloc, 
  ymin=\yminloc, ymax=\ymaxloc, 
  xlabel={$x$},
  ylabel={$z$}, 
% y dir=reverse,
  ytick={0.99,1},
  xtick={-.5,0,.5},
  xticklabels={\num{-0.005},0,\num{0.005}},
  yticklabels={\num{0.99},\num{1}},
ylabel style = {yshift=-2em, xshift=1.50em, rotate=90},
xlabel style = {yshift =0mm, xshift=0mm},
colormap/jet,colorbar,
colorbar style={title={$\rho$ (\si{\kg\per\meter\cubed})},
width=.5em,xshift=-.5em,ymode=log},
point meta min=\cmin,point meta max=\cmax,
label style={font=\scriptsize},
tick label style={font=\scriptsize},
legend style={font=\scriptsize\selectfont},
]
\addplot [forget plot] graphics [xmin=\xmin,xmax=\xmax,
                                 ymin=\ymin,ymax=\ymax] 
                                 {{\modelfile}.png};
\end{axis}
\end{tikzpicture}%
}} 

  \caption{Computational mesh and solar-like density 
           represented on the mesh near surface.
           The piecewise-constant representation is 
           unable to preserve the radial symmetry, as 
           shown in the last  (radial) layers. }
  \label{figure:2D-solar:models}
\end{figure}
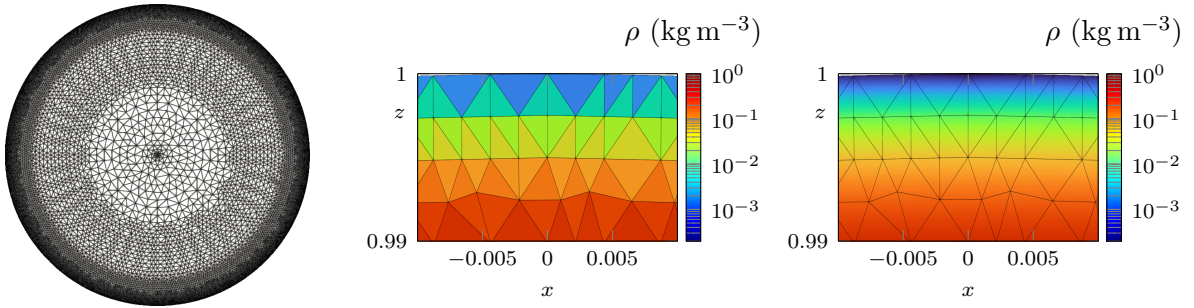

% ------------------------------------------------------------------------------- %
% \renewcommand{\experiment}{\texttt{benchmark2Diso-He}} % heterogeneous / helio
% \subsubsection{Elastic isotropy (\experiment)}
% \renewcommand{\experiment}{\texttt{Configuration HeIso}} % heterogeneous / helio
\subsubsection{Elastic isotropic medium}
\label{subsection:elastic-iso:sun}
% isotropy
\renewcommand{\plottauSuIda}{tau6}
\renewcommand{\plottauSuIdb}{tau4}
\renewcommand{\plottauSuKCa}{tau16}
\renewcommand{\plottauSuKCb}{tau18} % tau13 is using s-wavespeed
\renewcommand{\plottauSuMGa}{tau23}
\renewcommand{\plottauSvIda}{tau0}
\renewcommand{\plottauSvIdb}{tau1}
\renewcommand{\plottauSvKCb}{tau11}
\renewcommand{\plottauSvMGa}{tau20}
% ------------------------------------------------------------------------------- %

We first consider an isotropic medium, with the P-wavespeed
and density given by the solar background, and S-wavespeed
selected such that $\cs=\num{0.70}\cp$. 
The Dirac point-source is positioned at the center of 
the unit disk, and absorbing conditions are implemented
on the boundary, \cite{Givoli1990,Higdon1991}, such that
\begin{equation} \label{eq:numerics:abc}
  \tensS\cdot\n  \,-\, 
  \ii \omega \Big( \sqrt{(\lambda+2\mu)\rho} \,\,\, \n\otimes\n \,+\, \sqrt{\mu\,\rho}\,\,
  (\mathbb{I_d} \,-\, \n\otimes\n )\Big)
  \displacement \,=\, 0 \,, \qquad\quad \text{on $\partial\Omega$.}
\end{equation}
% In the case of (strongly) heterogeneous properties, 
% there is no analytic solution for the elastic wave 
% equation. 
\flo{Since there are no analytic solutions in this case, 
to investigate accuracy, a reference solution is constructed 
with a refined mesh of}
%To investigate the accuracy associated with each stabilization, 
%we build a reference solution which is computed with a refined 
%mesh composed of 
\num{80000} cells (while for the simulation we have
\num{50000} cells, \cref{figure:2D-solar:models_a}),
polynomials of order 7, and with stabilization
$\tauSuMGunit$. 
The reference solutions are shown in \cref{figure:2D-solar:waves}
% where we also compare between piecewise-constant and polynomial representation.
where we also compare with the solution 
\flo{obtained with a piecewise-constant model representation}.
Due to the exponentially decreasing nature of the density,
we scale the displacement fields by $\sqrt{\rho}\,$ to 
\flo{give a better visualization of} the solution, \cite{Pham2020Siam,Pham2019Esaim}.
For \flo{a} similar reason, the components of \flo{the} stress tensor are 
scaled by $1/\sqrt{\rho}$.

\setlength{\modelwidth} {4.10cm}
\setlength{\modelheight}{4.10cm}
% ----------------------------
\begin{figure}[ht!] \centering
  \pgfkeys{/pgf/fpu=true} \pgfmathsetmacro{\cmin}{-1} \pgfmathsetmacro{\cmax}{1} \pgfkeys{/pgf/fpu=false}
  \renewcommand{\cbtitle}{\scriptsize $\cdot$\num{e3}}
  \renewcommand{\modelfile}{figures/numerics/2D-solar/wavefield/UxsqrtRho_4mHz_real_tau2p6_scale1e3}
  \subfloat[][$\sqrt{\rho}\,\, \Real(u_x)$.]
             {\begin{tikzpicture}
\pgfmathsetmacro{\xminloc}{-1.00} 
\pgfmathsetmacro{\xmaxloc}{ 1.00}
\pgfmathsetmacro{\yminloc}{-1.00} 
\pgfmathsetmacro{\ymaxloc}{ 1.00}
\pgfmathsetmacro{\xmin}{-1} 
\pgfmathsetmacro{\xmax}{ 1} 
\pgfmathsetmacro{\ymin}{-1} 
\pgfmathsetmacro{\ymax}{ 1}

\begin{axis}[%
  width=\modelwidth, height=\modelheight,
  axis on top, separate axis lines,
  xmin=\xminloc, xmax=\xmaxloc, 
  ymin=\yminloc, ymax=\ymaxloc, 
  xlabel={$x$},
  ylabel={$z$},
  xtick={-.5,0,0.5}, 
  y dir=reverse,
ylabel style = {yshift=-1em, xshift=0mm},
xlabel style = {yshift =1.10em, xshift=2.80em},
%colormap/jet,colorbar,
%colorbar style={title={\cbtitle},
%width=.5em,xshift=-.80em,% ymode=log
%},
%point meta min=\cmin,point meta max=\cmax,
label style={font=\scriptsize},
tick label style={font=\scriptsize},
legend style={font=\scriptsize\selectfont},
]
\addplot [forget plot] graphics [xmin=\xmin,xmax=\xmax,
                                 ymin=\ymin,ymax=\ymax] 
                                 {{\modelfile}.png};
\end{axis}
\end{tikzpicture}%
} \hspace*{-1.0em}
  \renewcommand{\modelfile}{figures/numerics/2D-solar/wavefield/UzsqrtRho_4mHz_real_tau2p6_scale1e3}
  \subfloat[][$\sqrt{\rho}\,\, \Real(u_z)$.]
             {\begin{tikzpicture}
\pgfmathsetmacro{\xminloc}{-1.00} 
\pgfmathsetmacro{\xmaxloc}{ 1.00}
\pgfmathsetmacro{\yminloc}{-1.00} 
\pgfmathsetmacro{\ymaxloc}{ 1.00}
\pgfmathsetmacro{\xmin}{-1} 
\pgfmathsetmacro{\xmax}{ 1} 
\pgfmathsetmacro{\ymin}{-1} 
\pgfmathsetmacro{\ymax}{ 1}

\begin{axis}[%
  width=\modelwidth, height=\modelheight,
  axis on top, separate axis lines,
  xmin=\xminloc, xmax=\xmaxloc, 
  ymin=\yminloc, ymax=\ymaxloc, 
  xlabel={$x$},
  xtick={-.5,0,0.5},
  xlabel style = {yshift =1.10em, xshift=2.80em},
  y dir=reverse,
  ytick={},
  yticklabels={,,},
colormap/jet,colorbar,
colorbar style={title={\cbtitle},
width=.5em,xshift=-0.60em,% ymode=log
title style={yshift=-.60em},
},
point meta min=\cmin,point meta max=\cmax,
label style={font=\scriptsize},
tick label style={font=\scriptsize},
legend style={font=\scriptsize\selectfont},
]
\addplot [forget plot] graphics [xmin=\xmin,xmax=\xmax,
                                 ymin=\ymin,ymax=\ymax] 
                                 {{\modelfile}.png};
\end{axis}
\end{tikzpicture}%
} \hspace*{-2.0em}
  \renewcommand{\cbtitle}{\scriptsize $\cdot$\num{e-3}}
  \pgfkeys{/pgf/fpu=true} \pgfmathsetmacro{\cmin}{-5} \pgfmathsetmacro{\cmax}{5} \pgfkeys{/pgf/fpu=false}
  \renewcommand{\modelfile}{figures/numerics/2D-solar/wavefield/Sxx-over-sqrtRho_4mHz_real_tau2p6_scale5e-3} 
  \subfloat[][$\Real(\sigma_{xx}) \,/\, \sqrt{\rho}$.]
             {\begin{tikzpicture}
\pgfmathsetmacro{\xminloc}{-1.00} 
\pgfmathsetmacro{\xmaxloc}{ 1.00}
\pgfmathsetmacro{\yminloc}{-1.00} 
\pgfmathsetmacro{\ymaxloc}{ 1.00}
\pgfmathsetmacro{\xmin}{-1} 
\pgfmathsetmacro{\xmax}{ 1} 
\pgfmathsetmacro{\ymin}{-1} 
\pgfmathsetmacro{\ymax}{ 1}

\begin{axis}[%
  width=\modelwidth, height=\modelheight,
  axis on top, separate axis lines,
  xmin=\xminloc, xmax=\xmaxloc, 
  ymin=\yminloc, ymax=\ymaxloc, 
  xlabel={$x$},
  xtick={-.5,0,0.5},
  y dir=reverse,
  ytick={},
  yticklabels={,,},
  xlabel style = {yshift =1.10em, xshift=2.80em},
% ytick={\ymin,\ymax},
% xtick={\xmin,1,\xmax},
% xticklabels={\num{\xmin},1,\num{\xmax}},
% yticklabels={\num{\ymin},\ymax},
% ylabel style = {yshift=-1em, xshift=0mm},
%xlabel style = {yshift =0mm, xshift=0mm},
%colormap/jet,colorbar,
%colorbar style={title={\cbtitle},
%width=.5em,xshift=-0.20em,% ymode=log
%},
%point meta min=\cmin,point meta max=\cmax,
label style={font=\scriptsize},
tick label style={font=\scriptsize},
legend style={font=\scriptsize\selectfont},
]
\addplot [forget plot] graphics [xmin=\xmin,xmax=\xmax,
                                 ymin=\ymin,ymax=\ymax] 
                                 {{\modelfile}.png};
\end{axis}
\end{tikzpicture}%
} \hspace*{-1.0em}
  \renewcommand{\modelfile}{figures/numerics/2D-solar/wavefield/Szz-over-sqrtRho_4mHz_real_tau2p6_scale5e-3} 
  \subfloat[][$\Real(\sigma_{zz}) \,/\, \sqrt{\rho}$.]
             {\begin{tikzpicture}
\pgfmathsetmacro{\xminloc}{-1.00} 
\pgfmathsetmacro{\xmaxloc}{ 1.00}
\pgfmathsetmacro{\yminloc}{-1.00} 
\pgfmathsetmacro{\ymaxloc}{ 1.00}
\pgfmathsetmacro{\xmin}{-1} 
\pgfmathsetmacro{\xmax}{ 1} 
\pgfmathsetmacro{\ymin}{-1} 
\pgfmathsetmacro{\ymax}{ 1}

\begin{axis}[%
  width=\modelwidth, height=\modelheight,
  axis on top, separate axis lines,
  xmin=\xminloc, xmax=\xmaxloc, 
  ymin=\yminloc, ymax=\ymaxloc, 
  xlabel={$x$},
  xtick={-.5,0,0.5},
  y dir=reverse,
  ytick={},
  yticklabels={,,},
  xlabel style = {yshift =1.10em, xshift=2.80em},
% ytick={\ymin,\ymax},
% xtick={\xmin,1,\xmax},
% xticklabels={\num{\xmin},1,\num{\xmax}},
% yticklabels={\num{\ymin},\ymax},
% ylabel style = {yshift=-1em, xshift=0mm},
%xlabel style = {yshift =0mm, xshift=0mm},
%colormap/jet,colorbar,
%colorbar style={title={\cbtitle},
%width=.5em,xshift=-0.20em,% ymode=log
%},
%point meta min=\cmin,point meta max=\cmax,
label style={font=\scriptsize},
tick label style={font=\scriptsize},
legend style={font=\scriptsize\selectfont},
]
\addplot [forget plot] graphics [xmin=\xmin,xmax=\xmax,
                                 ymin=\ymin,ymax=\ymax] 
                                 {{\modelfile}.png};
\end{axis}
\end{tikzpicture}%
} \hspace*{-1.0em}
  \renewcommand{\modelfile}{figures/numerics/2D-solar/wavefield/Sxz-over-sqrtRho_4mHz_real_tau2p6_scale5e-3} 
  \subfloat[][$\Real(\sigma_{xz}) \,/\, \sqrt{\rho}$.]
             {\begin{tikzpicture}
\pgfmathsetmacro{\xminloc}{-1.00} 
\pgfmathsetmacro{\xmaxloc}{ 1.00}
\pgfmathsetmacro{\yminloc}{-1.00} 
\pgfmathsetmacro{\ymaxloc}{ 1.00}
\pgfmathsetmacro{\xmin}{-1} 
\pgfmathsetmacro{\xmax}{ 1} 
\pgfmathsetmacro{\ymin}{-1} 
\pgfmathsetmacro{\ymax}{ 1}

\begin{axis}[%
  width=\modelwidth, height=\modelheight,
  axis on top, separate axis lines,
  xmin=\xminloc, xmax=\xmaxloc, 
  ymin=\yminloc, ymax=\ymaxloc, 
  xlabel={$x$},
  xtick={-.5,0,0.5},
  xlabel style = {yshift =1.10em, xshift=2.80em},
  y dir=reverse,
  ytick={},
  yticklabels={,,},
colormap/jet,colorbar,
colorbar style={title={\cbtitle},
width=.5em,xshift=-0.60em,% ymode=log
title style={yshift=-.60em},
},
point meta min=\cmin,point meta max=\cmax,
label style={font=\scriptsize},
tick label style={font=\scriptsize},
legend style={font=\scriptsize\selectfont},
]
\addplot [forget plot] graphics [xmin=\xmin,xmax=\xmax,
                                 ymin=\ymin,ymax=\ymax] 
                                 {{\modelfile}.png};
\end{axis}
\end{tikzpicture}%
}  \\[-1em]

  \pgfkeys{/pgf/fpu=true} \pgfmathsetmacro{\cmin}{-1} \pgfmathsetmacro{\cmax}{1} \pgfkeys{/pgf/fpu=false}
  \renewcommand{\cbtitle}{\scriptsize $\cdot$\num{e3}}
  \renewcommand{\modelfile}{figures/numerics/2D-solar/wavefield/UxsqrtRho_4mHz_real_tau2p6_scale1e3_pconstant}
  \subfloat[][$\sqrt{\rho}\,\, \Real(u_x)$.]
             {\begin{tikzpicture}
\pgfmathsetmacro{\xminloc}{-1.00} 
\pgfmathsetmacro{\xmaxloc}{ 1.00}
\pgfmathsetmacro{\yminloc}{-1.00} 
\pgfmathsetmacro{\ymaxloc}{ 1.00}
\pgfmathsetmacro{\xmin}{-1} 
\pgfmathsetmacro{\xmax}{ 1} 
\pgfmathsetmacro{\ymin}{-1} 
\pgfmathsetmacro{\ymax}{ 1}

\begin{axis}[%
  width=\modelwidth, height=\modelheight,
  axis on top, separate axis lines,
  xmin=\xminloc, xmax=\xmaxloc, 
  ymin=\yminloc, ymax=\ymaxloc, 
  xlabel={$x$},
  ylabel={$z$},
  xtick={-.5,0,0.5}, 
  y dir=reverse,
ylabel style = {yshift=-1em, xshift=0mm},
xlabel style = {yshift =1.10em, xshift=2.80em},
%colormap/jet,colorbar,
%colorbar style={title={\cbtitle},
%width=.5em,xshift=-.80em,% ymode=log
%},
%point meta min=\cmin,point meta max=\cmax,
label style={font=\scriptsize},
tick label style={font=\scriptsize},
legend style={font=\scriptsize\selectfont},
]
\addplot [forget plot] graphics [xmin=\xmin,xmax=\xmax,
                                 ymin=\ymin,ymax=\ymax] 
                                 {{\modelfile}.png};
\end{axis}
\end{tikzpicture}%
} \hspace*{-1.0em}
  \renewcommand{\modelfile}{figures/numerics/2D-solar/wavefield/UzsqrtRho_4mHz_real_tau2p6_scale1e3_pconstant}
  \subfloat[][$\sqrt{\rho}\,\, \Real(u_z)$.]
             {\begin{tikzpicture}
\pgfmathsetmacro{\xminloc}{-1.00} 
\pgfmathsetmacro{\xmaxloc}{ 1.00}
\pgfmathsetmacro{\yminloc}{-1.00} 
\pgfmathsetmacro{\ymaxloc}{ 1.00}
\pgfmathsetmacro{\xmin}{-1} 
\pgfmathsetmacro{\xmax}{ 1} 
\pgfmathsetmacro{\ymin}{-1} 
\pgfmathsetmacro{\ymax}{ 1}

\begin{axis}[%
  width=\modelwidth, height=\modelheight,
  axis on top, separate axis lines,
  xmin=\xminloc, xmax=\xmaxloc, 
  ymin=\yminloc, ymax=\ymaxloc, 
  xlabel={$x$},
  xtick={-.5,0,0.5},
  xlabel style = {yshift =1.10em, xshift=2.80em},
  y dir=reverse,
  ytick={},
  yticklabels={,,},
colormap/jet,colorbar,
colorbar style={title={\cbtitle},
width=.5em,xshift=-0.60em,% ymode=log
title style={yshift=-.60em},
},
point meta min=\cmin,point meta max=\cmax,
label style={font=\scriptsize},
tick label style={font=\scriptsize},
legend style={font=\scriptsize\selectfont},
]
\addplot [forget plot] graphics [xmin=\xmin,xmax=\xmax,
                                 ymin=\ymin,ymax=\ymax] 
                                 {{\modelfile}.png};
\end{axis}
\end{tikzpicture}%
} \hspace*{-2.0em}
  \renewcommand{\cbtitle}{\scriptsize $\cdot$\num{e-3}}
  \pgfkeys{/pgf/fpu=true} \pgfmathsetmacro{\cmin}{-5} \pgfmathsetmacro{\cmax}{5} \pgfkeys{/pgf/fpu=false}
  \renewcommand{\modelfile}{figures/numerics/2D-solar/wavefield/Sxx-over-sqrtRho_4mHz_real_tau2p6_scale5e-3_pconstant} 
  \subfloat[][$\Real(\sigma_{xx}) \,/\, \sqrt{\rho}$.]
             {\begin{tikzpicture}
\pgfmathsetmacro{\xminloc}{-1.00} 
\pgfmathsetmacro{\xmaxloc}{ 1.00}
\pgfmathsetmacro{\yminloc}{-1.00} 
\pgfmathsetmacro{\ymaxloc}{ 1.00}
\pgfmathsetmacro{\xmin}{-1} 
\pgfmathsetmacro{\xmax}{ 1} 
\pgfmathsetmacro{\ymin}{-1} 
\pgfmathsetmacro{\ymax}{ 1}

\begin{axis}[%
  width=\modelwidth, height=\modelheight,
  axis on top, separate axis lines,
  xmin=\xminloc, xmax=\xmaxloc, 
  ymin=\yminloc, ymax=\ymaxloc, 
  xlabel={$x$},
  xtick={-.5,0,0.5},
  y dir=reverse,
  ytick={},
  yticklabels={,,},
  xlabel style = {yshift =1.10em, xshift=2.80em},
% ytick={\ymin,\ymax},
% xtick={\xmin,1,\xmax},
% xticklabels={\num{\xmin},1,\num{\xmax}},
% yticklabels={\num{\ymin},\ymax},
% ylabel style = {yshift=-1em, xshift=0mm},
%xlabel style = {yshift =0mm, xshift=0mm},
%colormap/jet,colorbar,
%colorbar style={title={\cbtitle},
%width=.5em,xshift=-0.20em,% ymode=log
%},
%point meta min=\cmin,point meta max=\cmax,
label style={font=\scriptsize},
tick label style={font=\scriptsize},
legend style={font=\scriptsize\selectfont},
]
\addplot [forget plot] graphics [xmin=\xmin,xmax=\xmax,
                                 ymin=\ymin,ymax=\ymax] 
                                 {{\modelfile}.png};
\end{axis}
\end{tikzpicture}%
} \hspace*{-1.0em}
  \renewcommand{\modelfile}{figures/numerics/2D-solar/wavefield/Szz-over-sqrtRho_4mHz_real_tau2p6_scale5e-3_pconstant} 
  \subfloat[][$\Real(\sigma_{zz}) \,/\, \sqrt{\rho}$.]
             {\begin{tikzpicture}
\pgfmathsetmacro{\xminloc}{-1.00} 
\pgfmathsetmacro{\xmaxloc}{ 1.00}
\pgfmathsetmacro{\yminloc}{-1.00} 
\pgfmathsetmacro{\ymaxloc}{ 1.00}
\pgfmathsetmacro{\xmin}{-1} 
\pgfmathsetmacro{\xmax}{ 1} 
\pgfmathsetmacro{\ymin}{-1} 
\pgfmathsetmacro{\ymax}{ 1}

\begin{axis}[%
  width=\modelwidth, height=\modelheight,
  axis on top, separate axis lines,
  xmin=\xminloc, xmax=\xmaxloc, 
  ymin=\yminloc, ymax=\ymaxloc, 
  xlabel={$x$},
  xtick={-.5,0,0.5},
  y dir=reverse,
  ytick={},
  yticklabels={,,},
  xlabel style = {yshift =1.10em, xshift=2.80em},
% ytick={\ymin,\ymax},
% xtick={\xmin,1,\xmax},
% xticklabels={\num{\xmin},1,\num{\xmax}},
% yticklabels={\num{\ymin},\ymax},
% ylabel style = {yshift=-1em, xshift=0mm},
%xlabel style = {yshift =0mm, xshift=0mm},
%colormap/jet,colorbar,
%colorbar style={title={\cbtitle},
%width=.5em,xshift=-0.20em,% ymode=log
%},
%point meta min=\cmin,point meta max=\cmax,
label style={font=\scriptsize},
tick label style={font=\scriptsize},
legend style={font=\scriptsize\selectfont},
]
\addplot [forget plot] graphics [xmin=\xmin,xmax=\xmax,
                                 ymin=\ymin,ymax=\ymax] 
                                 {{\modelfile}.png};
\end{axis}
\end{tikzpicture}%
} \hspace*{-1.0em}
  \renewcommand{\modelfile}{figures/numerics/2D-solar/wavefield/Sxz-over-sqrtRho_4mHz_real_tau2p6_scale5e-3_pconstant} 
  \subfloat[][$\Real(\sigma_{xz}) \,/\, \sqrt{\rho}$.]
             {\begin{tikzpicture}
\pgfmathsetmacro{\xminloc}{-1.00} 
\pgfmathsetmacro{\xmaxloc}{ 1.00}
\pgfmathsetmacro{\yminloc}{-1.00} 
\pgfmathsetmacro{\ymaxloc}{ 1.00}
\pgfmathsetmacro{\xmin}{-1} 
\pgfmathsetmacro{\xmax}{ 1} 
\pgfmathsetmacro{\ymin}{-1} 
\pgfmathsetmacro{\ymax}{ 1}

\begin{axis}[%
  width=\modelwidth, height=\modelheight,
  axis on top, separate axis lines,
  xmin=\xminloc, xmax=\xmaxloc, 
  ymin=\yminloc, ymax=\ymaxloc, 
  xlabel={$x$},
  xtick={-.5,0,0.5},
  xlabel style = {yshift =1.10em, xshift=2.80em},
  y dir=reverse,
  ytick={},
  yticklabels={,,},
colormap/jet,colorbar,
colorbar style={title={\cbtitle},
width=.5em,xshift=-0.60em,% ymode=log
title style={yshift=-.60em},
},
point meta min=\cmin,point meta max=\cmax,
label style={font=\scriptsize},
tick label style={font=\scriptsize},
legend style={font=\scriptsize\selectfont},
]
\addplot [forget plot] graphics [xmin=\xmin,xmax=\xmax,
                                 ymin=\ymin,ymax=\ymax] 
                                 {{\modelfile}.png};
\end{axis}
\end{tikzpicture}%
}

  \caption{Comparison of the reference solutions at \num{4} \si{\milli\Hz} 
           of the solar-like isotropic elastic background  models, % of \cref{figure:2D-solar:models-1D}
           obtained using piecewise-constant representation (bottom) or 
           with piecewise polynomial representation (Lagrange basis \modif{of} order 2) 
           (top).
           The computations use HDG formulation $\formulationU$
           with stabilization $\tauSuMGunit$ and polynomial order 7.}
  \label{figure:2D-solar:waves}
       
\end{figure}
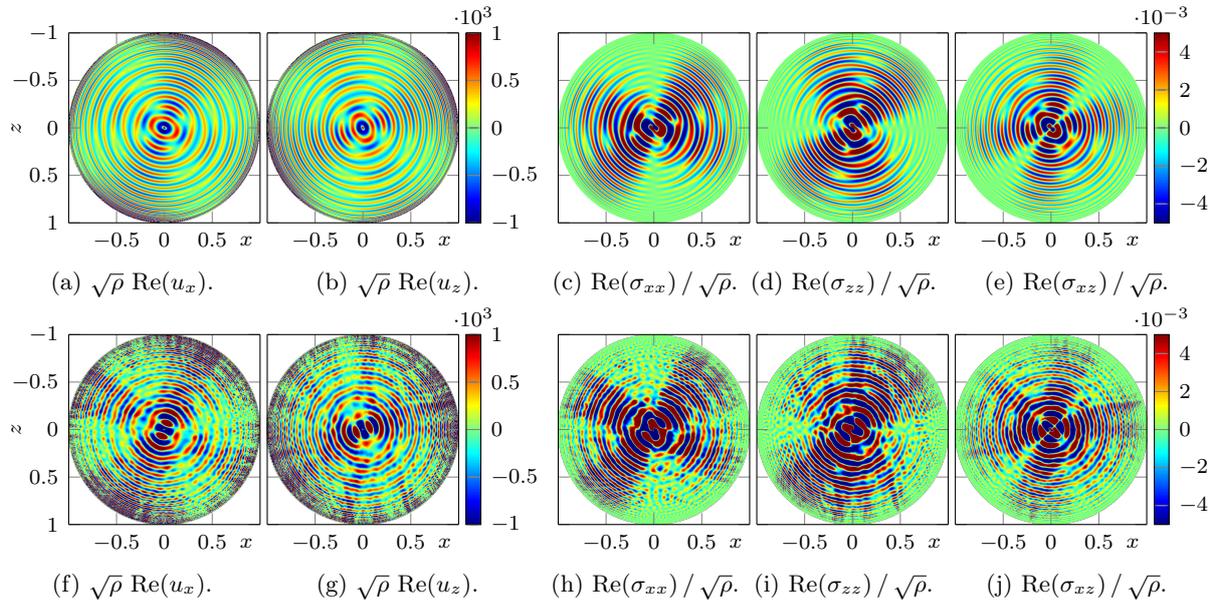

% Comparing simulations using a piecewise-constant representation
% for the model parameters or allowing them to vary within the cell,
%we see drastic difference. 
\flo{The solutions corresponding to a piecewise-constant model representation}
% Using a piecewise-constant representation,
% the solutions 
(both displacement $\displacement$ and stress tensor $\tensS$) 
show \flo{ripples} and artifacts, in addition, 
their radial nature is not preserved.
On the contrary, \flo{results employing the Lagrange basis model representation}
capture well the spherical pattern
and provide `clean' solutions.
With such high variations in the background models, we see that
it is mandatory to design an efficient representation, and that
piecewise-constant model representation is not appropriate. 
In this way, the HDG formulation based on the compliance tensor $\tensCS$ is 
useful as it allows us to vary easily the models within each cell.
% without having to compute their derivatives (as it would be the
% case using a formulation with the stiffness tensor $\tensC$).

In \cref{figure:2D-solar:Su-error-freq-order_01}, we plot the 
relative error $\errall$ \cref{eq:relative-error_sum} 
in terms of the frequency and polynomial orders, for the 
stabilization coefficients of \cref{eq:numerics:stabilization-u:2Dhe}.
Note that for some frequencies and orders, stabilization 
$\tauSuIda$ results in \flo{an} ill-conditioned matrix that 
cannot be factorized or leads to erroneous solutions.
In these cases, the discretized wave system with $\tauSuIda$ cannot 
be solved. \flo{However, when it works, using $\tauSuIda$ 
still gives worst performance with high levels of error.}
\begin{itemize}\setlength{\itemsep}{-2pt}
  \item We observe that the error increases with increasing frequency, 
        \cref{figure:2D-solar:Su-error-freq-order_01a}, which is common 
        as the wavelength is reduced. 
  \item As expected, the error also decreases with increasing polynomial
        orders, cf. \cref{figure:2D-solar:Su-error-freq-order_01b}.
%  \item Comparing the stabilizations, $\tauSuKCunit$ and $\tauSuIda$ 
%        (when working)  are the worst, with higher levels of error in 
%        all cases.
  \item The stabilization $\tauSuMGunit$ gives accurate results, 
        with an accuracy that can be match\flo{ed} by using the appropriate scaling
        in other stabilizations, 
        with either $\tauSuIdb$ or $\tauSuKCb$.
\end{itemize}
Overall, this test confirms %heterogeneous test confirms the observation 
that $\tauSuMGunit$ is the most versatile choice. 
While \flo{an} equivalent accuracy can be obtained with other stabilizations, 
it requires a judicial choice of scaling.

% ---------------------------------------------------------------------- %
\setlength{\plotwidth} {6.50cm}
\setlength{\plotheight}{3.10cm}
\renewcommand{\myfreqplot} {0.006}
\renewcommand{\myorderplot}{5}
% ----------------------------
\begin{figure}[ht!] \centering

  \renewcommand{\myxlabel}{frequency (\si{\milli\Hz})}
  \pgfkeys{/pgf/fpu=true}
  \pgfmathsetmacro{\xmin}{1.30} \pgfmathsetmacro{\xmax}{10.20} 
  \pgfmathsetmacro{\ymin}{2e-6} \pgfmathsetmacro{\ymax}{7e-3} 
  \pgfmathsetmacro{\scalefreq}{1e3} % cs/freq
  \pgfkeys{/pgf/fpu=false}
  \renewcommand{\datafile}{figures/numerics/2D-solar/data_cs0.70cp_error0.05to0.99_Su_reference80kTau23p7_mesh50k/error-scaled_mean_per-freq_u-and-sigmas.txt}

  \renewcommand{\dataA}{sumU_\plottauSuIdb_p\myorderplot} \renewcommand{\legendA}{$\tauSuIdb$} 
  \renewcommand{\dataB}{sumU_\plottauSuKCb_p\myorderplot} \renewcommand{\legendB}{$\tauSuKCb$}
  \renewcommand{\dataC}{sumU_\plottauSuMGa_p\myorderplot} \renewcommand{\legendC}{$\tauSuMGunit$}
  \renewcommand{\dataD}{sumU_\plottauSuIda_p\myorderplot} \renewcommand{\legendD}{$\tauSuIda$} 
  \renewcommand{\dataE}{sumU_\plottauSuKCa_p\myorderplot} \renewcommand{\legendE}{$\tauSuKCunit$}
  \renewcommand{\legendpos}{south east}
  \renewcommand{\myylabel}{$\errall(\displacement)$}
  \subfloat[][$\errall(\displacement)$ for polynomial order \myorderplot.]
             {\makebox[.35\linewidth][c]{%%%%%%%%%%%%%%%%%%%%%%%%%%%%%%%%%%%%%%%%%%%%%%%%%%%%%%%%%%%%%%%%%%%%%%%%%%%%%%%%
%% begin the TikZ picture
\begin{tikzpicture}
%%%%%%%%%%%%%%%%%%%%%%%%%%%%%%%%%%%%%%%%%%%%%%%%%%%%%%%%%%%%%%%%%%%%%%%%%%%%%%%%
\begin{axis}[
             enlargelimits=false, 
             xlabel={\scriptsize \myxlabel },
             ylabel={\scriptsize \myylabel },
             enlarge y limits=false,
             enlarge x limits=false,
             xmin=\xmin,xmax=\xmax, 
             ymin=\ymin,ymax=\ymax,
             yminorticks=true,
             height=\plotheight,width=\plotwidth,
             scale only axis,
             ylabel style = {yshift=-0.2em, xshift=0mm},
             xlabel style = {yshift =0.5em, xshift=0mm},
             ytick={1e-6,1e-5,1e-4,1e-3,1e-2,1e-1},
             xtick={0,1,2,3,4,5,6,7,8,9,10,11,12,13,14,15,16,17,18,19,20},
             xticklabels={,,2,,4,,6,,8,,10,,12,,14,,16,,18,,20},
             clip mode=individual,
             minor grid style={line width=.5pt, draw=gray!50, densely dotted},
             major grid style={line width=.5pt, draw=gray!50, densely dashed},
             xmajorgrids=true,
             xminorgrids=true,
             yminorgrids=true,
             ymajorgrids=true,
             ymode=log,
            %xmode=log,
             %y  tick label style={/pgf/number format/fixed},
             %y  tick label style={/pgf/number format/fixed zerofill},
             %y  tick label style={/pgf/number format/precision=1},
             %y  tick label style={/pgf/number format/precision=2},
             % y  tick label style={/pgf/number format/sci},
             %y tick label style={/pgf/number format/.cd},
             label style={font=\scriptsize},
             tick label style={font=\scriptsize},
             legend style={font=\scriptsize\selectfont},
             legend pos={\legendpos}, 
             legend columns=2,
             ]  %% foreground marks

     \pgfmathsetmacro{\scale}{1}

     %% load current data
     \ifthenelse{\equal{\dataA}{}}
     {}
     {
     \addplot[color=black,densely dotted,line width=0.2, mark=x, 
              mark size=4.50,mark options={solid,fill=black, line width=1.0}]
              table[x expr = \scalefreq*\thisrow{freq}, 
                    y expr = \scale*\thisrow{\dataA}, % select coords between index={1}{5}
                   ]
              {\datafile}; \addlegendentry{\legendA}    
     }
     \ifthenelse{\equal{\dataD}{}}
     {}
     {
     \addplot[color=blue,densely dashed, line width=0.2, mark=+,
              mark size=3.50,mark options={solid,fill=blue, line width=0.70}]
              table[x expr = \scalefreq*\thisrow{freq}, 
                    y expr = \scale*\thisrow{\dataD}]
              {\datafile}; \addlegendentry{\legendD}
     }

     \ifthenelse{\equal{\dataC}{}}
     {}
     {
     \addplot[color=\mygreen,densely dashed, line width=0.2, mark=star,
              mark size=4.50,mark options={solid,fill=\mygreen, line width=.75}]
              table[x expr = \scalefreq*\thisrow{freq}, 
                    y expr = \scale*\thisrow{\dataC}]
              {\datafile}; \addlegendentry{\legendC}
     }

     \ifthenelse{\equal{\dataE}{}}
     {}
     {
     \addplot[color=magenta,densely dotted, line width=0.2, mark=triangle,
              mark size=3.50,mark options={solid,fill=magenta, line width=0.50}]
              table[x expr = \scalefreq*\thisrow{freq}, 
                    y expr = \scale*\thisrow{\dataE}]
              {\datafile}; \addlegendentry{\legendE}
     }
     
     \ifthenelse{\equal{\dataB}{}}
     {}
     {
     \addplot[color=red,densely dashed, line width=0.2, mark=o,
              mark size=2.50,mark options={solid,fill=\myred, line width=.60}]
              table[x expr = \scalefreq*\thisrow{freq}, 
                    y expr = \scale*\thisrow{\dataB}]
              {\datafile}; \addlegendentry{\legendB}
     }

%%%%%%%%%%%%%%%%%%%%%%%%%%%%%%%%%%%%%%%%%%%%%%%%%%%%%%%%%%%%%%%%%%%%%%%%%%%%%%%%
\end{axis}
%%%%%%%%%%%%%%%%%%%%%%%%%%%%%%%%%%%%%%%%%%%%%%%%%%%%%%%%%%%%%%%%%%%%%%%%%%%%%%%%
\end{tikzpicture}
             \label{figure:2D-solar:Su-error-freq-order_01a}}}
              \hspace*{.50em}
  \renewcommand{\myxlabel}{polynomial order}
  \renewcommand{\legendpos}{south west}
  \pgfkeys{/pgf/fpu=true}
  \pgfmathsetmacro{\xmin}{0.60} \pgfmathsetmacro{\xmax}{7.20} 
  \pgfmathsetmacro{\ymin}{1e-5} \pgfmathsetmacro{\ymax}{5e-2} 
  \pgfkeys{/pgf/fpu=false}
  \renewcommand{\datafile}{figures/numerics/2D-solar/data_cs0.70cp_error0.05to0.99_Su_reference80kTau23p7_mesh50k/error-scaled_mean_per-order_u-and-sigmas.txt}
  \renewcommand{\dataA}{sumU_\plottauSuIdb_freq\myfreqplot} \renewcommand{\legendA}{$\tauSuIdb$} 
  \renewcommand{\dataB}{sumU_\plottauSuKCb_freq\myfreqplot} \renewcommand{\legendB}{$\tauSuKCb$}
  \renewcommand{\dataC}{sumU_\plottauSuMGa_freq\myfreqplot} \renewcommand{\legendC}{$\tauSuMGunit$}
  \renewcommand{\dataD}{sumU_\plottauSuIda_freq\myfreqplot} \renewcommand{\legendD}{$\tauSuIda$}
  \renewcommand{\dataE}{sumU_\plottauSuKCa_freq\myfreqplot} \renewcommand{\legendE}{$\tauSuKCunit$}   
  \renewcommand{\myylabel}{$\errall(\displacement)$}
  \subfloat[][$\errall(\displacement)$ at frequency 6~\si{\milli\Hz}.]
             {\makebox[.35\linewidth][c]{%%%%%%%%%%%%%%%%%%%%%%%%%%%%%%%%%%%%%%%%%%%%%%%%%%%%%%%%%%%%%%%%%%%%%%%%%%%%%%%%
%% begin the TikZ picture
\begin{tikzpicture}
%%%%%%%%%%%%%%%%%%%%%%%%%%%%%%%%%%%%%%%%%%%%%%%%%%%%%%%%%%%%%%%%%%%%%%%%%%%%%%%%
\begin{axis}[
             enlargelimits=false, 
             xlabel={\scriptsize \myxlabel },
             ylabel={\scriptsize \myylabel },
             enlarge y limits=false,
             enlarge x limits=false,
             xmin=\xmin,xmax=\xmax, 
             ymin=\ymin,ymax=\ymax,
             yminorticks=true,
             height=\plotheight,width=\plotwidth,
             scale only axis,
             ylabel style = {yshift=-0.2em, xshift=0mm},
             xlabel style = {yshift =0.5em, xshift=0mm},
             ytick={1e-6,1e-5,1e-4,1e-3,1e-2,1e-1},
             xtick={0,1,2,3,4,5,6,7,8,9,10,11,12,13,14,15,16,17,18,19,20},
             xticklabels={,,2,,4,,6,,8,,10,,12,,14,,16,,18,,20},
             clip mode=individual,
             minor grid style={line width=.5pt, draw=gray!50, densely dotted},
             major grid style={line width=.5pt, draw=gray!50, densely dashed},
             xmajorgrids=true,
             xminorgrids=true,
             yminorgrids=true,
             ymajorgrids=true,
             ymode=log,
             %y  tick label style={/pgf/number format/fixed},
             %y  tick label style={/pgf/number format/fixed zerofill},
             %y  tick label style={/pgf/number format/precision=1},
             %y  tick label style={/pgf/number format/precision=2},
             % y  tick label style={/pgf/number format/sci},
             %y tick label style={/pgf/number format/.cd},
             label style={font=\scriptsize},
             tick label style={font=\scriptsize},
             legend style={font=\scriptsize\selectfont},
             legend pos={\legendpos}, 
             legend columns=2,
             ]  %% foreground marks

     \pgfmathsetmacro{\scale}{1}

     %% load current data
     \ifthenelse{\equal{\dataA}{}}
     {}
     {
     \addplot[color=black,densely dotted,line width=0.2, mark=x, 
              mark size=4.50,mark options={solid,fill=black, line width=1.0}]
              table[x expr = \thisrow{order}, 
                    y expr = \scale*\thisrow{\dataA}, % select coords between index={1}{5}
                   ]
              {\datafile}; \addlegendentry{\legendA}    
     }

     \ifthenelse{\equal{\dataD}{}}
     {}
     {
     \addplot[color=blue,densely dashed, line width=0.2, mark=+,
              mark size=3.50,mark options={solid,fill=blue, line width=0.70}]
              table[x expr = \thisrow{order}, 
                    y expr = \scale*\thisrow{\dataD}]
              {\datafile}; \addlegendentry{\legendD}
     }
     \ifthenelse{\equal{\dataC}{}}
     {}
     {
     \addplot[color=\mygreen,densely dashed, line width=0.2, mark=star,
              mark size=4.50,mark options={solid,fill=\mygreen, line width=.75}]
              table[x expr = \thisrow{order}, 
                    y expr = \scale*\thisrow{\dataC}]
              {\datafile}; \addlegendentry{\legendC}
     }

     \ifthenelse{\equal{\dataE}{}}
     {}
     {
     \addplot[color=magenta,densely dotted, line width=0.2, mark=triangle,
              mark size=3.50,mark options={solid,fill=magenta, line width=0.50}]
              table[x expr = \thisrow{order}, 
                    y expr = \scale*\thisrow{\dataE}]
              {\datafile}; \addlegendentry{\legendE}
     }

     \ifthenelse{\equal{\dataB}{}}
     {}
     {
     \addplot[color=red,densely dashed, line width=0.2, mark=o,
              mark size=2.50,mark options={solid,fill=\myred, line width=.60}]
              table[x expr = \thisrow{order}, 
                    y expr = \scale*\thisrow{\dataB}]
              {\datafile}; \addlegendentry{\legendB}
     }
%%%%%%%%%%%%%%%%%%%%%%%%%%%%%%%%%%%%%%%%%%%%%%%%%%%%%%%%%%%%%%%%%%%%%%%%%%%%%%%%
\end{axis}
%%%%%%%%%%%%%%%%%%%%%%%%%%%%%%%%%%%%%%%%%%%%%%%%%%%%%%%%%%%%%%%%%%%%%%%%%%%%%%%%
\end{tikzpicture}
             \label{figure:2D-solar:Su-error-freq-order_01b}}} 
  \caption{Relative error $\errall$ of \ref{eq:relative-error_sum} with 
           frequency (left) and polynomial order (right) for heterogeneous 
           elastic isotropic medium, % using HDG discretization of formulation $\formulationU$ and 
           with \flo{the} stabilization coefficients of \cref{eq:numerics:stabilization-u:2Dhe}.}
  \label{figure:2D-solar:Su-error-freq-order_01}

\end{figure}
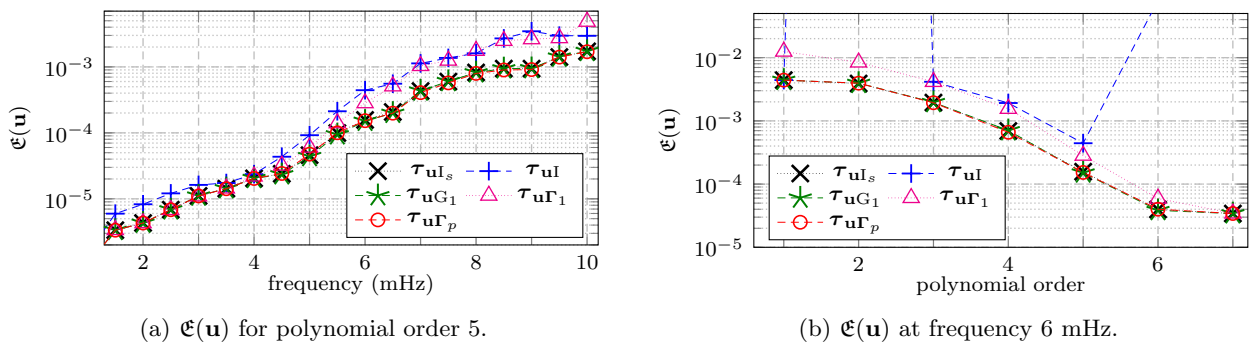
\subsubsection{Elastic tilted transverse isotropic medium}
\label{subsection:TTI-2D-helio}
% anisotropy
\renewcommand{\plottauSuIda}{tau9}
\renewcommand{\plottauSuIdb}{tau5}
\renewcommand{\plottauSuKCa}{tau17}
\renewcommand{\plottauSuKCb}{tau18}  % tau15 is using s-wavespeed
\renewcommand{\plottauSuMGa}{tau23}
\renewcommand{\plottauSuMGb}{tau33}
% ------------------------------------------------------------------------------- %

We now consider the propagation of waves in 
a tilted transverse isotropic (TTI) medium, 
which is a rotated version of a VTI one \flo{by an angle $\theta$}, % rotating the symmetry axis, 
cf.~\cite[Section 3.5.3]{Pham2023hdgRR}.
\flo{The} anisotropy is represented with constant Thomsen's parameters % such that
$\epsilon \,=\, \num{0.25}$, $\delta=\num{0.15}$, and % we use an angle 
$\theta\,=\, 45\si{\degree}$.
Note that compared to \cite[Section 3.5.3]{Pham2023hdgRR} written for 3D, 
the 2D case only has one angle, and no anisotropic coefficient $\gamma$.
Here, the heterogeneous solar-like background velocity ${c}$ 
and density $\rho$ are used to define 
$\mu^{\mathrm{TI}}=(\num{0.70}\mathrm{c})^2 \rho$
and 
$\lambda^{\mathrm{TI}}=\mathrm{c}^2 \rho - 2\mu^{\mathrm{TI}}$.
In \cref{figure:2D-TTI-solar:waves}, we show the reference solution at 
frequency 4 \si{\milli\Hz}, computed on a mesh with \num{80000} triangles
(while following simulations use \num{50000} elements), with polynomial 
approximation of order 7
and Godunov stabilization.
Compared to the isotropic case of \cref{figure:2D-solar:waves}, the 
solution is less smooth and contains more ripples.
Note that for the TTI absorbing boundary conditions, we follow \cite{bonnasse2018hybridizable}.
It is out of the scope to discuss the accuracy of the anisotropic boundary conditions here.

% ------------------------------
\setlength{\modelwidth} {4.10cm}
\setlength{\modelheight}{4.10cm}
% ------------------------------
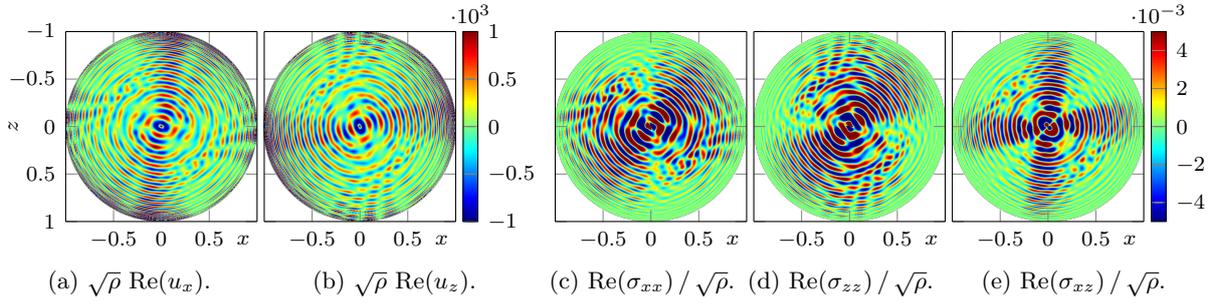
\begin{figure}[ht!] \centering
  \pgfkeys{/pgf/fpu=true} \pgfmathsetmacro{\cmin}{-1} \pgfmathsetmacro{\cmax}{1} \pgfkeys{/pgf/fpu=false}
  \renewcommand{\cbtitle}{\scriptsize $\cdot$\num{e3}}
  \renewcommand{\modelfile}{figures/numerics/2D-TTI_solar/wavefield/UxsqrtRho_4mHz_real_tau23p7_scale1e3}
  \subfloat[][$\sqrt{\rho}\,\, \Real(u_x)$.]
             {\begin{tikzpicture}
\pgfmathsetmacro{\xminloc}{-1.00} 
\pgfmathsetmacro{\xmaxloc}{ 1.00}
\pgfmathsetmacro{\yminloc}{-1.00} 
\pgfmathsetmacro{\ymaxloc}{ 1.00}
\pgfmathsetmacro{\xmin}{-1} 
\pgfmathsetmacro{\xmax}{ 1} 
\pgfmathsetmacro{\ymin}{-1} 
\pgfmathsetmacro{\ymax}{ 1}

\begin{axis}[%
  width=\modelwidth, height=\modelheight,
  axis on top, separate axis lines,
  xmin=\xminloc, xmax=\xmaxloc, 
  ymin=\yminloc, ymax=\ymaxloc, 
  xlabel={$x$},
  ylabel={$z$},
  xtick={-.5,0,0.5}, 
  y dir=reverse,
ylabel style = {yshift=-1em,  xshift=0mm},
xlabel style = {yshift =1.10em, xshift=2.80em},
%colormap/jet,colorbar,
%colorbar style={title={\cbtitle},
%width=.5em,xshift=-.80em,% ymode=log
%},
%point meta min=\cmin,point meta max=\cmax,
label style={font=\scriptsize},
tick label style={font=\scriptsize},
legend style={font=\scriptsize\selectfont},
]
\addplot [forget plot] graphics [xmin=\xmin,xmax=\xmax,
                                 ymin=\ymin,ymax=\ymax] 
                                 {{\modelfile}.png};
\end{axis}
\end{tikzpicture}%
} \hspace*{-1.0em}
  \renewcommand{\modelfile}{figures/numerics/2D-TTI_solar/wavefield/UzsqrtRho_4mHz_real_tau23p7_scale1e3}
  \subfloat[][$\sqrt{\rho}\,\, \Real(u_z)$.]
             {\begin{tikzpicture}
\pgfmathsetmacro{\xminloc}{-1.00} 
\pgfmathsetmacro{\xmaxloc}{ 1.00}
\pgfmathsetmacro{\yminloc}{-1.00} 
\pgfmathsetmacro{\ymaxloc}{ 1.00}
\pgfmathsetmacro{\xmin}{-1} 
\pgfmathsetmacro{\xmax}{ 1} 
\pgfmathsetmacro{\ymin}{-1} 
\pgfmathsetmacro{\ymax}{ 1}

\begin{axis}[%
  width=\modelwidth, height=\modelheight,
  axis on top, separate axis lines,
  xmin=\xminloc, xmax=\xmaxloc, 
  ymin=\yminloc, ymax=\ymaxloc, 
  xlabel={$x$},
  xtick={-.5,0,0.5},
  y dir=reverse,
  ytick={},
  yticklabels={,,},
  xlabel style = {yshift =1.10em, xshift=2.80em},
xlabel style = {yshift =0em, xshift=0mm},
colormap/jet,colorbar,
colorbar style={title={\cbtitle},
width=.5em,xshift=-0.60em,% ymode=log
title style={yshift=-.60em},
},
point meta min=\cmin,point meta max=\cmax,
label style={font=\scriptsize},
tick label style={font=\scriptsize},
legend style={font=\scriptsize\selectfont},
]
\addplot [forget plot] graphics [xmin=\xmin,xmax=\xmax,
                                 ymin=\ymin,ymax=\ymax] 
                                 {{\modelfile}.png};
\end{axis}
\end{tikzpicture}%
} \hspace*{-2.0em}
  \renewcommand{\cbtitle}{\scriptsize $\cdot$\num{e-3}}
  \pgfkeys{/pgf/fpu=true} \pgfmathsetmacro{\cmin}{-5} \pgfmathsetmacro{\cmax}{5} \pgfkeys{/pgf/fpu=false}
  \renewcommand{\modelfile}{figures/numerics/2D-TTI_solar/wavefield/Sxx-over-sqrtRho_4mHz_real_tau23p7_scale5e-3} 
  \subfloat[][$\Real(\sigma_{xx}) \,/\, \sqrt{\rho}$.]
             {\begin{tikzpicture}
\pgfmathsetmacro{\xminloc}{-1.00} 
\pgfmathsetmacro{\xmaxloc}{ 1.00}
\pgfmathsetmacro{\yminloc}{-1.00} 
\pgfmathsetmacro{\ymaxloc}{ 1.00}
\pgfmathsetmacro{\xmin}{-1} 
\pgfmathsetmacro{\xmax}{ 1} 
\pgfmathsetmacro{\ymin}{-1} 
\pgfmathsetmacro{\ymax}{ 1}

\begin{axis}[%
  width=\modelwidth, height=\modelheight,
  axis on top, separate axis lines,
  xmin=\xminloc, xmax=\xmaxloc, 
  ymin=\yminloc, ymax=\ymaxloc, 
  xlabel={$x$},
  xtick={-.5,0,0.5},
  y dir=reverse,
  ytick={},
  yticklabels={,,},
  xlabel style = {yshift =1.10em, xshift=2.80em},
% ytick={\ymin,\ymax},
% xtick={\xmin,1,\xmax},
% xticklabels={\num{\xmin},1,\num{\xmax}},
% yticklabels={\num{\ymin},\ymax},
% ylabel style = {yshift=-1em, xshift=0mm},
%xlabel style = {yshift =0mm, xshift=0mm},
%colormap/jet,colorbar,
%colorbar style={title={\cbtitle},
%width=.5em,xshift=-0.20em,% ymode=log
%},
%point meta min=\cmin,point meta max=\cmax,
label style={font=\scriptsize},
tick label style={font=\scriptsize},
legend style={font=\scriptsize\selectfont},
]
\addplot [forget plot] graphics [xmin=\xmin,xmax=\xmax,
                                 ymin=\ymin,ymax=\ymax] 
                                 {{\modelfile}.png};
\end{axis}
\end{tikzpicture}%
} \hspace*{-1.0em}
  \renewcommand{\modelfile}{figures/numerics/2D-TTI_solar/wavefield/Szz-over-sqrtRho_4mHz_real_tau23p7_scale5e-3} 
  \subfloat[][$\Real(\sigma_{zz}) \,/\, \sqrt{\rho}$.]
             {\begin{tikzpicture}
\pgfmathsetmacro{\xminloc}{-1.00} 
\pgfmathsetmacro{\xmaxloc}{ 1.00}
\pgfmathsetmacro{\yminloc}{-1.00} 
\pgfmathsetmacro{\ymaxloc}{ 1.00}
\pgfmathsetmacro{\xmin}{-1} 
\pgfmathsetmacro{\xmax}{ 1} 
\pgfmathsetmacro{\ymin}{-1} 
\pgfmathsetmacro{\ymax}{ 1}

\begin{axis}[%
  width=\modelwidth, height=\modelheight,
  axis on top, separate axis lines,
  xmin=\xminloc, xmax=\xmaxloc, 
  ymin=\yminloc, ymax=\ymaxloc, 
  xlabel={$x$},
  xtick={-.5,0,0.5},
  y dir=reverse,
  ytick={},
  yticklabels={,,},
  xlabel style = {yshift =1.10em, xshift=2.80em},
% ytick={\ymin,\ymax},
% xtick={\xmin,1,\xmax},
% xticklabels={\num{\xmin},1,\num{\xmax}},
% yticklabels={\num{\ymin},\ymax},
% ylabel style = {yshift=-1em, xshift=0mm},
%xlabel style = {yshift =0mm, xshift=0mm},
%colormap/jet,colorbar,
%colorbar style={title={\cbtitle},
%width=.5em,xshift=-0.20em,% ymode=log
%},
%point meta min=\cmin,point meta max=\cmax,
label style={font=\scriptsize},
tick label style={font=\scriptsize},
legend style={font=\scriptsize\selectfont},
]
\addplot [forget plot] graphics [xmin=\xmin,xmax=\xmax,
                                 ymin=\ymin,ymax=\ymax] 
                                 {{\modelfile}.png};
\end{axis}
\end{tikzpicture}%
} \hspace*{-1.0em}
  \renewcommand{\modelfile}{figures/numerics/2D-TTI_solar/wavefield/Sxz-over-sqrtRho_4mHz_real_tau23p7_scale5e-3} 
  \subfloat[][$\Real(\sigma_{xz}) \,/\, \sqrt{\rho}$.]
             {\begin{tikzpicture}
\pgfmathsetmacro{\xminloc}{-1.00} 
\pgfmathsetmacro{\xmaxloc}{ 1.00}
\pgfmathsetmacro{\yminloc}{-1.00} 
\pgfmathsetmacro{\ymaxloc}{ 1.00}
\pgfmathsetmacro{\xmin}{-1} 
\pgfmathsetmacro{\xmax}{ 1} 
\pgfmathsetmacro{\ymin}{-1} 
\pgfmathsetmacro{\ymax}{ 1}

\begin{axis}[%
  width=\modelwidth, height=\modelheight,
  axis on top, separate axis lines,
  xmin=\xminloc, xmax=\xmaxloc, 
  ymin=\yminloc, ymax=\ymaxloc, 
  xlabel={$x$},
  xtick={-.5,0,0.5},
  y dir=reverse,
  ytick={},
  yticklabels={,,},
  xlabel style = {yshift =1.10em, xshift=2.80em},
xlabel style = {yshift =0em, xshift=0mm},
colormap/jet,colorbar,
colorbar style={title={\cbtitle},
width=.5em,xshift=-0.60em,% ymode=log
title style={yshift=-.60em},
},
point meta min=\cmin,point meta max=\cmax,
label style={font=\scriptsize},
tick label style={font=\scriptsize},
legend style={font=\scriptsize\selectfont},
]
\addplot [forget plot] graphics [xmin=\xmin,xmax=\xmax,
                                 ymin=\ymin,ymax=\ymax] 
                                 {{\modelfile}.png};
\end{axis}
\end{tikzpicture}%
}

  \caption{
           Reference solutions for 
           the elastic TTI wave equation \flo{in a solar-like background  model}
           at \num{4} \si{\milli\Hz}
           obtained with stabilization $\tauSuMGunit$ and polynomial approximation 
           of order 7.
           \flo{The background models are represented 
                by piecewise polynomials on each cell with Lagrange basis \modif{of} order 2}.
          }
          
  \label{figure:2D-TTI-solar:waves}
       
\end{figure}

In \cref{figure:2D-TTI-solar:Su-error-freq-order_02}, 
we plot the relative error $\errall$ depending on the 
frequency and polynomial orders. 
%% for the HDG formulation $\formulationU$. 
Similar to the isotropic case, stabilization $\tauSuIda$ 
sometimes results in a\flo{n} ill-conditioned matrix leading to
a linear system that cannot be solved.
The results confirm the behaviour observed in the isotropic
case: the Godunov stabilization gives the most accurate
results, and the level \flo{of} accuracy can be met with other stabilization
involving a well-chosen scaling factor.
% Namely, % it appears that 
% Thus, considering transverse isotropy does not modify the behaviour 
% and efficiency of the discretization scheme.

% ---------------------------------------------------------------------- %
\setlength{\plotwidth} {6.50cm}
\setlength{\plotheight}{3.10cm}
\renewcommand{\myfreqplot} {0.005}
\renewcommand{\myorderplot}{4}
% ----------------------------
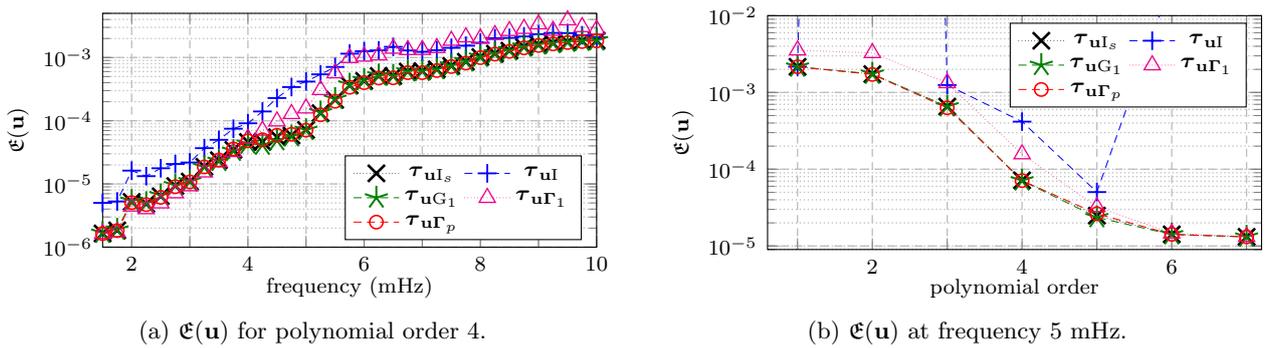
\begin{figure}[ht!] \centering

  \renewcommand{\myxlabel}{frequency (\si{\milli\Hz})}
  \pgfkeys{/pgf/fpu=true}
  \pgfmathsetmacro{\xmin}{1.50} \pgfmathsetmacro{\xmax}{10} 
  \pgfmathsetmacro{\ymin}{1e-6} \pgfmathsetmacro{\ymax}{5e-3} 
  \pgfmathsetmacro{\scalefreq}{1e3} % cs/freq
  \pgfkeys{/pgf/fpu=false}
  \renewcommand{\datafile}{figures/numerics/2D-TTI_solar/data/error-scaled_mean_per-freq_u-and-sigmas.txt}

  \renewcommand{\dataA}{sumU_\plottauSuIdb_p\myorderplot} \renewcommand{\legendA}{$\tauSuIdb$} 
  \renewcommand{\dataB}{sumU_\plottauSuKCb_p\myorderplot} \renewcommand{\legendB}{$\tauSuKCb$}
  \renewcommand{\dataC}{sumU_\plottauSuMGa_p\myorderplot} \renewcommand{\legendC}{$\tauSuMGunit$}
  \renewcommand{\dataD}{sumU_\plottauSuIda_p\myorderplot} \renewcommand{\legendD}{$\tauSuIda$} 
  \renewcommand{\dataE}{sumU_\plottauSuKCa_p\myorderplot} \renewcommand{\legendE}{$\tauSuKCunit$}
  \renewcommand{\legendpos}{south east}
  \renewcommand{\myylabel}{$\errall(\displacement)$}
  \subfloat[][$\errall(\displacement)$ for polynomial order \myorderplot.]
             {\makebox[.35\linewidth][c]{%%%%%%%%%%%%%%%%%%%%%%%%%%%%%%%%%%%%%%%%%%%%%%%%%%%%%%%%%%%%%%%%%%%%%%%%%%%%%%%%
%% begin the TikZ picture
\begin{tikzpicture}
%%%%%%%%%%%%%%%%%%%%%%%%%%%%%%%%%%%%%%%%%%%%%%%%%%%%%%%%%%%%%%%%%%%%%%%%%%%%%%%%
\begin{axis}[
             enlargelimits=false, 
             xlabel={\scriptsize \myxlabel },
             ylabel={\scriptsize \myylabel },
             enlarge y limits=false,
             enlarge x limits=false,
             xmin=\xmin,xmax=\xmax, 
             ymin=\ymin,ymax=\ymax,
             yminorticks=true,
             height=\plotheight,width=\plotwidth,
             scale only axis,
             ylabel style = {yshift=-0.2em, xshift=0mm},
             xlabel style = {yshift =0.5em, xshift=0mm},
             ytick={1e-6,1e-5,1e-4,1e-3,1e-2,1e-1},
             xtick={0,1,2,3,4,5,6,7,8,9,10,11,12,13,14,15,16,17,18,19,20},
             xticklabels={,,2,,4,,6,,8,,10,,12,,14,,16,,18,,20},
             clip mode=individual,
             minor grid style={line width=.5pt, draw=gray!50, densely dotted},
             major grid style={line width=.5pt, draw=gray!50, densely dashed},
             xmajorgrids=true,
             xminorgrids=true,
             yminorgrids=true,
             ymajorgrids=true,
             ymode=log,
            %xmode=log,
             %y  tick label style={/pgf/number format/fixed},
             %y  tick label style={/pgf/number format/fixed zerofill},
             %y  tick label style={/pgf/number format/precision=1},
             %y  tick label style={/pgf/number format/precision=2},
             % y  tick label style={/pgf/number format/sci},
             %y tick label style={/pgf/number format/.cd},
             label style={font=\scriptsize},
             tick label style={font=\scriptsize},
             legend style={font=\scriptsize\selectfont},
             legend pos={\legendpos}, 
             legend columns=2,
             ]  %% foreground marks

     \pgfmathsetmacro{\scale}{1}

     %% load current data
     \ifthenelse{\equal{\dataA}{}}
     {}
     {
     \addplot[color=black,densely dotted,line width=0.2, mark=x, 
              mark size=4.50,mark options={solid,fill=black, line width=1.0}]
              table[x expr = \scalefreq*\thisrow{freq}, 
                    y expr = \scale*\thisrow{\dataA}, % select coords between index={1}{5}
                   ]
              {\datafile}; \addlegendentry{\legendA}    
     }
     \ifthenelse{\equal{\dataD}{}}
     {}
     {
     \addplot[color=blue,densely dashed, line width=0.2, mark=+,
              mark size=3.50,mark options={solid,fill=blue, line width=0.70}]
              table[x expr = \scalefreq*\thisrow{freq}, 
                    y expr = \scale*\thisrow{\dataD}]
              {\datafile}; \addlegendentry{\legendD}
     }

     \ifthenelse{\equal{\dataC}{}}
     {}
     {
     \addplot[color=\mygreen,densely dashed, line width=0.2, mark=star,
              mark size=4.50,mark options={solid,fill=\mygreen, line width=.75}]
              table[x expr = \scalefreq*\thisrow{freq}, 
                    y expr = \scale*\thisrow{\dataC}]
              {\datafile}; \addlegendentry{\legendC}
     }

     \ifthenelse{\equal{\dataE}{}}
     {}
     {
     \addplot[color=magenta,densely dotted, line width=0.2, mark=triangle,
              mark size=3.50,mark options={solid,fill=magenta, line width=0.50}]
              table[x expr = \scalefreq*\thisrow{freq}, 
                    y expr = \scale*\thisrow{\dataE}]
              {\datafile}; \addlegendentry{\legendE}
     }
     
     \ifthenelse{\equal{\dataB}{}}
     {}
     {
     \addplot[color=red,densely dashed, line width=0.2, mark=o,
              mark size=2.50,mark options={solid,fill=\myred, line width=.60}]
              table[x expr = \scalefreq*\thisrow{freq}, 
                    y expr = \scale*\thisrow{\dataB}]
              {\datafile}; \addlegendentry{\legendB}
     }

%%%%%%%%%%%%%%%%%%%%%%%%%%%%%%%%%%%%%%%%%%%%%%%%%%%%%%%%%%%%%%%%%%%%%%%%%%%%%%%%
\end{axis}
%%%%%%%%%%%%%%%%%%%%%%%%%%%%%%%%%%%%%%%%%%%%%%%%%%%%%%%%%%%%%%%%%%%%%%%%%%%%%%%%
\end{tikzpicture}}}
              \hspace*{.50em}
  \pgfkeys{/pgf/fpu=true}
  \pgfmathsetmacro{\ymin}{9e-6} \pgfmathsetmacro{\ymax}{1e-2} 
  \pgfmathsetmacro{\xmin}{0.60} \pgfmathsetmacro{\xmax}{7.20} 
  \pgfmathsetmacro{\scalefreq}{1e3} % cs/freq
  \pgfkeys{/pgf/fpu=false}
  \renewcommand{\myxlabel}{polynomial order}
  \renewcommand{\legendpos}{north east}
  \renewcommand{\datafile}{figures/numerics/2D-TTI_solar/data/error-scaled_mean_per-order_u-and-sigmas.txt}
  \renewcommand{\dataA}{sumU_\plottauSuIdb_freq\myfreqplot} \renewcommand{\legendA}{$\tauSuIdb$} 
  \renewcommand{\dataB}{sumU_\plottauSuKCb_freq\myfreqplot} \renewcommand{\legendB}{$\tauSuKCb$}
  \renewcommand{\dataC}{sumU_\plottauSuMGa_freq\myfreqplot} \renewcommand{\legendC}{$\tauSuMGunit$}
  \renewcommand{\dataD}{sumU_\plottauSuIda_freq\myfreqplot} \renewcommand{\legendD}{$\tauSuIda$}
  \renewcommand{\dataE}{sumU_\plottauSuKCa_freq\myfreqplot} \renewcommand{\legendE}{$\tauSuKCunit$}   
  \renewcommand{\myylabel}{$\errall(\displacement)$}
  \subfloat[][$\errall(\displacement)$ at frequency 5~\si{\milli\Hz}.]
             {\makebox[.35\linewidth][c]{%%%%%%%%%%%%%%%%%%%%%%%%%%%%%%%%%%%%%%%%%%%%%%%%%%%%%%%%%%%%%%%%%%%%%%%%%%%%%%%%
%% begin the TikZ picture
\begin{tikzpicture}
%%%%%%%%%%%%%%%%%%%%%%%%%%%%%%%%%%%%%%%%%%%%%%%%%%%%%%%%%%%%%%%%%%%%%%%%%%%%%%%%
\begin{axis}[
             enlargelimits=false, 
             xlabel={\scriptsize \myxlabel },
             ylabel={\scriptsize \myylabel },
             enlarge y limits=false,
             enlarge x limits=false,
             xmin=\xmin,xmax=\xmax, 
             ymin=\ymin,ymax=\ymax,
             yminorticks=true,
             height=\plotheight,width=\plotwidth,
             scale only axis,
             ylabel style = {yshift=-0.2em, xshift=0mm},
             xlabel style = {yshift =0.5em, xshift=0mm},
             ytick={1e-6,1e-5,1e-4,1e-3,1e-2,1e-1},
             xtick={0,1,2,3,4,5,6,7,8,9,10,11,12,13,14,15,16,17,18,19,20},
             xticklabels={,,2,,4,,6,,8,,10,,12,,14,,16,,18,,20},
             clip mode=individual,
             minor grid style={line width=.5pt, draw=gray!50, densely dotted},
             major grid style={line width=.5pt, draw=gray!50, densely dashed},
             xmajorgrids=true,
             xminorgrids=true,
             yminorgrids=true,
             ymajorgrids=true,
             ymode=log,
             %y  tick label style={/pgf/number format/fixed},
             %y  tick label style={/pgf/number format/fixed zerofill},
             %y  tick label style={/pgf/number format/precision=1},
             %y  tick label style={/pgf/number format/precision=2},
             % y  tick label style={/pgf/number format/sci},
             %y tick label style={/pgf/number format/.cd},
             label style={font=\scriptsize},
             tick label style={font=\scriptsize},
             legend style={font=\scriptsize\selectfont},
             legend pos={\legendpos}, 
             legend columns=2,
             ]  %% foreground marks

     \pgfmathsetmacro{\scale}{1}

     %% load current data
     \ifthenelse{\equal{\dataA}{}}
     {}
     {
     \addplot[color=black,densely dotted,line width=0.2, mark=x, 
              mark size=4.50,mark options={solid,fill=black, line width=1.0}]
              table[x expr = \thisrow{order}, 
                    y expr = \scale*\thisrow{\dataA}, % select coords between index={1}{5}
                   ]
              {\datafile}; \addlegendentry{\legendA}    
     }

     \ifthenelse{\equal{\dataD}{}}
     {}
     {
     \addplot[color=blue,densely dashed, line width=0.2, mark=+,
              mark size=3.50,mark options={solid,fill=blue, line width=0.70}]
              table[x expr = \thisrow{order}, 
                    y expr = \scale*\thisrow{\dataD}]
              {\datafile}; \addlegendentry{\legendD}
     }
     \ifthenelse{\equal{\dataC}{}}
     {}
     {
     \addplot[color=\mygreen,densely dashed, line width=0.2, mark=star,
              mark size=4.50,mark options={solid,fill=\mygreen, line width=.75}]
              table[x expr = \thisrow{order}, 
                    y expr = \scale*\thisrow{\dataC}]
              {\datafile}; \addlegendentry{\legendC}
     }

     \ifthenelse{\equal{\dataE}{}}
     {}
     {
     \addplot[color=magenta,densely dotted, line width=0.2, mark=triangle,
              mark size=3.50,mark options={solid,fill=magenta, line width=0.50}]
              table[x expr = \thisrow{order}, 
                    y expr = \scale*\thisrow{\dataE}]
              {\datafile}; \addlegendentry{\legendE}
     }

     \ifthenelse{\equal{\dataB}{}}
     {}
     {
     \addplot[color=red,densely dashed, line width=0.2, mark=o,
              mark size=2.50,mark options={solid,fill=\myred, line width=.60}]
              table[x expr = \thisrow{order}, 
                    y expr = \scale*\thisrow{\dataB}]
              {\datafile}; \addlegendentry{\legendB}
     }
%%%%%%%%%%%%%%%%%%%%%%%%%%%%%%%%%%%%%%%%%%%%%%%%%%%%%%%%%%%%%%%%%%%%%%%%%%%%%%%%
\end{axis}
%%%%%%%%%%%%%%%%%%%%%%%%%%%%%%%%%%%%%%%%%%%%%%%%%%%%%%%%%%%%%%%%%%%%%%%%%%%%%%%%
\end{tikzpicture}}} 
  \caption{Relative error $\errall$ of \cref{eq:relative-error_sum} with 
           frequency (left) and the polynomial order (right) for heterogeneous 
           elastic TTI medium, % using HDG discretization of formulation $\formulationU$ and 
           with stabilization coefficients of \cref{eq:numerics:stabilization-u:2Dhe}.}
  \label{figure:2D-TTI-solar:Su-error-freq-order_02}

\end{figure}
% ---------------------------------------------------------------------- %

% ---------------------------------------------------------------------- %
\subsection{Summary of numerical experiments}
\label{subsection:numerics-summary}
% ---------------------------------------------------------------------- %

In our first experiments working with planewaves, 
we identify \flo{the} optimal scaling factor for each family 
of stabilization in \cref{eq:numerics:stabilization-u:3Dpw}.
We find that the Godunov stabilization is the most robust as it 
does not need scaling and $\tau=1$ gives the best results for
all types of waves. 
On the other hand, the Kelvin--Christoffel and identity-based
families need a suitable scaling factor in order to reach 
the same level of accuracy obtained with the Godunov matrix. 
The optimal value for the scaling does not seem to change for the 
Kelvin--Christoffel stabilization and leans towards the P-slowness.
However, \flo{the} identity-based stabilization shows a different optimal
scaling factor depending on the type of waves propagating.
Our experiments also compare the performance of the stabilizations 
in both isotropy and anisotropy, \flo{for which the same conclusions are obtained}.

In our second experiment, \flo{we consider waves generated by a point source 
in a heterogeneous medium, thus a superposition of different types of waves.}
\flo{In this case}, the Godunov stabilization is the most accurate. 
\flo{The identity-based stabilization with the S-impedance scaling factor and 
and the Kelvin-Christoffel stabilization with the P-slowness also yield comparable accuracy.}
\flo{In} the identity-based family, the better performance of \flo{the} S-impedance 
compared to \flo{the} P-impedance can be explained \flo{by} the equipartition 
phenomenon, \cite{shapiro2000partitioning,Snieder2002equilibration}, 
which says that S-waves are more energetic than P-waves.
\flo{In another word}, although all types of waves are present, the magnitude of 
S-waves is stronger, hence stabilization with S-impedance provides higher accuracy.

As  a common feature among the optimal form from each family, 
the components of the stabilization operators take \flo{the} value of 
an impedance (i.e., density multiplied by wavespeed). 
This is evident for the identity-based and Godunov stabilization
(cf. \cref{eq:MGodunov-iso,eq:MGodunov-aniso}).
For the Kelvin--Christoffel stabilization, the components of the 
matrix $\Gambold$ \flo{take the form of} $\rho \mathrm{c}^2$ \cref{KC_iso}, and 
the optimal scaling $\tau$ is a slowness, hence $\tau \Gambold$
has a magnitude of an impedance.

% The experiments also compared elastic isotropy and anisotropy, 
% a property that does not appear to create drastic differences 
% in the results. 

% ------------------------------------------------------------------------------- %
\section{Conclusion}
% ------------------------------------------------------------------------------- %

% In this work, we have investigated the selection of an efficient 
% stabilization of HDG discretization for elastic wave equations. 
% We have provided mathematical analysis to derive different choices
% of stabilization, for isotropic and anisotropic media.
% In particular we define the Godunov stabilization for anisotropic 
% elasticity. We have compared and evaluated the efficiency of 
% several stabilization coefficients and the Godunov stabilization 
% is shown to be the most versatile choice, giving the best accuracy.
% While similar accuracy can be obtained with an identity-based 
% stabilization, it requires some tuning coefficient and one cannot 
% find a universal value, contrary to the Godunov stabilization 
% which is entirely defined from the elastic properties.
% Furthermore, we have used the elastic formulation based upon the 
% compliance tensor that allows us to easily account for physical 
% properties that vary within a cell, hence to accurately represent 
% models having drastic variation, as highlighted by our implementation. 
% It also gives us the indispensable flexibility to now consider 
% the quantitative inverse problem for the reconstruction of elastic 
% properties.
% Future works also include the derivation of accurate absorbing 
% boundary conditions in elastic anisotropy.

In this work, we employed Voigt notation in the HDG method
to describe compactly the discrete problem for anisotropic elasticity.
Additionally, a first-order formulation with \flo{the} compliance 
tensor is used and allows for mesh-wise variation of parameters.
This, together with Voigt's notation which provides 
efficient bookkeeping of physical parameters,
% form indispensable features 
\flo{will be indispensable}
in quantitative reconstruction of
elastic parameters in inverse problems.
Secondly, to determine an optimal choice of stabilization, 
we constructed the hybridized Godunov-upwind flux
for anisotropic elasticity, which
offers a versatile choice and removes the need for 
scaling factor tuning. 
% This problem concerns in particular \flo{the} identity-based 
% stabilization, a popular choice due to its simple form, however, 
% one which lacks a universal scaling factor.
\flo{It is worth noting that the identity-based stabilization, 
which is very popular, lacks a universal scaling factor, making 
it less robust for wave simulation in a general setting.}
We have carried out numerical experiments in two and three
dimensions, considering isotropic elasticity and anisotropy, 
with constant backgrounds as well as one containing high variation 
and contrast.
They demonstrate the performance and versatility of the Godunov 
stabilization, which is well suited for generic anisotropic 
material and different types of waves.

% ------------------------
\section*{Acknowledgments}

  % We acknowledge the use of computational 
  % resources from CINES--Adastra with 
  % HPC Grand Challenge allocation gda2306.
  This project was provided with computer and storage resources 
  by GENCI at CINES thanks to the grant \texttt{gda2306} on the 
  supercomputer Adastra's GENOA partition.
  This work was partially supported by the EXAMA (Methods and Algorithms at Exascale) 
  project under grant ANR--22--EXNU--0002.
  FF acknowledges funding by the European Union 
  with ERC Project \textsc{Incorwave} -- grant 101116288. 
  Views and opinions expressed are however those of the author(s) 
  only and do not necessarily reflect those 
  of the European Union or the European Research Council 
  Executive Agency (ERCEA). Neither the European Union nor the
  granting authority can be held responsible for them.

% ---------------------------------------------------------------------

% \appendix
% \input{sections/appendix_green_kernel}
% \input{sections/appendix_planewave}
% \input{sections/appendix_scalarflux}
% \input{sections/05-stabilization}

% -----------------------------------
% Bibliography
% -----------------------------------
% \newpage \footnotesize 
\bibliographystyle{siamplain}
\bibliography{sections/bibliography}

\begin{thebibliography}{10}

\bibitem{arnold2002unified}
{\sc D.~N. Arnold, F.~Brezzi, B.~Cockburn, and L.~D. Marini}, {\em Unified
  analysis of discontinuous {G}alerkin methods for elliptic problems}, SIAM
  journal on numerical analysis, 39 (2002), pp.~1749--1779.

\bibitem{Pham2020Siam}
{\sc H.~Barucq, F.~Faucher, D.~Fournier, L.~Gizon, and H.~Pham}, {\em Efficient
  and accurate algorithm for the full modal {G}reen's kernel of the scalar wave
  equation in helioseismology}, SIAM Journal on Applied Mathematics, 80 (2020),
  pp.~2657--2683.

\bibitem{Pham2019Esaim}
{\sc H.~Barucq, F.~Faucher, and H.~Pham}, {\em Outgoing solutions and radiation
  boundary conditions for the ideal atmospheric scalar wave equation in
  helioseismology}, ESAIM: Mathematical Modelling and Numerical Analysis, 54
  (2020), pp.~1111--1138.

\bibitem{bonnasse2018hybridizable}
{\sc M.~Bonnasse-Gahot, H.~Calandra, J.~Diaz, and S.~Lanteri}, {\em
  Hybridizable discontinuous {G}alerkin method for the 2-{D} frequency-domain
  elastic wave equations}, Geophysical Journal International, 213 (2018),
  pp.~637--659.

\bibitem{bui2015godunov}
{\sc T.~Bui-Thanh}, {\em From {G}odunov to a unified hybridized discontinuous
  {G}alerkin framework for partial differential equations}, Journal of Comp.
  Physics, 295 (2015), pp.~114--146.

\bibitem{bui2015rankine}
{\sc T.~Bui-Thanh}, {\em From {R}ankine-{H}ugoniot condition to a constructive
  derivation of {HDG} methods}, in Spectral and High Order Methods for Partial
  Differential Equations ICOSAHOM 2014, Springer, 2015, pp.~483--491.

\bibitem{bui2016construction}
{\sc T.~Bui-Thanh}, {\em Construction and analysis of {HDG} methods for
  linearized shallow water equations}, SIAM Journal on Scientific Computing, 38
  (2016), pp.~A3696--A3719.

\bibitem{carcione2007wave}
{\sc J.~M. Carcione}, {\em Wave fields in real media: Wave propagation in
  anisotropic, anelastic, porous and electromagnetic media}, Elsevier, 2007.

\bibitem{Dalsgaard1996}
{\sc J.~Christensen-Dalsgaard, W.~D{\"a}ppen, S.~V. Ajukov, E.~R. Anderson,
  H.~M. Antia, S.~Basu, V.~A. Baturin, G.~Berthomieu, B.~Chaboyer, S.~M.
  Chitre, A.~N. Cox, P.~Demarque, J.~Donatowicz, W.~A. Dziembowski, M.~Gabriel,
  D.~O. Gough, D.~B. Guenther, J.~A. Guzik, J.~W. Harvey, F.~Hill, G.~Houdek,
  C.~A. Iglesias, A.~G. Kosovichev, J.~W. Leibacher, P.~Morel, C.~R. Proffitt,
  J.~Provost, J.~Reiter, E.~J. Rhodes, F.~J. Rogers, I.~W. Roxburgh, M.~J.
  Thompson, and R.~K. Ulrich}, {\em The current state of solar modeling},
  Science, 272 (1996), pp.~1286--1292.
\newblock 10.1126/science.272.5266.1286.

\bibitem{cockburn2016static}
{\sc B.~Cockburn}, {\em Static condensation, hybridization, and the devising of
  the {HDG} methods}, Building bridges: connections and challenges in modern
  approaches to numerical partial differential equations,  (2016),
  pp.~129--177.

\bibitem{cockburn2023hybridizable}
{\sc B.~Cockburn}, {\em Hybridizable discontinuous {G}alerkin methods for
  second-order elliptic problems: overview, a new result and open problems},
  Japan Journal of Industrial and Applied Mathematics,  (2023), pp.~1--40.

\bibitem{cockburn2016bridging}
{\sc B.~Cockburn, D.~A. Di~Pietro, and A.~Ern}, {\em Bridging the hybrid
  high-order and hybridizable discontinuous {G}alerkin methods}, ESAIM:
  Mathematical Modelling and Numerical Analysis, 50 (2016), pp.~635--650.

\bibitem{cockburn2008superconvergent}
{\sc B.~Cockburn, B.~Dong, and J.~Guzm{\'a}n}, {\em A superconvergent
  {LDG}-hybridizable {G}alerkin method for second-order elliptic problems},
  Math. of Comp., 77 (2008), pp.~1887--1916.

\bibitem{cockburn2009unified}
{\sc B.~Cockburn, J.~Gopalakrishnan, and R.~Lazarov}, {\em Unified
  hybridization of discontinuous {G}alerkin, mixed, and continuous {G}alerkin
  methods for second order elliptic problems}, SIAM Journal on Numerical
  Analysis, 47 (2009), pp.~1319--1365.

\bibitem{cockburn2013superconvergent}
{\sc B.~Cockburn and K.~Shi}, {\em Superconvergent {HDG} methods for linear
  elasticity with weakly symmetric stresses}, IMA Journal of Numerical
  Analysis, 33 (2013), pp.~747--770.

\bibitem{di2015hybrid}
{\sc D.~A. Di~Pietro and A.~Ern}, {\em A hybrid high-order locking-free method
  for linear elasticity on general meshes}, Comp. Meth. in App. Mechanics and
  Engineering, 283 (2015), pp.~1--21.

\bibitem{du2020new}
{\sc S.~Du and F.-J. Sayas}, {\em New analytical tools for {HDG} in elasticity,
  with applications to elastodynamics}, Mathematics of Computation, 89 (2020),
  pp.~1745--1782.

\bibitem{fabien2020gpu}
{\sc M.~S. Fabien}, {\em A {GPU}-accelerated hybridizable discontinuous
  {G}alerkin method for linear elasticity}, Commun Comput Phys, 27 (2020),
  pp.~513--545.

\bibitem{Hawen2021}
{\sc F.~Faucher}, {\em \texttt{hawen}: time-harmonic wave modeling and
  inversion using hybridizable discontinuous {G}alerkin discretization},
  Journal of Open Source Software, 6 (2021).

\bibitem{Faucher2020adjoint}
{\sc F.~Faucher and O.~Scherzer}, {\em Adjoint-state method for {H}ybridizable
  {D}iscontinuous {G}alerkin discretization, application to the inverse
  acoustic wave problem}, Computer Methods in Applied Mechanics and
  Engineering, 372 (2020), p.~113406.

\bibitem{fernandez2018hybridized}
{\sc P.~Fernandez, A.~Christophe, S.~Terrana, N.~C. Nguyen, and J.~Peraire},
  {\em Hybridized discontinuous {G}alerkin methods for wave propagation},
  Journal of Scientific Computing, 77 (2018), pp.~1566--1604.

\bibitem{fu2015analysis}
{\sc G.~Fu, B.~Cockburn, and H.~Stolarski}, {\em Analysis of an {HDG} method
  for linear elasticity}, International Journal for Numerical Methods in
  Engineering, 102 (2015), pp.~551--575.

\bibitem{giacomini2019discontinuous}
{\sc M.~Giacomini and R.~Sevilla}, {\em Discontinuous {G}alerkin approximations
  in computational mechanics: hybridization, exact geometry and degree
  adaptivity}, SN Applied Sciences, 1 (2019), p.~1047.

\bibitem{giacomini2021hdglab}
{\sc M.~Giacomini, R.~Sevilla, and A.~Huerta}, {\em {HDGlab}: An open-source
  implementation of the hybridisable discontinuous {G}alerkin method in
  matlab}, Archives of Computational Methods in Engineering, 28 (2021),
  pp.~1941--1986.

\bibitem{Givoli1990}
{\sc D.~Givoli and J.~B. Keller}, {\em Non-reflecting boundary conditions for
  elastic waves}, Wave motion, 12 (1990), pp.~261--279.

\bibitem{gopalakrishnan2015stabilization}
{\sc J.~Gopalakrishnan, S.~Lanteri, N.~Olivares, and R.~Perrussel}, {\em
  Stabilization in relation to wavenumber in {HDG} methods}, Advanced Modeling
  and Simulation in Engineering Sciences, 2 (2015), pp.~1--24.

\bibitem{Higdon1991}
{\sc R.~L. Higdon}, {\em Absorbing boundary conditions for elastic waves},
  Geophysics, 56 (1991), pp.~231--241.
\newblock 10.1190/1.1443035.

\bibitem{hungria2017hdg}
{\sc A.~Hungria, D.~Prada, and F.-J. Sayas}, {\em {HDG} methods for
  elastodynamics}, Computers \& Mathematics with Applications, 74 (2017),
  pp.~2671--2690.

\bibitem{nguyen2012hybridizable}
{\sc N.~C. Nguyen and J.~Peraire}, {\em Hybridizable discontinuous {G}alerkin
  methods for partial differential equations in continuum mechanics}, Journal
  of Computational Physics, 231 (2012), pp.~5955--5988.

\bibitem{nguyen2011high}
{\sc N.~C. Nguyen, J.~Peraire, and B.~Cockburn}, {\em High-order implicit
  hybridizable discontinuous {G}alerkin methods for acoustics and
  elastodynamics}, Journal of Computational Physics, 230 (2011),
  pp.~3695--3718.

\bibitem{Pham2023hdgRR}
{\sc H.~Pham, F.~Faucher, and H.~Barucq}, {\em {On the implementation of
  Hybridizable Discontinuous Galerkin discretization for linear anisotropic
  elastic wave equation: Voigt-notation and stabilization}}, Research Report
  RR-9533, {Inria Bordeaux Sud-Ouest; Project-Team Makutu}, December 2023.
\newblock https://hal.archives-ouvertes.fr/hal-04356602/file/RR-9533.pdf.

\bibitem{sevilla2018superconvergent}
{\sc R.~Sevilla, M.~Giacomini, A.~Karkoulias, and A.~Huerta}, {\em A
  superconvergent hybridisable discontinuous {G}alerkin method for linear
  elasticity}, International Journal for Numerical Methods in Engineering, 116
  (2018), pp.~91--116.

\bibitem{shapiro2000partitioning}
{\sc N.~M. Shapiro, M.~Campillo, L.~Margerin, S.~K. Singh, V.~Kostoglodov, and
  J.~Pacheco}, {\em {The Energy Partitioning and the Diffusive Character of the
  Seismic Coda}}, Bulletin of the Seismological Society of America, 90 (2000),
  pp.~655--665.

\bibitem{Snieder2002equilibration}
{\sc R.~Snieder}, {\em Coda wave interferometry and the equilibration of energy
  in elastic media}, Phys. Rev. E, 66 (2002), p.~046615.

\bibitem{soon2009hybridizable}
{\sc S.-C. Soon, B.~Cockburn, and H.~K. Stolarski}, {\em A hybridizable
  discontinuous {G}alerkin method for linear elasticity}, International journal
  for numerical methods in engineering, 80 (2009), pp.~1058--1092.

\bibitem{terrana2018spectral}
{\sc S.~Terrana, J.-P. Vilotte, and L.~Guillot}, {\em A spectral hybridizable
  discontinuous {G}alerkin method for elastic--acoustic wave propagation},
  Geophysical Journal International, 213 (2018), pp.~574--602.

\bibitem{thomsen1986weak}
{\sc L.~Thomsen}, {\em Weak elastic anisotropy}, Geophysics, 51 (1986),
  pp.~1954--1966.

\bibitem{tie2020systematic}
{\sc B.~Tie and A.-S. Mouronval}, {\em Systematic development of upwind
  numerical fluxes for the space discontinuous {G}alerkin method applied to
  elastic wave propagation in anisotropic and heterogeneous media with physical
  interfaces}, Computer Methods in Applied Mechanics and Engineering, 372
  (2020), p.~113352.

\bibitem{tie2018unified}
{\sc B.~Tie, A.-S. Mouronval, V.-D. Nguyen, L.~Series, and D.~Aubry}, {\em A
  unified variational framework for the space discontinuous {G}alerkin method
  for elastic wave propagation in anisotropic and piecewise homogeneous media},
  Computer Methods in Applied Mechanics and Engineering, 338 (2018),
  pp.~299--332.

\bibitem{tournier20223d}
{\sc P.-H. Tournier, P.~Jolivet, V.~Dolean, H.~S. Aghamiry, S.~Operto, and
  S.~Riffo}, {\em 3d finite-difference and finite-element frequency-domain wave
  simulation with multilevel optimized additive schwarz domain-decomposition
  preconditioner: A tool for full-waveform inversion of sparse node data sets},
  Geophysics, 87 (2022), pp.~T381--T402.

\bibitem{vila2021hybridisable}
{\sc J.~Vila-P{\'e}rez, M.~Giacomini, R.~Sevilla, and A.~Huerta}, {\em
  Hybridisable discontinuous {G}alerkin formulation of compressible flows},
  Archives of Computational Methods in Engineering, 28 (2021), pp.~753--784.

\bibitem{wang2018weak}
{\sc R.~Wang and R.~Zhang}, {\em A weak {G}alerkin finite element method for
  the linear elasticity problem in mixed form}, Journal of Computational
  Mathematics, 36 (2018), p.~469.

\bibitem{wilcox2010high}
{\sc L.~C. Wilcox, G.~Stadler, C.~Burstedde, and O.~Ghattas}, {\em A high-order
  discontinuous {G}alerkin method for wave propagation through coupled
  elastic--acoustic media}, Journal of Computational Physics, 229 (2010),
  pp.~9373--9396.

\bibitem{zhan2018exact}
{\sc Q.~Zhan, Q.~Ren, M.~Zhuang, Q.~Sun, and Q.~H. Liu}, {\em An exact
  {R}iemann solver for wave propagation in arbitrary anisotropic elastic media
  with fluid coupling}, Computer Methods in Applied Mechanics and Engineering,
  329 (2018), pp.~24--39.

\end{thebibliography}

\end{document}